\documentclass[a4paper,11pt]{article}
\usepackage[english]{babel}
\usepackage{amsfonts}
\usepackage{amssymb}
\usepackage{amsmath}
\usepackage{latexsym}
\usepackage{mathrsfs}
\usepackage{mathtools}
\usepackage{rotating}
\usepackage{pdflscape}
\usepackage[utf8]{inputenc} % Any characters can be typed directly from the keyboard
\usepackage{textcomp} % provide lots of new symbols
\usepackage{graphicx}  % Add graphics capabilities
\usepackage{epstopdf} % to include .eps graphics files with pdfLaTeX
\usepackage{bm}  % Define \bm{} to use bold math fonts
\usepackage{pdfsync}  
							
\begin{document}

\begin{center}
\Large\bf Analytic theory of finite 
asymptotic  expansions \\ in the real domain. Part II:\\ the factorizational theory for
Chebyshev asymptotic scales.
\end{center}

 \vspace{25pt}
\centerline{ANTONIO GRANATA}
 
\vspace{20pt}
 \centerline{\footnotesize\sl Dipartimento di Matematica e Informatica, Università della
Calabria,} \par\centerline{\footnotesize\sl 87036 Rende (Cosenza), Italy. email:
antonio.granata@unical.it}

 \vspace{50pt}
 {\footnotesize {\bf Abstract.} 
This paper contains a general theory for asymptotic expansions of type 
$$
 f(x)=a_1\phi_1(x)+\dots+a_n\phi_n(x)+o(\phi_n(x)),\ \  x\to x_0,\ n\ge3, \leqno(*)$$
 where the asymptotic scale 
$$
\phi_1(x)\gg\phi_2(x)\gg\cdots\gg\phi_n(x), \ \  x\to x_0, $$
 is assumed to be an extended complete Chebyshev system on a one-sided 
 neighborhood of $x_0$. ``Factorizational theory'' refers to proofs being  based  on  various types of factorizations of a  differential operator 
 associated to $(\phi_1,\cdots,\phi_n)$, and it is necessary to clearly understand the relationships between the nonvanishingness of the Wronskians $W(\phi_1,\dots,\phi_i)$, that of the Wronskians $W(\phi_n,\dots,\phi_{n-i})$ and the so-called canonical factorizations of disconjugate operators. 
 Next we focus on the guiding thread of our theory, which is
 the property of formal differentiation, aiming at characterizing
 some $n$-tuples of asymptotic expansions formed by $(*)$ and $n-1$ expansions 
 obtained by formal applications of suitable linear differential operators of orders 
 $1,2,\ldots,n-1$. Whereas for $n=2$ there are only two 
 such operators ``naturally'' suggested by the structure of the scale and the theory is comparatively simple, for $n\ge3$ 
 a result by Levin on the hierarchies of the Wronskians highlights a large class of  
 operators which preserve the hierarchy of the $\phi_i$'s and, as such, are a-priori 
 candidates to be formally applicable to (*). Our second preliminary step will be that of  noticing that the restricted class of the operators naturally 
 associated to ``canonical factorizations'' seems to be the most meaningful to be used in a  context of  formal differentiation. This gives rise to conjectures whose proofs build an analytic theory of finite asymptotic expansions in the real domain which, though not  elementary, parallels the familiar results about Taylor's formula.
  One of the results states that to each scale of the type under consideration it remains 
  associated an important class of functions (namely that of generalized convex 
  functions) enjoying the property that the expansion $(*)$, if valid, is automatically 
  formally differentiable $n-1$ times in two special senses.}
  
\vspace{30pt}

 \vspace{5pt} {\bf AMS subject classifications.}
41A60, 34E05, 26C99.

\newpage
\centerline{\bf  Contents}

1.\  Introduction

2.\  Canonical factorizations of  disconjugate operators and Chebyshev asymptotic scales

3.\  Applying differential operators to asymptotic scales

4.\  The first factorizational approach

5.\  The second factorizational approach and estimates of the remainder

 6.\  Absolute convergence and solutions of differential inequalities
 
 7.\  Proofs
  
  8.\ Appendix: algorithms for constructing canonical factorizations
  
\vspace{20pt}
\centerline{\bf 1. Introduction }
\vspace{10pt}

In this paper we develop a general analytic theory of asymptotic expansions of type
$$
 f(x)=a_1\phi_1(x)+\dots+
a_n\phi_n(x)+o\big(\phi_n(x)\big),\  x\to x_0; \  n\geq 3, \ \leqno(1.1)$$
where
$$
\phi_1(x)>>\phi_2(x)>>\dots >>\phi_n(x),\  x\to x_0.\ \leqno(1.2)$$

Though asymptotic expansions are since long a very useful tool in pure and applied mathematics, as far as asymptotic expansions in the {\it real domain} are concerned the general theory lacks basic results paralleling, for instance, (i) the classical Taylor's formula for polynomial expansions at a point $x_0\in\mathbb{R}$; (ii) the theory of polynomial expansions at $\infty$ systematized in [4]; (iii) the (not-too-trivial) case $n=2$ thoroughly investigated in [7]. Here we have in mind characterizations of (1.1) via integro-differential conditions useful for applications unlike the trivial characterization of (1.1) by means of the existence (as finite numbers) of the following $n$ limits defining the coefficients $a_i$:
 $$
 a_1 :=  \lim_{x\to x_0} f(x)/\phi_1(x) ,\ \ 
 a_i := \lim_{x\to x_0}\mbox{\Large$\frac{[f(x)-a_1\phi_1(x)-\ldots-a_{i-1}\phi_{i-1}(x)]}{\phi_i(x)}$},\ 2\leq i\leq n, \leqno (1.3)  $$
 the  $\phi_i$'s being supposed non-vanishing on a deleted neighborhood of $x_0$. 
 The three mentioned cases show that a proper approach to a satisfying theory consists in studying (1.1) not by itself but matched to other expansions obtained by formal application of certain differential operators. For this we need some preliminary material: a first part concerning factorizations of a linear ordinary differential operator and nonvanishingness of various Wronskians involving a basis of its kernel, a second part concerning those operators which act on the vector space ``$\textnormal{span}  \ (\phi_1,\dots ,\phi_n)$" preserving   asymptotic scales. 
The scale of comparison functions $(\phi_1,\dots ,\phi_n)$ is practically assumed to form an extended Chebyshev system on some left  deleted neighborhood of $x_0$ and two special types of factorizations, called canonical factorizations, are used: this is the content of  \S2 where various related properties are systematized about the concept of Chebyshev asymptotic scale.

 Now our theory revolves around the idea of  formal differentiation of an asymptotic expansion and, in general, applying an arbitrary  differential operator to an asymptotic expansion yields a meaningless result; so it is necessary to have some a-priori information on the differential operators which are most likely to be formally applicable to (1.1) in the sense of generating a new asymptotic expansion. A possible approach to obtain such an information consists in investigating the  case of an asymptotic expansion with an identically-zero remainder and this is done in \S3. For this case a deep result by Levin, already used in \S2, highlight certain differential operators, defined by means of Wronskians involving the $\phi_i$'s,  which preserve the hierarchy of the $\phi_i$'s and, as such, are a-priori candidates to be formally applicable to (1.1).  These operators may be too many (for $n\ge3$) to be included in a useful theory whereas canonical factorizations  automatically define two $(n-1)$-tuples of differential operators, of orders $1,2,\ldots,n-1$, which are practically more meaningful than a generic Levin's Wronskian. Their investigation gives rise to certain ``natural'' conjectures whose proofs are the core of our theory  called ``the factorizational theory'' and developed in \S\S4,5,6. All proofs are collected in \S7. The main features of this theory are:
	
	 (i)\ \  It yields applicable analytic characterizations of an expansion (1.1) matched to other asymptotic relations obtained by formal differentiations in suitable senses.
	 
	 (ii)\ \  For each Chebyshev asymptotic scale  there are at least two well-defined $(n-1)$-tuples of linear differential operators $(L_1,\dots,L_{n-1})$ and $(M_1,\dots,M_{n-1})$, of orders $1,2,\dots,n-1$ respectively, which can be formally applied to (1.1) under suitable integrability conditions. In one of the two circumstances useful representations of the remainders are also available.
		
	 (iii)\ \  A special family of functions is associated to each Chebyshev asymptotic scale namely that of generalized convex functions, for which the validity of the sole relation (1.1) automatically implies its formal differentiability $(n-1)$ times in the two senses involving the above-mentioned operators $(L_1,\dots,L_{n-1})$ and $(M_1,\dots,M_{n-1})$.
				
 The introductions in [4] and [7] contain other comments but the general theory to be developed in this paper is independent of  any previous results in these references: only the line of thought is the same. 
 
 In the appendix (\S8) we present two algorithms, admitting of asymptotic interpretations, for constructing canonical factorizations of disconjugate operators.
 
 Occasionally an asymptotic expansion 
$$
 f(x)=a_1\phi_1(x)+\dots+
a_i\phi_i(x)+o\big(\phi_i(x)\big),\  x\to x_0, \  i<n, \ \leqno(1.4)$$
 will be called ``incomplete" | with respect to the given scale $(\phi_1\dots,\phi_n)$, of course | whereas (1.1) will be called ``complete", and these locutions refer to the specified growth-order of the remainder and not to the terms effectively present in the expansion i.e. those with non-zero coefficients.

\vspace{10pt}
\centerline{\bf  Notations}

\vspace{10pt}
| $f\in AC^0(I)\equiv AC(I) \iff\ f$ is absolutely continuous on each compact subinterval of $I;\ \  f\in AC^{k}(I) \iff\ f^{(k)}\in AC(I)$;

| For $f\in AC^{k}(I)$ we write $\lim_{x\to x_0}f^{(k+1)}(x)$ meaning  
that $x$ runs through the points wherein $f^{(k+1)}$ exists as a finite number. Applying L'Hospital's rule in such a context means using Ostrowski's version [11] valid for absolutely continuous functions.

| $\overline{\mathbb{R}}:=\mathbb{R}\cup \{\pm\infty\}$ denotes the extended real line.

| If no ambiguity arises we use the following shorthand notations or similar ones:
$$\begin{cases}\displaystyle
 \int_T^x f_1 \int_T^{t_1}
f_2\dots\int_T^{t_{n-2}}f_{n-1}\int_T^{t_{n-1}}f_n(t) \, dt := \\ \\ \displaystyle
:=\int_T^x
f_1(t_1) dt_1 \int_T^{t_1}
f_2(t_2) dt_2\dots\int_T^{t_{n-2}}f_{n-1}(t_{n-1}) dt_{n-1}\int_T^{t_{n-1}}f_n(t_n) \, dt_n ; \end{cases}$$
$$\begin{cases}\displaystyle
 \int_x^{x_0}
f_1 \int_{t_1}^{x_0}
f_2\dots\int_{t_{n-2}}^{x_0}f_{n-1}\int_{t_{n-1}}^{x_0}f_n(t) \,dt :=\\ \\ \displaystyle
:=\int_x^{x_0}
f_1(t_1) dt_1 \int_{t_1}^{x_0}
f_2(t_2) dt_2\dots \int_{t_{n-2}}^{x_0}f_{n-1}(t_{n-1}) dt_{n-1}\int_{t_{n-1}}^{x_0}f_n(t_n) \,dt_n ; \end{cases}$$
wherein each integral $\int^{x_0}\equiv \int^{\to x_0}$ is to be understood as an improper integral.

| The acronyms we systematically use are: 

 T.A.S. := Chebyshev asymptotic scale: Def. 2.1; 
 
 C.F. := canonical factorization: Proposition 2.1-(iv) and (v).

| Propositions are numbered consecutively in each section irrespective of their labelling as lemma, theorem and so on.

 For later references we report here a classic fundamental formula for Wronskians:
$$
 W\big(\phi(x) \phi_1(x), \dots ,\phi(x)\phi_n(x)\big)=\big(\phi(x)\big)^n\cdot W\big(\phi_1(x), \dots ,\phi_n(x)\big),\leqno (1.5) $$
 valid under the required order of differentiability regardless of the sign of $\phi$.
 \vspace{20pt}
\begin{center}\bf 2.  Canonical factorizations of  disconjugate operators and \\ Chebyshev asymptotic scales
 \end{center}
 
 \vspace{10pt}
Our theory is built upon appropriate integral representations stemming from a special structure of the asymptotic scale $(\phi_1,\dots ,\phi_n)$: practically it forms a fundamental system of solutions of a disconjugate equation on a one-sided neighborhood of $x_0$ such that certain  Wronskians do not vanish thereon, a property granted by a result by Levin [9] which justifies our definition of Chebyshev asymptotic scale. We preliminarly recall some facts about factorizations of differential operators.

In this section $L_n,\ n\ge2$, denotes a linear ordinary differential operator of type
 $$
 L_nu:=u^{(n)}+\alpha_{n-1}(x)u^{(n-1)}+\ldots +\alpha_0(x)u \hspace{1pc}\forall \ u\in AC^{n-1}(J),\leqno (2.1)_1 $$
 $$
 \alpha_i\in L_{loc}^1(J),\  0 \leq i \leq n-1,\ J\ a\ generic\ interval\ of\ \mathbb{R},\leqno (2.1)_2 $$
  where $L_{loc}^1(J)$ denotes the class of functions Lebesgue-summable on every compact subinterval of $J$.  The matters to be discussed depend on the property of  disconjugacy and   several characterizations involving factorizations are collected in the next proposition where special locutions are defined in the statement itself.
  For general properties about disconjugacy we refer to the book by Coppel [1] and the paper by Levin [9], whereas for facts concerning canonical factorizations we refer to the papers by Trench [14] and the author [2; 3].
 
 \vspace{10pt}
 {\bf Proposition 2.1}\ (Disconjugacy on an open interval via factorizations). {\it For an operator $L_n$ of type $(2.1)_{1,2}$, $n\geq2$, on an open interval $]a,b[$, bounded or not, the following properties are equivalent:
 
  {\rm (i)}\ $L_n$ is disconjugate on $]a,b[$ in the sense that: every nontrivial solution of $L_nu=0$ has at most $(n-1)$ zeros on $]a,b[$ counting multiplicities or, equivalently, has at most $(n-1)$ distinct zeros on $]a,b[$.

   {\rm (ii)}\ $L_nu=0$ has a fundamental system of solutions on $]a,b[$, $(u_1, \dots ,u_n)$, satisfying P\'olya's W-property:
 $$
W\big(u_1(x),\dots,u_i(x)\big) >0 \hspace{1pc} \forall \ x\in  ]a,b[\ , \ \ 1\leq i\leq n; 
\leqno (2.2) $$
 or equivalently $L_nu=0$ has solutions $u_1, \dots ,u_{n-1}$ satisfying $(2.2)$ for $1\leq i\leq n-1$.

{\rm (iii)} \ $L_n$ has a P\'olya-Mammana factorization on $]a,b[$ i.e.
$$
L_nu\equiv r_n[r_{n-1}( \dots (r_1(r_0)')' \dots )']' \hspace{1pc} \forall \ u\in  AC^{n-1}]a,b[, \leqno(2.3) $$
where the $r_i$'s are suitable functions such that:
$$\begin{cases}
     r_i(x)>0 \ \ \forall \ x\in  ]a,b[; \ r_i\in  AC^{n-1-i}]a,b[\ , \ \ 0 \leq i\leq n-1; \cr
      r_n\in AC^0]a,b[. \end{cases} \leqno (2.4) $$ 

{\rm (iv)} \ $L_n$ has a ``canonical factorization $($C.F. for short $)$ of type {\rm (I)} at the endpoint a", i.e. a factorization of type {\rm (2.3)-(2.4)} with the additional conditions:
$$
\int_{\to a}(1/r_i) =\,+ \infty,\;1\leq i\leq n-1, \leqno (2.5)_a $$
and a similar ``C.F. of type {\rm (I)} at the endpoint b", i.e. with the $r_i$'s satisfying
$$
\int^{\to b}(1/r_i)=\,+\infty,\;1\leq i\leq n-1. \leqno (2.5)_b $$

{\rm (v)} \ For each $c,\, a<c<b,\, L_n $ has a "C.F. on the interval $]a,c[$ which is of type {\rm (II)} at the endpoint $a$", i.e. a factorization {\rm (2.3)-(2.4)} valid on the interval $]a,c[$ and with the $r_i$'s satisfying
$$
\int_{\to a}(1/r_i) <\,+ \infty,\;1\leq i\leq n-1. \leqno (2.6)_a$$

And $L_n$ has a ``C.F. on the interval $]c,b[$ which is of type {\rm (II)} at the endpoint $b$", i.e. a factorization {\rm (2.3)-(2.4)} valid on the interval $]c,b[$ and with the $r_i$'s satisfying}
$$
\int^{\to b}(1/r_i) <\,+ \infty,\;1\leq i\leq n-1.   \leqno (2.6)_b $$

\vspace{10pt}
{\it Remarks.}  
1. In the definition of a C.F. conditions (2.5) or (2.6) are required to hold for the index $i$ running from 1 to $(n-1)$: there are no conditions on $r_0$ and $r_n$.
 Factorizations in properties (iii)-(iv) are global i.e. valid on the whole given interval $]a,b[$, whereas property (v) claims the existence of local C.F.'s of type (II). The existence of a global  C.F. of type (II)  at $a$ or at $b$ is a special circumstance [2; Thm.  3.11, p. 163].

2. A global C.F. of type (I) at a specified endpoint does always exist for a disconjugate operator on an open interval and is ``essentially" unique in the sense that the functions $r_i$ are determined up to multiplicative constants with product 1: Trench [14].
The situation is quite different for C.F.'s of type (II). For example the operator $L_n \equiv u^{(n)}$ has no global C.F. on $(-\infty,+\infty)$ of type (II) at any of the endpoints for it has only ``one" (up to constant factors) P\'olya-Mammana factorization on $(-\infty,+\infty)$, namely
$$
u^{(n)}\equiv ( \dots (u')' \dots )' \,, \leqno(2.7) $$
which is a special contingency characterized in [2; Thm.  3.3] and in [3; Thm.  7.1]. But the operator $u^{(n)} $ thought of as acting on the space $AC^{n-1}]0,+\infty)$, or even on the space $C^\infty]0,+\infty)$, has infinitely many ``essentially" different C.F.'s of type (II), for instance the following ones
$$
u^{(n)}\equiv \frac{1}{(x-c)^{n-1}}\left[(x-c)^2\left(\dots\left((x-c)^2\left(\frac{u}{(x-c)^{n-1}}\right)'\right)'\dots\right)'\right]' \ , \leqno(2.8) $$
which are C.F.'s of type (II) at both the endpoints ``$0$" and ``$+\infty$" whatever the choice of the constant $c<0$. For $c=0$ we get a factorization on $]0,+\infty)$ which is a C.F.  of type (I) at ``$0$" and of type (II) at  ``$+\infty$"; for  $c>0$ we have nonglobal factorizations which are of type (II) at $+\infty$.

\vspace{5pt}
C.F.'s are naturally linked to bases of ker $L_n$ forming asymptotic scales at one or both endpoints and the following results, due to Levin [9; \S2], highlight important properties of the Wronskians constructed with an  asymptotic scale.

\vspace{10pt}
{\bf Proposition 2.2}\ (Wronskians of asymptotic scales and their hierarchies). 

(I)\ (Results involving a differential operator).\ {\it Let $L_n$ be an operator of type $(2.1)_{1,2}$ disconjugate on an open interval $]a,b[$. Then:

{\rm(i)}\ Its kernel has some basis $(\phi_1,\dots,\phi_n)$ satisfying\,$:$
$$\begin{cases}
\phi_i(x)>0 \ \textit{on some interval} \ ]b-\epsilon,b[\ ,\ 1\leq i \leq n; \cr
\phi_1(x) \gg \phi_2(x) \gg \dots \gg \phi_n(x)  , \ x \to b^-. 
\end{cases} \leqno (2.9) $$

{\rm(ii)} \ For each such basis:
$$ 
W\big(\phi_n(x),\phi_{n-1}(x),\dots,\phi_i(x)\big)>0 \ \textit{on the whole interval} \ \  ]a,b[\ , \ 1 \leq i \leq n, \leqno (2.10)$$
noticing the reversed order of the $\phi_i$'s in the Wronskians.

{\rm(iii)} \ For any strictly decreasing set of indexes $\{ i_1, \dots ,i_k \} $, i.e. such that
$$
n \geq i_1>i_2> \dots >i_k \geq 1, \ 1 \leq k \leq n-1, \leqno (2.11) $$
 we have:
$$
W\big(\phi_{i_{1}}(x), \dots ,\phi_{i_{k}}(x)\big)> 0 \ \textit{on a left deleted neighborhood of} \ b, \leqno (2.12) $$ 
and in particular we have the inequalities:
$$\begin{cases}
sign\ W\big(\phi_1(x), \dots ,\phi_i(x)\big)=(-1)^{i(i-1)/2}
 \\\mbox{on a left deleted neighborhood of} \ b,\ 1\le i\le n.\end{cases} \leqno (2.13) $$

{\rm(iv)}  \ For each $k$, $1 \leq k \leq n-1$, and for any two distinct and strictly increasing sets of indexes ${ i_1, \dots ,i_k}$  and ${ j_1, \dots ,j_k}$ such that $i_h \leq j_h,\ 1\leq h\leq k$,  we have
$$
 W\big(\phi_{i_{1}}(x), \dots ,\phi_{i_{k}}(x)\big) \gg W\big(\phi_{j_{1}}(x), \dots ,\phi_{j_{k}}(x)\big), \ x \to b^- . \leqno (2.14) $$
 
 Notice the ordering of the $\phi_i$'s and the $\phi_j$'s in {\rm (2.14)}: if each $\phi_i$ has an order of growth at $b^-$ greater than that of the corresponding $\phi_j$ then the same is true for the Wronskians. In the claim {\rm(iii)} we have a different ordering of the $\phi_i$'s as this grants the positivity of the Wronskians in {\rm (2.12)}.}
 
(II)\ (Results involving scales with less regularity).\ {\it Let $(\phi_1,\dots,\phi_n)$ be functions of class $C^{n-1}]a,b[$ satisfying conditions $(2.9)$ and condition 
$$ 
W\big(\phi_n(x),\phi_{n-1}(x),\dots,\phi_1(x)\big)\ \ either \ \ge0\ \ or \ \le0\  on \  ]a,b[\ , \ 1 \leq i \leq n; \leqno (2.15)$$
and let there exist an integer $r,\ 1\le r\le n,$ such that:
$$
W_r\big(\phi_n(x),\phi_{n-1}(x),\dots,\phi_i(x)\big)\ne0\ on \  ]a,b[\ , \ \forall\ i,\ 1 \leq i \leq n; \leqno (2.16)$$
where the symbol $\ W_r\big(\phi_n(x),\phi_{n-1}(x),\dots,\phi_i(x)\big)$ denotes the  Wronskian determinant wherein the column involving $\phi_r$ has been suppressed.
Then they hold true:
$$
W\big(\phi_n(x),\phi_{n-1}(x),\dots,\phi_1(x)\big) \ge0\  on \  ]a,b[\,; \leqno (2.17)$$
$$
W\big(\phi_n(x),\phi_{n-1}(x),\dots,\phi_i(x)\big)>0\ on \  ]a,b[\ , \ 1 \leq i \leq n-1; \leqno (2.18)$$
and the above-stated properties in {\rm(iii)} and {\rm(iv)}. Notice that in {\rm (2.17)-(2.18)} the signs of the Wronskians are well defined even if they remain undefined in the assumptions {\rm (2.15)-(2.16)}.}

 \vspace{5pt}
To visualize (2.14) we list a few asymptotic scales at $b^-$  constructed with the Wronskians:
 $$ \begin{cases}
 W(\phi_1, \phi_2)\gg W(\phi_1, \phi_3)\gg \dots \gg W(\phi_1, \phi_n) \cr
 W(\phi_2, \phi_3) \gg W(\phi_2, \phi_4)\gg \dots \gg W(\phi_2, \phi_n) \cr
 \dots \dots \dots \dots \dots \dots \dots\dots \dots \dots \dots \dots \dots \dots 
 \quad , \  \ x \to b^-,     \cr
 W(\phi_{n-2}, \phi_{n-1}) \gg W(\phi_{n-2}, \phi_n) \end{cases} \leqno (2.19)  $$
 $$
 W(\phi_1, \phi_2, \phi_3) \gg W(\phi_1, \phi_2, \phi_4) \gg \dots \gg W(\phi_1, \phi_2, \phi_n),
 \ x \to b^-.  \leqno (2.20) $$
  
\vspace{10pt}
It is quite important to note the order of the $\phi_i$'s forming the asymptotic scale in (2.9); if we mantain the same ordering in the analogous statement for $x\to a^+$, i.e.
$\phi_1(x) \gg \phi_2(x) \gg \dots \gg \phi_n(x)  , \ x \to a^+$, then the Wronskians in (2.10) and in  (2.12) to (2.18) are the same, the essential point being the relative growth-order of the $\phi_i$'s. From the point of view of asymptotic expansions the correct numbering is that adopted by us irrespective of the limiting process. 

The above results substantiate the following definition of special asymptotic scales wherein we merely fix the neighborhood of $b$ left undefined in Proposition 2.2 whose part (I) grants the existence of such scales whereas part (II) implies a lot of useful properties even for scales with less regularity. From now on the interval will be denoted as in the two-term theory [7].

\vspace{5pt}
\textbf{Definition 2.1}\ (Chebyshev asymptotic scales).\ {\it  The ordered n-tuple of real-valued functions $(\phi_1, \dots ,\phi_n), n\ge2,$ is termed a ``Chebyshev asymptotic scale'' {\rm(T.A.S. for short)} on the half-open interval $[T,x_0[\, ,\ T\in\mathbb{R},\  x_0\leq+\infty$, provided the following properties are satisfied:
$$
\phi_i\in C^{n-1}[T,x_0[\, ,\ 1\leq i\leq n;\ \leqno(2.21) $$
$$
\phi_i(x)\neq0 \ \ \mbox{on some left deleted neighborhood of}\ x_0, \ 1\leq i\leq n;\leqno(2.22)$$
$$
\phi_1(x)\gg\phi_2(x)\gg\cdots\gg\phi_n(x),\  x\to x_0^-; \leqno(2.23) $$
$$
W\big(\phi_1(x),\dots,\phi_i(x)\big) \ne0 \ on\ [T,x_0[\, , \ 1\leq i\leq n.\leqno(2.24)$$

 Whenever the $\phi_i$'s satisfy the stronger regularity condition
 $$
 \phi_i\in AC^{n-1}[T,x_0[\, ,\ 1\leq i\leq n,\leqno(2.25)$$
 they remain associated to  the operator:
 $$
  L_{\phi_1,\dots,\phi_n}u:=W\big(\phi_1(x),\dots,\phi_n(x),u\big)\left/ W\big(\phi_1(x),\dots,\phi_n(x)\big)\right.,\leqno(2.26)$$
which is the unique linear ordinary differential operator of type $(2.1)_{1,2}$, acting on the space $AC^{n-1}[T,x_0[$ and such that}
 $\textnormal {ker} \ L_{\phi_1,\dots,\phi_n}=\textnormal{span}  \, (\phi_1,\dots,\phi_n).$
 
\vspace{5pt}
{\it Remarks.}\ 1. Condition (2.21) is the usual regularity assumption in approximation theory (Chebyshev systems and the like), whereas in matters involving differential equations/inequalities it is natural to assume (2.25).

2. Choosing an half-open interval in this definition is a matter of convenience: the point $x_0$ involved in the asymptotic relations is characterized as the endpoint not belonging to the interval, possibly  $x_0=+\infty$, whereas the other endpoint marks off an interval whereon the inequalities involving the Wronskians are satisfied and these in turn allow certain integral representations valid on the whole given interval and essential to our theory. These remarks make evident the analogous definition for an interval $]x_0,T]$ where: $-\infty\le x_0$, and $T\in\mathbb{R}$. 

3. In the above definition we have merely supposed the nonvanishingness of various functions instead of specifying their signs as in Proposition 2.2; this avoids restrictions that are immaterial in asymptotic investigations. If the $\phi_i$'s are strictly positive near $x_0$ then Levin's theorem provides the exact signs of certain Wronskians.

4. As concrete examples of such asymptotic scales on $[T,+\infty )$ the reader may think of scales  whose non-identically zero and infinitely-differentiable functions are represented by linear combinations, products, ratios and compositions of a finite number of powers, exponentials and logarithms. As a rule such functions and their Wronskians have a principal part at $+\infty$ which can be expressed by products of similar functions, hence they do not vanish on a neighborhood of $+\infty$.

\vskip5pt
When comparing our notations with other authors' results the reader must carefully notice the numbering of the $\phi_i$'s in the asymptotic scale (2.23) and in the Wronskians (2.24); the next proposition contains various additional properties of a T.A.S. and, in particular, it claims that conditions (2.21)-(2.24) imply the nonvanishingness of the reversed Wronskians: 
$$
W\big(\phi_n(x),\phi_{n-1}(x),\dots,\phi_i(x)\big) \ne0 \ on\ [T,x_0[\, , 
\ 1,\leq i\leq n,\leqno(2.27)$$
though the converse generally fails as it may be easily checked  for the scale:
    $$1\gg cx+x^2\gg x^2,\ x\to 0^-,\ (c >0),\ on\ (-\infty,0[\ . \leqno (2.28) $$
    
  With our notations this scale satisfies (2.27) on $(-\infty,0[$ whereas:
  $$\begin{cases}
  \phi_1\ and\  W(\phi_1,\phi_2,\phi_3)\ne0 \ on\ (-\infty,0[\ ; \\
  W(\phi_1,\phi_2)\equiv W(1,cx+x^2)=c+2x\ne0\ on\ ]-c/2,0[\ but\ not\ on\ (-\infty,0[ \ .\end{cases}$$

 \vspace{10pt}
 {\bf Proposition 2.3}\ (Several characterizations and additional properties of T.A.S.'s). {\it Let the ordered n-tuple of real-valued functions $(\phi_1, \dots ,\phi_n), n\ge2,$ satisfy conditions {\rm (2.21)-(2.22)-(2.23)}.
 
  {\rm(I)}\ The following are equivalent properties:

{\rm (i)}\ $(\phi_1, \dots ,\phi_n)$ is a T.A.S. on $[T,x_0[$, i.e. $(2.24)$ hold true.

{\rm (ii)}\ Both sets of inequalities $(2.24)$ and $(2.27)$ hold true.

{\rm (iii)}\ The ordered n-tuple  $(\epsilon_1\phi_1(x), \dots ,\epsilon_n\phi_n(x))$, with proper choices of the constants $\epsilon_i=\pm1$, is  an extended complete Chebyshev system on $[T,x_0[$\ .

{\rm (iv)}\  The n-tuple $(\phi_1, \dots ,\phi_n)$ admits of an integral representation of the form
$$\begin{cases}\displaystyle
     \phi_1(x)=w_0(x); \ \ \phi_2(x)=w_0(x)\cdot \int_x^{x_0} w_1; 
     \\ \\ \displaystyle
   \phi_i(x)= w_0(x)\cdot \int_x^{x_0} w_1 \dots \int _{t_{i-2}}^{x_0} w_{i-1}, \ 2 \leq i \le n,\ x\in [T,x_0[, \end{cases} \leqno (2.29) $$
    with suitable functions  $w_i$  subjected to the following regularity conditions:
 $$\begin{cases}
 w_i(x)\ne0\ \ \forall \ x \in [T,x_0[;\  w_i\in C^{n-1-i}[T,x_0[\ ,\ \ 0 \leq i \leq n-1; \\
 \int^{x_0}|w_i|<+\infty, \ 1\leq i \leq n-1.\end{cases}\leqno (2.30)$$
 
 If this is the case the $w_i$'s are unique and may be expressed in terms of the $\phi_i$'s  on $[T,x_0[$ by the formulas:
  $$\begin{cases}
  w_0(x):=\phi_1(x);\  w_1 := -\big(\phi_2(x)/ \phi_1(x)\big)'=-\big(\phi_1(x)\big)^{-2}W\big(\phi_1(x),\phi_2(x)\big); 
  \\ \\ \displaystyle
    w_i(x):=-\left[\ \frac{W( \phi_1(x), \dots ,\phi_{i-1}(x),\phi_{i+1}(x))}{W( \phi_1(x), \dots ,\phi_{i-1}(x),\phi_i(x))}\  \right] ' \equiv \\ \\
    \equiv  -\big[W( \phi_1, \dots ,\phi_{i-1}) \cdot W( \phi_1, \dots ,\phi_{i+1})\big]\big/ \big[ W( \phi_1, \dots ,\phi_i) \big]^2  , \ 2 \le i \le n-1.\end{cases} \leqno (2.31) $$ 
     
  Conversely we have the following formulas for the Wronskians of the $\phi_i$'s:
  $$
     W\big(\phi_1,\dots,\phi_i\big)=(-1)^{i(i-1)/2}w_0^iw_1^{i-1} w_2^{i-2}\ldots w_{i-1} \ \ on\ [T,x_0[\, , \ 2\leq i\leq n.\leqno(2.32)$$
     
 {\rm(II)}\ For $(\phi_1, \dots ,\phi_n)$ a T.A.S. on $[T,x_0[$ we have the inequalities:
 $$
 \phi_i(x)\ne0\ on\ [T,x_0[\, ,\ 1\le i \le n,\ {\rm (implied\ by (2.29)-(2.30) );} \leqno (2.33)$$
  $$
W\big(\phi_{i_{1}}(x), \dots ,\phi_{i_{k}}(x)\big)\ne0 \ near\ x_0, \leqno (2.34) $$
for any set of indexes satisfying {\rm(2.11)} and we also have the hierarchies between the Wronskians stated in Proposition $2.2$-{\rm (iv)}  and referred to $x\to x_0^-$. Whenever the $\phi_i$'s are strictly positive then all the Wronskians in $(2.27)$ are strictly positive on $[T,x_0[$ by $(2.10)$, but not necessarily all the Wronskians in $(2.24)$; in this case the inverted $n$-tuple $(\phi_n, \dots ,\phi_1)$ is an extended complete Chebyshev system on $[T,x_0[$. On the contrary, if the given $n$-tuple $(\phi_1, \dots ,\phi_n)$ is an extended complete Chebyshev system on $[T,x_0[$, i.e. all the Wronskians in $(2.24)$ are  strictly positive on $[T,x_0[$, then $(2.29)$ and $(2.31)$ imply that the $\phi_i$'s have alternating signs, namely}:  $
  {\rm sign}\ \phi_i=(-1)^{i-1}\ on\ [T,x_0[$.

\vskip10pt
Part (I) of Proposition 2.3 generalizes a classical result, [8; Ch. XI, Th. 1.2, p. 379], which characterizes those special asymptotic scales formed by functions with zeros of increasing multiplicities   (namely $0,1,\ldots,n-1$) at an endopint of a compact interval; also refer to [9; Ch. 1] for locutions about Chebyshev systems.  Notice that formulas (2.31) in themselves are well defined if  the n-tuple $(\phi_1, \dots ,\phi_n)$ satisfies (2.21) and (2.24); under the addditional  assumption (2.23)  they establish a one-to-one correspondence between the $\phi_i$'s and the $w_i$'s. Formulas (2.31) for $i\ge3$ are not obvious consequences of (2.29): see the few introductory lines at the outset of \S8. For a T.A.S. on $]x_0,T]$ the integrals $\int_x^{x_0}$ in (2.29) are obviously replaced by $\int^x_{x_0}$, the $w_i$'s in (2.31) for $i\ge1$ are defined without the minus sign  and the coefficient $(-1)^{\ldots}$ is absent in (2.32). If all the Wronskians in $(2.24)$ are  strictly positive on $]x_0,T]$ then the same is true for all the $\phi_i$'s. 

\vskip5pt
Under condition (2.25)  formulas in Proposition 2.3-(ii) are related to C.F.'s of type (II) at $x_0$ and certain calculations used in our proof give quick proofs for the existence of both types of C.F.'s. We collect in the next proposition all the facts essential to  develop our theory of asymptotic expansions; here the focus is on C.F's rather than on integral representations of the given scale because we need both types of C.F.'s and the layout of Proposition 2.3 does not suit a C.F. of type (I) . 
 
 \vspace{10pt}
 {\bf Proposition 2.4}\ (Formulas concerning T.A.S.'s linked to differential operators). \textit{Let the ordered $n$-tuple $(\phi_1,\dots,\phi_n)$ satisfy conditions $(2.21)$ to $(2.25)$, hence the operator in $(2.26)$ is disconjugate on the open interval $]T,x_0[$ and enjoys the properties in Propositions {\rm2.1} and \ {\rm 2.2-(I)}. Moreover, as an operator acting on $AC^{n-1}[T,x_0[$ it has the following further properties:}
 
 (i)\ \ \textit{Define the following $(n+1)$ functions on $[T,x_0[$:}
  $$
   \begin{cases}
    q_0 := 1/ |\phi _1| ; \  \  q_1 := ( \phi_1)^2 \big/\ |W(\phi_1, \phi_2)|; \cr  
   q_i:= \big[ W( \phi_1, \dots ,\phi_i) \big]^2 \big/\ \big|W( \phi_1, \dots ,\phi_{i-1}) \cdot W( \phi_1, \dots ,\phi_{i+1})\big|, \ 2 \le i \le n-1; \cr 
   q_n := |q_0 q_1 \dots q_{n-1}|^{-1} \equiv \big|W( \phi_1, \dots ,\phi_n) / W( \phi_1, \dots ,\phi_{n-1})\big|.\end{cases}  \leqno (2.35) $$ 
   
   \textit{Then the $q_i$'s satisfy the following regularity conditions:}
   $$
\begin{cases}
 q_i(x)> 0\ \ \forall \ x \in [T,x_0[;\ q_i\in AC^{n-1-i}[T,x_0[\ ,\ \ 0 \leq i \leq n-1;  \\ q_n\in AC^0[T,x_0[. \end{cases} \leqno (2.36) $$

   \textit{Their reciprocals, left apart $q_0$ and $q_n$, may be expressed as derivatives of certain ratios}
    $$ \begin{cases}
    1/q_1(x)=\big|(\phi_2(x)/ \phi_1(x))'\big|, \\ \\ \displaystyle
    1/q_i(x)= \left|\left[\ \frac{W( \phi_1(x), \dots ,\phi_{i-1}(x),\phi_{i+1}(x))}{W( \phi_1(x), \dots ,\phi_{i-1}(x),\phi_i(x))}\  \right] ' \right| , \ 2 \le i \le n-1, \end{cases}  \leqno (2.37) $$ 
    \textit{on the interval $[T,x_0[$, and}
    $$
    \int _T^{x_0}(1/q_i)<+ \infty, \ 1 \leq i \leq n-1. \leqno (2.38) $$
    
    \textit{Our operator admits of the following factorization on $[T,x_0[$:}
   $$
        L_{\phi_1,\dots,\phi_n}u \equiv  q_n\big[q_{n-1}(\dots(q_0u)'\dots)'\big]' , 
        \leqno (2.39)  $$
   \textit{which is a global C.F. of type {\rm(II)} at both endpoints $T$ and $x_0$.}
 
 (ii) \ \textit{Our T.A.S. (apart from the signs) admits of the following integral representation in terms of the $q_i$'s:}
    $$\begin{cases}\displaystyle
     |\phi_1(x)|=\frac{1}{q_0(x)}; \ \ |\phi_2(x)|= \frac{1}{q_0(x)} \int_x^{x_0} \frac{1}{q_1}; \\ \\ \displaystyle
   |\phi_i(x)|= \frac{1}{q_0(x)} \int_x^{x_0} \frac{1}{q_1} \dots \int _{t_{i-2}}^{x_0} \frac{1}{q_{i-1}}, \ 2 \leq i \le n,\ x\in [T,x_0[; \end{cases} \leqno (2.40) $$
  \textit{hence the $\phi_i$'s, besides being everywhere non-zero on $[T,x_0[$, have the same order of growth at $T$, namely}
   $$
   \lim_{x \to T^+} \phi_i(x)/ \phi_j(x)=c_{ij} \in \mathbb{R} \backslash \{0\} \ \forall \  i \neq j.  \leqno (2.41) $$
  
  \textit {In the special case where all the Wronskians in $(2.24)$ are strictly positive, i.e. when   $(\phi_1,\dots,\phi_n)$ is an extended complete Chebyshev system on $[T,x_0[$, then the $\phi_i$'s have alternating signs, namely}
  $$
  {\rm sign}\ \phi_i=(-1)^{i-1}\ on\ [T,x_0[.\leqno (2.42) $$  
  
  (iii) \ \textit {Analogously we define the following $(n+1)$ functions on $[T,x_0[$:}
  $$ \begin{cases}
    p_0 := 1/|\phi _n| ; \  \  p_1 := ( \phi_n)^2 \big/\ |W(\phi_n,\phi_{n-1})|; \\  \\
   p_i:= \big[W( \phi_n,\phi_{n-1}, \dots ,\phi_{n-i+1}) \big]^2 \times \\
   \times \big| W(\phi_n,\phi_{n-1}, \dots ,\phi_{n-i+2}) \cdot W( \phi_n,\phi_{n-1}, \dots ,\phi_{n-i})\big|^{-1},\ \ 2 \le i \le n-1;  \\ \\
   p_n :=\big| p_0 p_1 \dots p_{n-1}\big|^{-1} \equiv  \big|W(\phi_n,\phi_{n-1}, \dots ,\phi_1) /W(\phi_n,\phi_{n-1}, \dots ,\phi_2)\big|\ .\end{cases}  \leqno (2.43) $$ 
   
  \textit {They satisfy the same regularity conditions on the half-open interval $[T,x_0[$ as the $q_i$'s do in $(2.36)$ and their reciprocals may be expressed as derivatives of the following ratios analogous to those in $(2.37)$:}
   $$ \begin{cases}
    1/p_1(x)=\big|(\phi_{n-1}(x)/ \phi_n(x))'\big|,  \\ \\ \displaystyle
    1/p_i(x)= \left| \left[\  \frac{W( \phi_n(x), \dots ,\phi_{n-i+2}(x),\phi_{n-i}(x))}{W( \phi_n(x), \dots ,\phi_{n-i+2}(x),\phi_{n-i+1}(x))} \ \right] '\right|  , \ \ 2 \le i \le n-1; \end{cases}  \leqno (2.44) $$
    \textit{on the interval $[T,x_0[$. Moreover:}
 $$
 \int ^{x_0}(1/p_i)=+ \infty, \ 1 \leq i \leq n-1, \leqno (2.45) $$
 $$
  \int _T(1/p_i)<+ \infty, \ 1 \leq i \leq n-1, \leqno (2.46) $$
 \textit{hence the associated factorization      
 $$
L_{\phi_1,\dots ,\phi_n}u \equiv  p_n\big[p_{n-1}(\dots(p_0u)' \dots )'\big]' ,
 \leqno (2.47) $$
    is (constant factors apart) ``the" global C.F. of $L_{\phi_1,\dots,\phi_n}$ of type {\rm (I)} at $x_0$ and it turns out to be of type  {\rm (II)} at $T$.}
     
  (iv) \ \textit{The special fundamental system  of solutions to 
    $L_{\phi_1,\dots,\phi_n}u=0$ defined by}
    $$ \begin{cases} \displaystyle
    P_0(x):= \frac{1}{p_0(x)}; \ \ P_1(x):= \frac{1}{p_0(x)} \int_T^x \frac{1}{p_1} ; \\ \\\displaystyle
    P_i(x):= \frac{1}{p_0(x)}  \int_T^x \frac{1}{p_1} \dots \int_T^{t_{i-1}} \frac{1}{p_i} ,\ 1 \leq i \leq n-1,\end{cases}  \leqno (2.48) $$
    \textit{satisfies the asymptotic relations:}
    $$ \begin{cases}
    P_0(x) \gg P_1(x) \gg \dots \gg P_{n-2}(x) \gg P_{n-1}(x), & x \to T^+, \cr 
   P_{n-1}(x) \gg P_{n-2}(x) \gg \dots \gg P_1(x)\gg P_0(x) , & x \to x_0^-. \end{cases} \leqno (2.49) $$
     
    \textit{Relations $(2.49)$ uniquely determine  the fundamental system $(P_0, \dots, P_{n-1})$ up to multiplicative constants.}  (In the terminology used by the author [2; 3] the $n$-tuple $(P_0, \dots, P_{n-1})$ is a "mixed hierarchical system" on $]T,x_0[$ whereas Levin [9; p. 80] would call it a "doubly hierarchical system" because he uses different arrangements for asymptotic scales at the left or right endpoints [9; p. 59].)
 {\it Whenever the $\phi_i$'s are strictly positive then the same is true for all the Wronskians appearing in $(2.43)$ hence the absolute values are redundant; in this case  it is the inverted $n$-tuple $(\phi_n,\dots,\phi_1)$ which forms  an extended complete Chebyshev system on $[T,x_0[$.}
 
  \vspace{10pt}
The construction of the two above factorizations starting from the given expressions of the coefficients $q_i$ or $p_i$ is the classical procedure by P\'olya [12]. Notice that the functions $p_i$'s in (2.47), which are unique (constant factors apart) by a mentioned result by Trench, may be recovered from many different asymptotic scales and not just from one! The main feature of the above proposition is that we can express all the properties of our basic operator (at least those needed in our theory) in terms of the a-priori given Chebyshev asymptotic scale.  The use of absolute values in the definitions of the $q_i$'s and $p_i$'s , though causing some incoveniences in the sequel, has the advantage of avoiding their use in the everywhere-present integral representations; and we must use them in at least one of the definitions as the two sets of Wronskians cannot have one and the same sign.
  
    \vspace{5pt}
{\it A quick proof  of the existence of C.F.'s } .\ \ The global existence of C.F.'s of type (I) was for the first time  proved by Trench [14]  by an original procedure which was subsequently adapted by the author [2] to show  the local existence of C.F.'s of type (II). Trench's result played a historical role as it had a great impact on the asymptotic theory of ordinary differential equations.  Levin's theorem easily implies Trench's result about global existence (but not uniqueness) in the case of disconjugate operators and the existence of a particular local C.F.'s of type (II): see the proof of Proposition 2.4. However we must point out that Trench's procedure, independent of properties of Wroskians, applies to a larger class of operators [14, \S1].
     As far as C.F.'s of type (II) are concerned the present quick approach does not yield a C.F. of type (II) at $b$ for each interval $]a+\epsilon,b[$, as asserted in Proposition 2.1-(v).
     
\vspace{20pt}  
   \centerline{\bf 3.  Applying differential operators to asymptotic scales}
   
   \vspace{10pt}
     In the elementary case of Taylor's formula the simple condition
     $$\exists\ f^{(n)}(x_0)\leqno(3.1)$$
     is not a mere sufficient condition for the validity of the asymptotic expansion
     $$
     f(x)=\sum_{i=0}^na_i(x-x_0)^i+o\big((x-x_0)^n\big)\equiv 
     T_n(x)+o\big((x-x_0)^n\big),\ x\to x_0;\leqno(3.2)$$
     it in fact characterizes the set of the $n$ asymptotic expansions
     $$
     f^{(k)}(x)=\sum_{i=0}^{n-k}T_n^{(k)}(x)+o\big((x-x_0)^{n-k}\big),\ x\to x_0,
     \ \ 0\le k\le n-1,\leqno(3.3)$$
     which is formed by (3.2) together with the relations obtained by formal differentiation $1,2,\ldots n-1$ times. In this case we have the known formulas for the coeffficients:
     $$a_i=f^{(i)}(x_0)/i!\ ,\ \  0\le i\le n.\leqno(3.4)$$
     
     If we strenghten condition (3.1) by assuming
$$f\in AC^n\big(I_{x_0}\big),\ I_{x_0}:a\ neighborhood\ of\ x_0,\leqno(3.5)$$
we have representation
$$f^{(n)}(x)=f^{(n)}(x_o)+\int_{x_0}^xf^{(n+1)}(t)dt,\leqno(3.6)$$
which, besides implying the validity of (3.3) for $k=n$ as well, gives rise to the integral representation formulas of all the remainders in (3.3).
     
          A similar situation occurs in the factorizational theory of polynomial asymptotic expansions at $+\infty$, [4], where the standard operators of differentiation 
$D^k:=d^k/dx^k$ happen to be formally applicable $n$ times to the expansion
$$
f(x)=a_nx^n+\ldots+a_1x+a_0+o(1),\ x\to+\infty,\leqno(3.7)$$
in two quite different senses and under suitable integral conditions. But in the analogous theory for expansions in arbitrary real powers
$$
f(x)=a_1x^{\alpha_1}+\ldots+a_nx^{\alpha_n}+o\big(x^{\alpha_n}\big),\ x\to+\infty,
\ (\alpha_1>\ldots>\alpha_n),\leqno(3.8)$$
developed in [6], it turns out that the most natural operators on which to build a satisfying theory are those linked to the C.F.'s of the differential operator in (2.26) with $\phi_i(x):=x^{\alpha_i}$ and not the operators $D^k$ though, in this special instance, the set of the formally-differentiated expansions may be equivalently expressed by expansions involving the standad derivatives. In the present general context  it is good  to preliminarly investigate which differential operators are likely to be formally applicable to an expansion (1.1) and a possible approach consists in investigating the case of an asymptotic expansion with a zero remainder i.e. a relation of type
     $$
      f(x)=a_1\phi_1(x)+\dots+a_n\phi_n(x). \leqno(3.9)$$

A first very general answer comes from the hierarchies of the Wronskians; a second, less general but practically more meaningful, answer comes from the use of C.F.'s; third, a C.F. of type (II) turns out to play a special role in computing the coefficients of an asymptotic expansion. We have used locutions such as ``formal application of an operator'' in an intuitive way but we give here a precise definition to avoid possible incongruences arising from identically-zero terms.

\vspace{10pt }
 {\bf Definition 3.1}\ (Asymptotically-admissible operators).\ \textit{Let $\mathcal{L}$ be a linear operator acting between two linear spaces of real- or complex-valued functions of one real variable. If  $\mathcal{L}[\phi_i(x)]\equiv 0$ on some neighborhood of $x_0 \ \forall \ i$ then the concept in question is not defined, otherwise we put:}
$$
m:=\max\{i\in \{1,\dots ,n\}: \mathcal{L}[\phi_i(x)]\not \equiv 0 \  \textit{on a neighborhood of} \ x_0\},\leqno(3.10)$$
\textit{and say that $\mathcal{L}$ is asymptotically admissible with respect to a given asymptotic expansion}
 $$
      f(x)=a_1\phi_1(x)+\dots+
a_n\phi_n(x)+o(\phi_n(x)),\  x\to x_0,\leqno(3.11)$$
\textit{if its formal application to both sides of $(3.11)$ yields}  
   $$
 \mathcal{L}[f(x)]=a_1\mathcal{L}[\phi_1(x)]+\dots+
a_m\mathcal{L}[\phi_m(x)]+o\big(\mathcal{L}[\phi_m(x)]\big),\  x\to x_0,\leqno(3.12)_1$$
\textit{wherein}
$$\begin{cases}
      \mathcal{L}[\phi_1(x)]\gg \dots \gg \mathcal{L}[\phi_n(x)],\  x\to x_0,\ \textit{after suppression}\\
       \textit{of all the terms $\equiv 0$ on some neighborhood of}\  x_0.
\end{cases}\leqno(3.12)_2$$

\textit{An alternative locution for an asymptotically-admissible $\mathcal{L}$ is ``$\mathcal{L}$ is formally applicable to the  asymptotic expansion". The spoken-of neighborhood of $x_0$ may well be one-sided.}

\vspace{5pt}
 We exhibit two simple examples clarifying the above definition;  in each of them the standard operator $d/dx$ is asymptotically admissible according to Definition 3.1 and inconsistencies  would occur without suppression of the identically-zero terms:
$$\begin{cases}
    \ \  f_1(x):=x^2+\log x +1+x^{-1}+e^{-x}, \ x>0, \\ \\
     \begin{cases}
      f_1(x)=x^2+\log x +1+x^{-1}+o(x^{-1}),\ x\to +\infty, \\
      f'_1(x)=2x+x^{-1}-x^{-2}+o(x^{-2}),\ x\to +\infty. 
\end{cases}
      \end{cases} \leqno (3.13)$$
         $$\begin{cases}
    \ \  f_2(x):=\log x +1+\sqrt{x}+x^2, \ x>0, \\ \\
     \begin{cases}
      f_2(x)=\log x +1+\sqrt{x}+o(x),\ x\to 0^+, \\
      f'_2(x)=x^{-1}+\frac{1}{2}x^{-1/2}+o(1),\ x\to 0^+. 
\end{cases}
      \end{cases} \leqno (3.14)$$

       \vspace{10pt}
    \centerline{\bf 3-A.\ The approach through the Wronskians, based on} 
     \centerline{\bf Levin's theorem on hierarchies}   
     
    \vspace{10pt}
   A mere rereading of Proposition 2.2 gives
   
   \vspace{10pt}
  {\bf Proposition 3.1.} \ \textit{Referring to a T.A.S. of class $AC^{n-1}[T,x_0[$ consider the operators}
    $$
    \mathcal{L}_{\phi_{i_1},\dots,\phi_{i_k}}u:=W( \phi_{i_1}, \dots ,\phi_{i_k},u),\
    1\leq i_1<i_2<\dots <i_k\leq n;\ 1\leq k\leq n-1,  \leqno (3.15) $$
   \textit{which are $k$th-order linear differential operators whose leading coefficients never vanish on a left deleted neighborhood of $x_0$. Then all these operators are asymptotically admissible with respect  to relation $(3.9)$ viewed as an asymptotic expansion with zero remainder; and this means that each relation}
    $$
    \mathcal{L}_{\phi_{i_1},\dots,\phi_{i_k}}f=\ \  \sum^{i=1,\dots ,n}_{i \neq i_j \ \forall j}\  a_iW( \phi_{i_1}, \dots ,\phi_{i_k},\phi_i)  \leqno (3.16) $$
 \textit{is again an  asymptotic expansion at $x_0$ with zero remainder. For instance we have the identities:}
      $$\begin{aligned}
   \mathcal{L}_{\phi_k}f & =a_1W(\phi_k, \phi_1)+\dots +a_{k-1}W(\phi_k, \phi_{k-1})+    \\     & + a_{k+1}W(\phi_k, \phi_{k+1})+\dots +a_nW(\phi_k, \phi_n), 
 \end{aligned}  \leqno (3.17)$$
\textit{wherein}
   $$\begin{aligned}
\label{}
  W(\phi_k, \phi_1) &\gg W(\phi_k, \phi_2)\gg \dots \gg W(\phi_k, \phi_{k-1}) \gg \\  &\gg W(\phi_k, \phi_{k+1})\gg \dots \gg W(\phi_k, \phi_n), \ x \to x_0^-,  
  \end{aligned} \leqno (3.18) $$
\textit{for each fixed $k, 1\leq k\leq n-1, n\geq 3.$}  (For $n=2$ the chain $(3.18)$ has only one term).
 \textit{And we also have the identities:}
 $$
    \mathcal{L}_{\phi_h,\phi_k}f=\ \  \sum ^{i=1,\dots ,n}_{ i \neq h,k} \ a_iW( \phi_h,\phi_k,\phi_i),  \leqno(3.19) $$
 \textit{wherein}
     $$\begin{cases}
    W(\phi_h,\phi_k, \phi_1)\gg W(\phi_h,\phi_k, \phi_2)\gg \dots \gg \cr 
    \gg W(\phi_h,\phi_k, \phi_{h-1}) \gg  W(\phi_h,\phi_k, \phi_{h+1})\gg \dots \gg \cr 
    \gg  W(\phi_h,\phi_k, \phi_{k-1})\gg  W(\phi_h,\phi_k, \phi_{k+1})\gg \dots \gg W(\phi_h,\phi_k, \phi_n),\ x \to x_0^- ,\end{cases} \leqno (3.20)$$
    \textit{for fixed $h,k: 1\leq h<k\leq n, n\geq 4.$}  (For $n=3$ the chain $(3.20)$ has only one term).
  
    \vspace{5pt}
    Proposition 3.1 gives rise to a first conjecture:
   
    \vspace{10pt}
 {\bf Conjecture A.} \textit{Referring to the asymptotic expansion $(1.1)$, or $(1.4)$ with $i<n$, there are many linear differential operators, namely $(3.15)$, which are likely to be formally applicable under reasonable hypotheses.}
 
 \vspace{10pt}
   \centerline{\bf 3-B. \ The special operators associated to canonical factorizations}
   
   \vspace{10pt}
   In the Wronskians in (3.15) a permutation of $(\phi_{i_1},\dots,\phi_{i_k})$ seems to be immaterial a sign apart, hence there are exactly $(2^n-2)$ essentially different operators of type (3.15): the number of the distinct subsets of $(\phi_1\dots,\phi_n)$ with cardinality $\tilde{k}: 1\leq \tilde{k} \leq n-1$. Now for $n\geq 3$ the object of our study, in a general formulation, involves a sequence of ``nested" operators: 
   $$
   \mathcal{L}_{\phi_{i_1}},\ \mathcal{L}_{\phi_{i_1},\phi_{i_2}},\ \dots 
    \ \mathcal{L}_{\phi_{i_1},\phi_{i_2},\dots ,\phi_{i_k}}, \leqno (3.21) $$
   where ``nested" refers to the inclusions of their kernels and the problem consists in finding sufficient, and possibly necessary, conditions for the validity of the set of asymptotic relations
   $$\begin{cases}
     f(x)=\sum_{i=1}^na_i\phi_i(x) + o\big(\phi_n(x)\big),  \\ \\
     \mathcal{L}_{\phi_{i_1}}[f(x)]=\sum_{i\neq i_1}^{i=1,\dots ,n}
     a_iW( \phi_{i_1},\phi_i; x)+o\big(\psi_1(x)\big), \\
     \dots \dots  \dots\\
     \mathcal{L}_{\phi_{i_1},\dots,\phi_{i_k}}[f(x)]= \sum^{i=1,\dots ,n}_{i \neq i_j \ \forall j}\  a_iW( \phi_{i_1}, \dots ,\phi_{i_k},\phi_i; x)+o\big(\psi_k(x)\big),
\end{cases} \leqno (3.22) $$
with  proper choices of the $\psi_i$'s. Once a subset $(\phi_{i_1},\ldots,,\phi_{i_k})$ has been fixed there is no a-priori reason to prefer one permutation of the $\phi_i$'s to another but it turns out that each {\it ordered} $k$-tuple $(\phi_{i_1},\ldots,,\phi_{i_k})$ is linked to a special factorization of $\mathcal{L}_{\phi_{i_1},\dots,\phi_{i_k}}$, possibly valid on a neighborhood of $x_0$ smaller than $[T, x_0[$ and calculations can be successfully carried out only under proper integrability assumptions on the coefficients of the factorization, hence the order of the $\phi_i$'s is not immaterial. Now a generic factorization of  $\mathcal{L}_{\phi_1,\dots ,\phi_n}$, say (2.3), assumed valid on $[T,x_0[$,  involves the differential operators
 $$
 r_0(x)u; \ r_1(x)( r_0(x)u)'; \ r_2(x)\big[r_1(x)( r_0(x)u)'\big]'  \ \dots \leqno (3.23) $$
which we label as ``\textit{weighted derivatives of orders $0,1,2$ etc. with respect to the weights $(r_0,r_1,\dots, r_n)$"} in preference to the (some-times used) generic locutions of ``\textit{quasi-derivatives or generalized derivatives"} with no reference to the $n$-tuples of weights. For convenience we include the operator of order zero. Operators (3.23) are not always linked to operators of the type in (3.15) nor they preserve the hierarchy of the $\phi_i$'s but the two C.F.'s highlighted in Proposition 2.1 yield two sequences of differential operators of orders $0,1,2,\dots ,n-1$ which are strictly related to operators in (3.15) and preserve the hierarchy; these operators were the core of the asymptotic theory in the case of real-power expansions [5; 6] hence they deserve a special attention and, as a matter of fact, the most meaningful results of our theory are based on them.

Referring to the factorization of type (I) in (2.47), with the $p_i$'s in (2.43), we define the differential operators acting on $AC^{n-1}[T,x_0[$: 
    $$\begin{cases}
  L_0u:=p_0(x)u;\  L_ku:=p_k\big[p_{k-1}(\dots(p_0u)'\dots)'\big]',\ 1\leq k\leq n;  \\     L_nu\equiv  L_{\phi_1,\dots,\phi_n}u,\end{cases} \leqno (3.24) $$
     which satisfy the recursive formula 
     $$
   L_ku:=p_k(x)(L_{k-1}u)',\  1\leq k\leq n. \leqno (3.25) $$
      
      And referring to the factorization of type (II) in (2.39), with the $q_i$'s in (2.35), we  define the differential operators acting on $AC^{n-1}[T,x_0[$: 
      $$\begin{cases}
 M_0u:=q_0(x)u;\ M_ku:=q_k(x)\big[q_{k-1}(x)(\dots (q_0(x)u)'\dots )'\big]',\ 1\leq k \leq n; \\  M_nu\equiv  L_{\phi_1,\dots,\phi_n}u,\end{cases} \leqno (3.26)$$
   which satisfy the recursive formula 
  $$
       M_ku:=q_k(x)(M_{k-1}u)',\ 1\leq k\leq n. \leqno (3.27) $$

Now representations (2.40) and (2.47) imply that:
$$\begin{cases}
{\rm ker}\ L_k={\rm span}\  (\phi_n,\phi_{n-1}, \dots ,\phi_{n-k+1}),\\ 
{\rm ker}\ M_k={\rm span}\  (\phi_1, \dots ,\phi_k),\ 1\le k\le n-1;
\end{cases}\leqno (3.28)$$
 hence there exist never-vanishing functions $\widetilde{p_k},\ \widetilde{q_k}$ such that:
 $$\begin{cases}
L_ku=\widetilde{p_k}\cdot W(\phi_n,\phi_{n-1}, \dots ,\phi_{n-k+1},u),\\ 
 M_ku=\widetilde{q_k}\cdot W(\phi_1, \dots ,\phi_k,u),\ 1\le k\le n-1.
 \end{cases}\leqno (3.29)$$
   
It follows that $L_k$ and $M_k$ preserve the hierarchy (2.23), namely we have the following asymptotic scales
$$
L_k[\phi_1(x)]\gg L_k[\phi_2(x)]\gg\dots \gg L_k[\phi_{n-k}(x)], x\to x_0^-,  \leqno (3.30) $$
$$
M_k[\phi_{k+1}(x)]\gg M_k[\phi_{k+2}(x)]\gg \dots \gg M_k[\phi_n(x)],\ x\to x_0^-,
\leqno (3.31) $$
for each fixed $k, 0\leq k\leq n-2$. For $k=0$ they respectively reduce to
$$
p_0(x)\phi_1(x)\gg p_0(x)\phi_2(x)\gg \dots \gg p_0(x)\phi_n(x),\ x\to x_0^-,
\leqno (3.32) $$
$$
q_0(x)\phi_1(x)\gg q_0(x)\phi_2(x)\gg \dots \gg q_0(x)\phi_n(x),\ x\to x_0^-,
\leqno (3.33) $$
both equivalent to (2.23). Hence, applying each $n$-tuple of operators $L_k$ and $M_k$, $0\leq k\leq n-1$, to (3.9) yields again asymptotic expansions with zero remainders and in this sense we may say that ``\textit{the asymptotic expansion $(3.9)$ is formally differentiable $(n-1)$ times with respect to the $n$-tuples of weights $(p_0,\dots ,p_{n-1})$ and $(q_0,\dots ,q_{n-1})$}" neglecting the $n$th-order weighted derivatives which yield identically-zero expressions. The above discussion leads to the following 
 
  \vspace{10pt}
  {\bf Conjecture B} \  (Particularization of Conjecture A). 
 \textit{For each chosen C.F. of $L_{\phi_1,\dots,\phi_n}$ of type either {\rm (I)} or {\rm (II)} at $x_0$,}
 $$  
L_{\phi_1,\dots,\phi_n}u\equiv r_n\big[r_{n-1}(\dots (r_0u)'\dots)'\big]' \  \ \forall \ u\in AC^{n-1}[T,x_0[\ ,\leqno (3.34) $$
\textit{there exists a linear subspace $\mathcal{D}\in AC^{n-1}[T,x_0[$, such that:}

(i)\  $\mathcal{D}\underset {\neq}{\supset} \textnormal{span}\  (\phi_1,\dots,\phi_n),$ 

(ii)\ \textit {each $f \in \mathcal{D}$ has an asymptotic expansion of type $(1.1)$ which is formally differentiable $(n-1)$ times with respect to the n-tuples of weights  $(r_0,r_1,\dots, r_{n-1})$.}

  \vspace{10pt}
The problem consists in finding out  analytic conditions characterizing the elements of $\mathcal{D}$ for a C.F. of type (I) or (II) separately. The foregoing approach suggests a smallness condition involving the quantity 
$L_{\phi_1,\dots,\phi_n}[f(x)]$ which is $\equiv0$ whenever the remainder in the  expansion is.           

\vspace{10pt}
   \centerline{\bf 3-C. \ The coefficients of an asymptotic expansion  with zero remainder}
   
  \vspace{10pt}
A third fact we wish to investigate is the possible expressions of the coefficients of an asymptotic expansion alternatively to the recurrent  formulas (1.3), so generalizing (3.4). It is clear from the study of  polynomial expansions in [4] that the C.F. of type (I) is of no use to this end  whereas the right approach is via a C.F. of type (II) by establishing a link between the coefficients of (3.9) and the limits of the weighted derivatives.

 \vspace{10pt}
{\bf Proposition 3.2}\  (The coefficients of an asymptotic expansion with zero remainder).\ \textit{Referring to the T.A.S. in Proposition $2.4$ and  to the special factorization $(2.39)$ the following facts hold true for the differential operators $M_k$ in $(3.26)$: }

  (I)\ \textit {The $M_k$'s satisfy the following relations:}
 $$
 \textnormal{ker}\ M_k=\ \textnormal{span}\  ( \phi_1, \dots ,\phi_k),\ 1\leq k\leq n; \leqno (3.35)$$
 $$
 M_k[\phi_{k+1}(x)]\equiv \epsilon_k=constant=\pm1,\ 1\leq k\leq n-1; \leqno (3.36)$$
 $$\begin{cases}\displaystyle
 M_k[\phi_h(x)]=\epsilon_{h,k}\cdot \int_x^{x_0}\frac{1}{q_{k+1}}\dots \int^{x_0}\frac{1}{q_{h-1}} =o(1),\ x\to x_0^-; \\ \\
 \epsilon_{h,k}=constant=\pm1,\ 1\leq k \leq h-2,\ h\leq n. \end{cases} \leqno(3.37)$$
$$
M_ku\equiv \epsilon_k\frac{W(\phi_1, \dots ,\phi_k,u)}{W(\phi_1, \dots ,\phi_k,\phi_{k+1})},\ 1\leq k\leq n-1. \leqno (3.38)$$

 (II)\  \textit{For a fixed k, $1\leq k \leq n$, we have the logical equivalence:}
$$M_{k-1}[f(x)]\equiv \epsilon_{k-1}\cdot a_k=\  \textit{constant on some interval J} \leqno (3.39)$$
 \textit{iff}
$$f(x)=a_1\phi_1(x)+\dots +a_k\phi_k(x) \  \textit{on J for some constants}\  a_i, \leqno (3.40)$$
$a_k$ \textit{being the same as in $(3.39)$ and $\epsilon_{k-1}$ as in $(3.36)$.}

\textit{If} (3.39)-(3.40) {\it hold true on a left neighborhood of $x_0$ then the following limits exist as finite numbers and}
$$\epsilon_{h-1}\cdot a_h=\lim_{x\to x_0^-} M_{h-1}[f(x)],\ 1\leq h\leq k,\leqno (3.41)$$
\textit{where, for h=k, $(3.41)$ is the identity $(3.39)$.}

\vskip5pt
(III)\ \textit {In the special case where all the Wronskians in $(2.24)$ are strictly positive then the constants in $(3.36)$-$(3.37)$ have the values:}
$$
\epsilon_k=1,\ \ \epsilon_{h,k}=(-1)^{h+k+1}.\leqno(3.42)$$

\vspace{5pt}
 We stress that the equivalence ``$(3.39)\Leftrightarrow (3.40)"$ is an algebraic fact based on (3.35)-(3.36) whereas the inference ``(3.39)-(3.40) $\Rightarrow$ (3.41)" is an asymptotic property whose validity requires that $( \phi_1, \dots ,\phi_k)$ be an asymptotic scale at $x_0$ and that the operators $M_k$ be defined as specified.
 Proposition 3.2 suggests the following
   
    \vspace{10pt}
 {\bf Conjecture C.} \ \textit{ If all the limits in $(3.41)$ 
exist as finite numbers for some function $f$ sufficiently regular on a left deleted neighborhood of $x_0$ then an asymptotic expansion}
  $$
   f(x)=a_1\phi_1(x)+\dots +a_k\phi_k(x)+o\big(\phi_k(x)\big),\ x\to x_0^-, 
   \leqno (3.43)$$
  \textit{holds true matched to other expansions obtained by formal applications of the operators $M_1,\ldots,M_{k-1}$. Moreover it is worth investigating if the validity of the sole last relation in $(3.41)$, i.e. for $h=k$, implies the validity of the other relations.} 
   
    \vspace{5pt}
We shall give complete answers to Conjectures B and C in \S\S4,5.

\vspace{10pt}
   \centerline{\bf 3-D.\ An heuristic approach via L'Hospiral's rule.}
    
\vspace{10pt}
There is another way to arrive at Conjecture C by the elementary use of L'Hospital's rule. The following calculations on the limits in (1.3), if legitimate, would yield:
$$\begin{cases}\displaystyle
a_1=\lim_{x\to x_0^-}f/\phi_1; \\ \displaystyle
a_2=\lim_{x\to x_0^-}{f-a_1\phi_1\over \phi_2}\equiv\lim_{x\to x_0^-}
\frac{(f/\phi_1)-a_1\ \big(=o(1)\big)}
{\phi_2/\phi_1\hskip30pt \big(=o(1)\big)}\ \overset{H}{=}\ \lim_{x\to x_0^-}
\frac{(f/\phi_1)'}{(\phi_2/\phi_1)'}\equiv \lim_{x\to x_0^-}M_1[f(x)];\\ \\ \displaystyle
a_3=\lim_{x\to x_0^-}{f-a_1\phi_1-a_2\phi_2\over \phi_3}\equiv\lim_{x\to x_0^-}
\frac{(f/\phi_1)-a_1-a_2(\phi_2/\phi_1)}{\phi_3/\phi_1}\ \overset{H}{=}\ \\ \\ \displaystyle \overset{H}{=}\ \lim_{x\to x_0^-}\frac{(f/\phi_1)'-a_2(\phi_2/\phi_1)'}
{(\phi_3/\phi_1)'}\equiv\lim_{x\to x_0^-}\frac{(f/\phi_1)'\big/(\phi_2/\phi_1)'-a_2}
{(\phi_3/\phi_1)'\big/(\phi_2/\phi_1)'}\ \overset{H}{=}\ \\ \\ \displaystyle \overset{H}{=}\ \lim_{x\to x_0^-}\frac{\left(\ (f/\phi_1)'\big/(\phi_2/\phi_1)'\ \right)'}
{\left(\ (\phi_3/\phi_1)'\big/(\phi_2/\phi_1)'\ \right)'}\equiv \lim_{x\to x_0^-}M_2[f(x)];\ \ \ and\ so\ on.
\end{cases}\leqno(3.44)$$

Such kind of manipulations may seem artificial and awkward from an elementary viewpoint and it is by no means obvious that iterating the procedure yields the relations in (3.41) for $h\ge4$ as well. In one of the two algorithms presented in \S8 (Proposition 8.1) it will be shown that the procedure is quite natural in the context of formal differentiation of an asymptotic expansion and that it actually leads to (3.41) for all values of $h$.

\vspace{20pt}
\centerline{\textbf {4. The first factorizational approach}}

\vspace{10pt}
We start from the "unique" C.F. of our operator $L_{\phi_1,\dots,\phi_n}$  on the interval $[T,x_0[$ of type (I) at $x_0$, i.e. identity (2.47) with conditions (2.45)-(2.46) and the $p_i$'s satisfying the same conditions as do the $q_i$'s in (2.36). We consider the fundamental system (2.48).  By (2.49) the ordered $n$-tuple $(P_{n-1},\dots ,P_0)$ is an asymptotic scale at $x_0^-$ but it cannot coincide (constant factors apart) with the given scale $(\phi_1,\dots,\phi_n)$ as (2.41) and (2.49) are incompatible. However  (2.23) and (2.49) imply that the two scales are linked by the following relations
   $$
   \phi_i(x)\sim b_iP_{n-i}(x),\ x\to x_0^-,\ 1\leq i\leq n, \leqno (4.1) $$
  with suitable nonzero constants $b_i$, hence 
  $$\begin{cases}
       \phi_i(x)=b_iP_{n-i}(x)+\sum_{j=i+1}^n \beta_{i,j}P_{n-j}(x),\ 1\leq i\leq n-1, \\  
       \phi_n(x)=b_nP_0(x),\end{cases} \leqno (4.2)$$ 
and viceversa 
   $$
       P_0(x)=\frac{1}{b_n}\phi_n(x), \ P_i(x)=\frac{1}{b_{n-i}} \phi_{n-i}(x)+\sum_{j=n-i+1}^n {\widetilde{\beta}}_{i,j}\phi_j(x),\ 1\leq i\leq n-1,   \leqno (4.3)$$
with suitable constants $\beta_{i,j},{\widetilde{\beta}}_{i.j}$.
  
 In this approach the appropriate differential operators to be used are the $L_k$'s defined in (3.24) and here are some elementary properties of these operators.
    
    \vspace{10pt}
   {\bf Lemma 4.1.}  \textit{The following relations are checked at once:}
    $$
 \textnormal{ker}\  L_k= \begin{cases}
 \textnormal{span}\  (P_0,P_1,\dots , P_{k-1})    \\
  \textnormal{span}\ (\phi_n,\phi_{n-1}, \dots ,\phi_{n-k+1}),\ 1\leq k\leq n; 
\end{cases} \leqno (4.4)$$
$$
L_k[P_k(x)]\equiv 1,\  0\leq k\leq n-1; \leqno (4.5) $$
$$
L_k[P_i(x)]\equiv \int_T^x \frac{dt_{k+1}}{p_{k+1}(t_{k+1})} \dots \int_T^{t_{i-1}}\frac{dt_i}{p_i(t_i)},
\ 0\leq k<i \leq n-1; \leqno (4.6) $$
$$
L_k[P_i(x)] \ll L_k[P_{i+1}(x)], \ x\to x_0^-,\ 0\leq k\leq i \leq n-2. \leqno (4.7) $$

\textit{Hence we have the following chains of asymptotic relations:}
$$\begin{cases}
       L_0[P_0(x)] \ll L_0[P_1(x)] \ll \dots \ll L_0[P_{n-1}(x)], \\
       L_1[P_1(x)] \ll L_1[P_2(x)] \ll \dots \ll L_1[P_{n-1}(x)], \\
      L_2[P_2(x)] \ll L_2[P_3(x)] \ll \dots \ll L_2[P_{n-1}(x)], \ \ \ x\to x_0^-,  \\
       \dots \dots \dots\\
      L_{n-2}[P_{n-2}(x)] \ll L_{n-2}[P_{n-1}(x)].  \end{cases} \leqno(4.8)$$
      
      \textit{The first chain in $(4.8)$ coincides with the second chain in $(2.49)$ apart from the ordering and the multiplicative factor $p_0(x)$. As the first term in each chain is the constant "$1$"  all the other terms diverge to $+\infty$.}
      
    \vspace{10pt}
{\bf Lemma 4.2.}\ \textit{If a solution $\phi$ of $ L_{\phi_1,\dots,\phi_n}u=0$ 
     satisfies the asymptotic relation} 
     $$\phi(x) \sim cP_i(x),\ x\to x_0^-, \leqno (4.9)$$
     \textit{for some $i\in \{ 0,1,\dots ,n-1\}$ and some nonzero constant $c$ then the following relations hold true:}
     $$
     L_k[\phi(x)]  \sim cL_k[P_i(x)],\ x\to x_0^-,\ 0\leq k\leq i\leq n-1; 
     \leqno (4.10)$$
     $$
      L_k[\phi(x)] \equiv 0,\ i+1\leq k\leq n. \leqno (4.11)$$
      
      \textit{Moreover:}
      $$
       L_k[\phi_{n-i}(x)] \equiv \begin{cases}
       b_{n-k},\ 0\leq i=k\leq n-1,  \\
       0,\ \ \ \ \ \ 0\leq i<k,\end{cases} \leqno (4.12) $$
      \textit{with the $b_i$'s defined in $(4.1)$. It follows from $(4.1)$ and $(4.10)$ that all relations in $(4.8)$ hold true after replacing $P_i$ by $\phi_{n-i}$ hence we have the asymptotic scales:}
      $$\begin{cases}
        L_0[\phi_1(x)] \gg L_0[\phi_2(x)] \gg \dots \gg L_0[\phi_n(x)], \\
       L_1[\phi_1(x)] \gg L_1[\phi_2(x)] \gg \dots \gg L_1[\phi_{n-1}(x)], \\
       L_2[\phi_1(x)] \gg L_2[\phi_2(x)] \gg \dots \gg L_2[\phi_{n-2}(x)],
        \\    \dots \dots  \ldots\\
      L_{n-2}[\phi_1(x)] \gg L_{n-2}[\phi_2(x)].  \end{cases}\ x\to x_0^-,  \leqno(4.13)$$
      
     \textit{Last, with the $b_i$'s defined in $(4.1)$, we have the identity:}
     $$
     L_ku\equiv b_{n-k}\frac{W( \phi_n,\phi_{n-1}, \dots ,\phi_{n-k+1},u)}{W( \phi_n,\phi_{n-1}, \dots ,\phi_{n-k})},\ 1\leq k\leq n-1. \leqno (4.14) $$
 
 \vspace{10pt}
{\bf Lemma 4.3.}\  \textit{Any function $f\in AC^{n-1}[T,x_0[$ admits of a  representation of type:}
     $$\begin{cases}
      f(x)=c_1\phi_1(x)+\dots +c_n\phi_n(x)+ \\ \\ \displaystyle
     + \frac{1}{p_0(x)}  \int_T^x \frac{1}{p_1} \dots \int_T^{t_{n-2}} \frac{1}{p_{n-1}}\int_T ^{t_{n-1}}
     \frac {L_{\phi_1,\dots,\phi_n}[f(t)]}{p_n(t)} dt, \ x\in [T,x_0[\ , 
     \end{cases} \leqno (4.15) $$  
     \textit{with suitable constants $c_i$. From $(4.6), (4.12)$ and $(4.15)$ we infer the following representations of the weighted derivatives of $f$ with respect to the weight functions $(p_0,\dots ,p_n)$}:
 $$\begin{cases}\displaystyle
L_k[f(x)] =c_1L_k[\phi_1(x)]+\dots +c_{n-k}L_k[\phi_{n-k}(x)]+  \\ \\ \displaystyle
   + \int_T^x \frac{dt_{k+1}}{p_{k+1}(t_{k+1})} \dots \int_T^{t_{n-2}} \frac{dt_{n-1}}{p_{n-1}(t_{n-1})}\int_T^{t_{n-1}} \frac{L_{\phi_1,\dots,\phi_n}[f(t)]}{p_n(t)}  dt, \\ \\
   for\ x\in [T,x_0[\ ;\ 0\leq k\leq n-2;
  \end{cases} \leqno (4.16)$$ 
  $$\begin{cases}\displaystyle
 L_{n-1}[f(x)]=c_1L_{n-1}[\phi_1(x)]+\int_T^x \frac {L_{\phi_1,\dots,\phi_n}[f(t)]}{p_n(t)} dt = \\ \\ \displaystyle
 \overset {(4.12)}{=}c_1b_1+ \int_T^x \frac{ L_{\phi_1,\dots,\phi_n}[f(t)]}{p_n(t)} dt,\ 
 x\in [T,x_0[\ . \end{cases}  \leqno (4.17) $$
 
 \textit{By $(4.13)$ the linear combination $\sum_{i=1}^{n-k}c_iL_k[\phi_i(x)]$ in the right-hand side of $(4.16)$ is in itself an asymptotic expansion at $x_0^-$ for each fixed $k$.}
 
  \vspace{5pt}
We shall now characterize various situations wherein relations (4.16)-(4.17) become asymptotic expansions. In the following two theorems we state separately three cases of a single claim lest a unified statement be obscure. The reader is referred to the first remark after next theorem to grasp the meaning of the differentiated asymptotic expansions which exhibit a special non-common phenomenon.
  
 \vspace{10pt}
{\bf Theorem 4.4}\ (Asymptotic expansions formally differentiable according to the C.F. of type (I) ).\ \textit{Let $f\in AC^{n-1}[T,x_0[$.}
 
 (I)\ \textit{The following are equivalent properties for a suitable constant $a_1$}:
 
 (i)\ \textit{The set of asymptotic relations}
 $$
 L_k[f(x)] =a_1L_k[\phi_1(x)]+o\big(L_k[\phi_1(x)]\big),\ x\to x_0^-;\ 0\leq k\leq n-1. 
 \leqno (4.18)$$
 
(ii)\ \textit{The single asymptotic relation}
  $$
  L_{n-1}[f(x)]=a_1b_1+o(1),\ x\to x_0^-,\  \textit{ with $b_1$ defined in $(4.1)$,} 
  \leqno (4.19)$$
  \textit{which is the explicit form of the relation in $(4.18)$ for $k=n-1$.}
  
(iii)\ \textit{The improper integral}
  $$
   \int_T^{x_0} \frac {L_{\phi_1,\dots,\phi_n}[f(t)]}{p_n(t)} dt \ \ \ \textit{converges.}\leqno (4.20)$$
   
   \textit{Under condition $(4.20)$ we have the representation formula:}
   $$
    L_{n-1}[f(x)]=a_1b_1- \int_x^{x_0} \frac{L_{\phi_1,\dots,\phi_n}[f(t)]}{p_n(t)} dt,\ x\in [T,x_0[\ . \leqno(4.21) $$
    
  (II)\ \textit{For a fixed $i\in \{2,\dots ,n\}$ the following are equivalent properties for suitable constants $a_i$} (the same in each set of conditions):
    
(iv)\ \textit{The set of asymptotic expansions as $x\to x_0^-$}:
    $$ \begin{cases}
      L_k[f(x)]=a_1L_k[\phi_1(x)]+\dots +a_iL_k[\phi_i(x)]+
     o\big( L_k[\phi_i(x)]\big),\ 0\leq k \leq n-i;  \\ \\
   \begin{aligned}
          L_{n-i+h}[f(x)]&= a_1L_{n-1+h}[\phi_1(x)]+\dots +a_{i-h}L_{n-i+h}[\phi_{i-h}(x)]+ \\   &  +o(1);\ 0\leq h \leq i-1.\end{aligned}
\end{cases} \leqno(4.22) $$

(v)\ \textit{The second group of asymptotic expansions in $(4.22)$, i.e.}
 $$\begin{cases}\begin{aligned}
          L_{n-i+h}[f(x)]&= a_1L_{n-1+h}[\phi_1(x)]+\dots +a_{i-h}L_{n-i+h}[\phi_{i-h}(x)]+ \\ &+o(1),\ x\to x_0^-;\ 0\leq h \leq i-1,\end{aligned}
          \end{cases} \leqno(4.23) $$
\textit{where we point out that the last meaningful term in the right-hand side is a constant.}

 (vi)\ \textit{The following  improper integral, involving `` $i$'' iterated integrations,}
 $$
  \int_T^{x_0}\frac{1}{p_{n-i+1}}\int_{t_{n-i+1}}^{x_0}\frac{1}{p_{n-i+2}}\dots
 \int_{t_{n-2}}^{x_0}\frac{1}{p_{n-1}}    
  \int_{t_{n-1}}^{x_0} \frac{ L_{\phi_1,\dots,\phi_n}[f(t)]}{p_n(t)} dt \  \ \ \textit{converges.}\leqno (4.24)$$

\textit{Under condition $(4.24)$  we have the representation formula}:
$$\begin{cases}
  L_{n-i}[f(x)]= a_1L_{n-i}[\phi_1(x)]+\dots +a_iL_{n-i}[\phi_i(x)]+\\ \\ \displaystyle  
    +(-1)^i \int_x^{x_0}\frac{1}{p_{n-i+1}}\int_{t_{n-i+1}}^{x_0}\frac{1}{p_{n-i+2}}\dots
 \int_{t_{n-2}}^{x_0}\frac{1}{p_{n-1}} \int_{t_{n-1}}^{x_0} \frac{ L_{\phi_1,\dots,\phi_n}[f(t)]}{p_n(t)} dt,  \end{cases} \leqno (4.25)$$
 \textit{for $x\in [T,x_0[$, as well as the corresponding formulas for the functions $L_{n-i+h}[f(x)]$ with $0\leq h\leq i-1,$ obtained by suitable differentiations of $(4.25)$}: see remark 3 below.

\vspace{10pt}
{\it Remarks}.\ 1. Relations in (4.22) may be read as follows. The first relation, involving $L_0$, is equivalent to the asymptotic expansion
$$
 f(x)=a_1\phi_1(x)+\dots +a_i\phi_i(x)+o\big(\phi_i(x)\big), \ x\to x_0^-, \leqno (4.26) $$
and the relations involving $L_k$, with $1\leq k\leq n-i$, state that (4.26) can be formally differentiated $(n-i)$ times in the sense of formally applying the operators $L_k$ to the remainder in (4.26). In so doing one arrives at the expansion
$$
L_{n-i}[f(x)]= a_1L_{n-i}[\phi_1(x)]+\dots +a_iL_{n-i}[\phi_i(x)]+o(1), \ x\to x_0^-, \leqno (4.27) $$
where $L_{n-i}[\phi_i(x)] \equiv $ constant. The process of formal differentiation, from the order $(n-i+1)$ up to $(n-1)$, goes on according to the following rule: in (4.27) and in each expansion in (4.23) the \underline{last} term is constant and is \underline{lost} after one further weighted differentiation while the remainder preserves its simple growth-order estimate of ``$o(1)"$. So the first $(n-i+1)$ expansions, i.e. those involving $L_0,L_1,\dots ,L_{n-i},$ have the same number of meaningful terms whereas  each of the other $(i-1)$ expansions is deprived of the last meaningful term at each successive differentiation. We rewrite more explicitly the expansions in (4.22) to better highlight the dynamics of this process:
$$\begin{cases}
       f(x)=a_1\phi_1(x)+\dots +a_i\phi_i(x)+o\big(\phi_i(x)\big), \\ \\
        L_1[f(x)]=a_1L_1[\phi_1(x)]+\dots +a_iL_1[\phi_i(x)]+o\big(L_1[\phi_i(x)]\big), \\
       \dots \dots \\
       L_{n-i}[f(x)]= a_1L_{n-i}[\phi_1(x)]+\dots +
      a_i  \underset {constant} {\underbrace{L_{n-i}[\phi_i(x)]}}+o(1),  \\
   L_{n-i+1}[f(x)]= a_1L_{n-i+1}[\phi_1(x)]+\dots +a_{i-1} \underset {constant}
  {\underbrace{L_{n-i+1}[\phi_{i-1}(x)]}}+o(1),   \\
  \dots \dots \\
  L_{n-2}[f(x)]= a_1L_{n-2}[\phi_1(x)]+a_2\underset {constant}
  {\underbrace{L_{n-2}[\phi_2(x)]}}+o(1), \\ 
  L_{n-1}[f(x)]= a_1 \underset {constant} {\underbrace{L_{n-1}[\phi_1(x)]}} +o(1). 
\end{cases} \leqno (4.28)$$

The loss of the last meaningful term, where it occurs, is caused by formula (4.12) for $i=k-1$ which, after renaming the indexes, reads
$$
L_{n-i+h}[\phi_{i-h+1}(x)]\equiv 0. \leqno (4.29)$$

Notice that in the second group of expansions in (4.28) the meaningful terms disappear one after one in reversed order if compared with Taylor's formula.

2.\ It is shown in \S7, after the proof of Theorem 4.4, that the set (4.23) is not equivalent in general to the single relation
$$
L_{n-i}[f(x)]= a_1L_{n-i}[\phi_1(x)]+\dots +
      a_i L_{n-i}[\phi_i(x)]+o(1),\  x\to x_0^-, \leqno (4.30) $$
as in part (I) of the theorem (case $i=1$).

3.\ Suitable weighted differentiations of (4.25) yield integral representations of the remainders in the differentiated expansions of orders greater than $(n-i)$ and these representations are numerically meaningful. On the contrary, if $i<n,$ then successive integrations of (4.25) contain some constants not uniquely defined hence the corresponding representations are of no numerical use without additional information on $f.$

\vspace{5pt}
For $i=n$ the subset of (4.22) involving the operators $L_k, 1\leq k\leq n-i,$ is empty and here is an explicit and expanded statement.
  
  \vspace{10pt}
{\bf Theorem 4.5}\ (The case $i=n$ in Theorem 4.4). \  \textit{For $f\in AC^{n-1}[T,x_0[$ the following are equivalent properties:}
 
 (i)\   \textit{The set of asymptotic expansions as $x\to x_0^-$ for suitable constants 
 $a_1,\ldots,a_n$}:
 $$\begin{cases}
       f(x)=a_1\phi_1(x)+\dots +a_n\phi_n(x)+o\big(\phi_n(x)\big),  \\ \\
        L_k[f(x)]= a_1L_k[\phi_1(x)]+\dots +
      a_{n-k}  \underset {constant} {\underbrace{L_k[\phi_{n-k}(x)]}}+ o(1),
       \ 1\leq k\leq n-1,\end{cases}\leqno (4.31)$$
{\it where the \underline{last} term in each expansion is \underline{lost} in the successive expansion.}

(ii)\  \textit{The improper integral}
 $$
  \int_T^{x_0} \frac{1}{p_1} \dots \int _{t_{n-2}}^{x_0} \frac{1}{p_{n-1}}\int_{t_{n-1}}^{x_0}
     \frac{ L_{\phi_1,\dots,\phi_n}[f(t)]}{p_n(t)} dt \ \ \ \textit{converges.} 
     \leqno (4.32)$$
 
 (iii) \textit {There exist $n$ real numbers $a_1,\dots ,a_n$ and a function 
$\Phi_n$ Lebesgue-summable on $[T,x_0[$ such that}
 $$\begin{cases}     
f(x)=a_1\phi_1(x)+\dots +a_n\phi_n(x)+ \\ \\ \displaystyle
+\frac{(-1)^n}{p_0(x)} \int_x^{x_0} \frac{1}{p_1}\int_{t_1}^{x_0} \frac{1}{p_2} \dots \int _{t_{n-2}}^{x_0} \frac{1}{p_{n-1}}\int_{t_{n-1}}^{x_0}
    \Phi_n(t) dt,\  x\in [T,x_0[\ .\end{cases}\leqno (4.33)$$
 
 \textit{If this is the case $\Phi_n$ is determined up to a set of measure zero and}
$$
\Phi_n(x)=\frac{1}{p_n(x)} L_{\phi_1,\dots,\phi_n}[f(x)] \ \textit{a.e. on} \ \  [T,x_0[\ . 
\leqno (4.34)  $$
    
     \vspace{10pt}
  The phenomenon described in the above theorems is intrinsic in the theory; it occurs even in the seemingly elementary case of real-power expansions, [6; Thm. 4.2-(ii), p. 181, and formula (7.2), p. 195], where the asymptotic scale enjoys the most favourable algebraic properties. This type of formal differentiation of an asymptotic expansion does not frequently occur in the literature though the results in this section show that it is one of the possible natural situations. An instance (not inserted in a general theory) is to be found in a paper by Schoenberg [13; Thm. 3, p. 258] and refers to the asymptotic expansion
   $$
   f(x)=a_1x^{-1}+a_2x^{-2}+\dots +a_nx^{-n}+O(x^{-n-1}),\ x\to +\infty. \leqno (4.35)$$  

\vspace{20pt}
\centerline{\textbf{5. The second factorizational approach and estimates of the remainder}}  

 \vspace{10pt} 
Now we face our problem starting from a C.F. of type (II) at $x_0$. Referring to Proposition 2.4 the most natural choice is the special C.F. of  $L_{\phi_1,\dots,\phi_n}$ in (2.39), with the $q_i$'s in (2.35) and satisfying conditions (2.36). According to Conjectures B and C in \S3 we shall characterize a set of asymptotic expansions, involving the operators $M_k$ defined in (3.26), wherein each coefficient of the first expansion may be found by an independent limiting process instead of the recursive formulas (1.3), and the existence of the sole last coefficient implies the existence of all the preceding coefficients. 
 
In this new context a representation of the following type is appropriate for any function  $f\in AC^{n-1}[T,x_0[:$
$$\begin{cases}
    f(x)=c_1\phi_1(x)+\dots +c_n\phi_n(x)+ \\ \\ \displaystyle
     + \frac{1}{q_0(x)}  \int_T^x \frac{1}{q_1} \dots \int_T^{t_{n-2}} \frac{1}{q_{n-1}}\int_T ^{t_{n-1}}
     \frac{L_{\phi_1,\dots,\phi_n}[f(t)]}{q_n(t)}dt\ , \ x\in [T,x_0[\ , 
     \end{cases}\leqno (5.1) $$  
with suitable constants $c_i$. 
 Applying the operators $M_k$ to (5.1) we get the following representations of the weighted derivatives of $f$ with respect to the weight functions $(q_0,\dots ,q_n):$  
     $$\begin{cases}
  M_k[f(x)] =c_{k+1}M_k[\phi_{k+1}(x)]+\dots +c_nM_k[\phi_n(x)]+\\ \\ \displaystyle 
   + \int_T^x \frac{dt_{k+1}}{q_{k+1}(t_{k+1})} \dots \int_T^{t_{n-2}} \frac{dt_{n-1}}{q_{n-1}(t_{n-1})}\int_T^{t_{n-1}} \frac{L_{\phi_1,\dots,\phi_n}[f(t)]}{q_n(t)}dt, \ 0\leq k\leq n-1.\end{cases}\leqno (5.2)$$

\vskip5pt
 \underline{Warning!}\ To simplify formulas and to leave no ambiguity about the signs of the involved quantities we assume throughout this section that  the Wronskians in (2.24) are strictly positive. 
 
 Hence, by (3.42)  $\epsilon_k=1$ and the last relation in (5.2) explicitly is
 $$
 M_{n-1}[f(x)]=c_n+\int_T^x \frac{L_{\phi_1,\dots,\phi_n}[f(t)]}{q_n(t)}dt,\ 
 x\in [T,x_0[\ . \leqno (5.3)  $$

By (3.31) the ordered linear combination in (5.2),
$$
\sum^n_{i=k+1}c_iM_k[\phi_i(x)], \leqno (5.4) $$
is an asymptotic expansion at $x_0^-$ for each fixed $k,\ 0\leq k\leq n-1.$

Unlike \S4 we first state here the result concerning a complete  asymptotic expansion, i.e. of type (1.1), because it is the most expressive result in this paper and characterizes the simple circumstance that $M_{n-1}[f(x)]=a_n+o(1)$ via a set of $n$ asymptotic expansions. 
Always refer to Proposition 3.2 for properties of the $M_k$'s.

 \vspace{10pt}
{\bf Theorem 5.1}\ (Complete asymptotic expansions formally differentiable according to a C.F. of type (II) ).\ \textit{Let our T.A.S. be such that all the Wronskians in $(2.24)$ are strictly positive and let $f\in AC^{n-1}[T,x_0[$.}
 
 (I)\ \textit{The following are equivalent properties:}
 
(i) \  \textit{There exist $n$ real numbers $a_1,\dots ,a_n$ such that:} 
 $$
 f(x)=a_1\phi_1(x)+\dots +a_n\phi_n(x)+o\big(\phi_n(x)\big),\ x\to x_0^-; 
 \leqno (5.5) $$
  $$\begin{cases}\begin{aligned}
   M_k[f(x)] &=a_{k+1}M_k[\phi_{k+1}(x)]+\dots +a_nM_k[\phi_n(x)]+ \\ 
   &+o\big(M_k[\phi_n(x)]),\ x\to x_0^-,\ 1\leq k\leq n-1;\end{aligned}\end{cases}\leqno (5.6)$$
  \textit{where the \underline{first} term in each expansion is \underline{lost} in the successive expansion.} 
  
  Notice that the relation that would be obtained in (5.6) for $k=0$ differs from relation in (5.5) by the common factor $q_0(x)$.
 
(ii) \ \textit{All the following limits exist as finite numbers:}
  $$
  \lim_{x\to x_0^-} M_k[f(x)]\equiv a_{k+1}, \ 0\leq k\leq n-1, \leqno (5.7) $$
  \textit{where the $a_k$'s coincide with those in $(5.5)$.}
 
(iii) \ \textit{The single last limit in $(5.7)$ exists as a finite number, i.e.}
 $$
  \lim_{x\to x_0^-} M_{n-1}[f(x)]\equiv a_n,\leqno (5.8) $$
  \textit{and $(5.8)$ is nothing but the relation in $(5.6)$ for $k=n-1$ which reads  $M_{n-1}[f(x)]=a_n+o(1), x\to x_0^-.$}
  
(iv) \ \textit{The improper integral}
 $$
  \int_T^{\to x_0} \frac{L_{\phi_1,\dots,\phi_n}[f(t)]}{q_n(t)}dt \ \ \  
 \textit{converges,} \leqno (5.9)  $$
\textit{and automatically also the iterated improper integral} 
$$
    \int_T^{x_0} \frac{1}{q_1}\int^{x_0} \frac{1}{q_2} \dots \int^{x_0} \frac{1}{q_{n-1}}\int
    ^{x_0}\frac{L_{\phi_1,\dots,\phi_n}[f(t)]}{q_n(t)}dt  \ \ \ \textit{converges.} 
    \leqno (5.10)  $$
    
(v) \ \textit{There exist $n$ real numbers $a_1,\dots ,a_n$ and a function 
$\Psi_n$ Lebesgue-summable on $[T,x_0[$ such that}
$$
f(x)=\sum_{i=1}^na_i\phi_i(x)+\frac{(-1)^n}{q_0(x)}\int_x^{x_0} \frac{1}{q_1}\dots \int_{t_{n-2}}^{x_0} \frac{1}{q_{n-1}}\int_{t_{n-1}}^{x_0}\Psi_n(t)dt,\ x\in [T,x_0[\ , 
\leqno (5.11)  $$
\textit{where we remind that, by $(2.35)$, $1/q_0(x)=\phi_1(x)$. In this case $\Psi_n$ is determined up to a set of measure zero and}
$$
\Psi_n(x)=\frac{1}{q_n(x)} L_{\phi_1,\dots,\phi_n}[f(x)] \ \textit{a.e. on} \  [T,x_0[\ . \leqno (5.12)  $$

 (II)\   \textit{Whenever properties in part $(I)$ hold true we have integral representation formulas for the remainders}
$$\begin{cases}\displaystyle
 R_0(x):=f(x)-\sum_{i=1}^na_i\phi_i(x), \\  \displaystyle
  R_k(x):=M_k[f(x)] -\sum _{i=1}^{n-k}a_{k+i}M_k[\phi_{k+i}(x)], \ 1\leq k\leq n-1, 
\end{cases} \leqno (5.13) $$
\textit{namely}:
$$
R_0(x)=\frac{(-1)^n}{q_0(x)}\int_x^{x_0} \frac{1}{q_1}\dots \int_{t_{n-2}}^{x_0} \frac{1}{q_{n-1}}\int_{t_{n-1}}^{x_0}\frac{L_{\phi_1,\dots,\phi_n}[f(t)]}{q_n(t)}dt,
\leqno (5.14) $$
$$
R_k(x)=(-1)^{n+k}\int_x^{x_0} \frac{1}{q_{k+1}}\dots \int_{t_{n-2}}^{x_0} \frac{1}{q_{n-1}}\int_{t_{n-1}}^{x_0}\frac{L_{\phi_1,\dots,\phi_n}[f(t)]}{q_n(t)}dt, 
 \leqno (5.15)$$
\textit{for $x\in [T,x_0[,  1\leq k\leq n-1.$ From $(5.14)$ we get the following estimate of $R_0$ wherein the order of smallness with respect to $\phi_n$ is made more explicit than in Theorem $4.5$} (formula in (2.40) for $i=n$ is used):
$$
|R_0(x)|\leq|\phi_n(x)|\cdot \underset{t\geq x}{\textnormal{sup}} \ 
\left| \int_t^{x_0}\frac{L_{\phi_1,\dots,\phi_n}[f(\tau)]}{q_n(\tau)}d\tau \right|, \ x\in [T,x_0[\ .
 \leqno (5.16)  $$
 
 \textit{Under the stronger hypothesis of absolute convergence for the improper integral we get:}
 $$
 |R_0(x)|\leq|\phi_n(x)|\cdot 
 \int_x^{x_0}\frac{\left| L_{\phi_1,\dots,\phi_n}[f(t)] \right|}{|q_n(t)|}dt , \ x\in [T,x_0[\ .
 \leqno (5.17)  $$
 
 \textit{Similar estimates can be obtained for the $R_k$'s.}
 
 \vspace{10pt}
 {\it Remarks. 1.}\ As noticed in [6; Remark 1 after Thm. 4.1, pp. 179-180] the remarkable inference $"(iii)\Rightarrow(ii)"$ is true for the special operator $M_{n-1}$ stemming out from a C.F. of type (II) at $x_0$ but not for any $(n-1)$th-order differential operator originating from an arbitrary factorization of $L_{\phi_1,\dots,\phi_n}$. 
 
2.\ Condition (5.9) involves the sole coefficient $q_n$ which admits of the explicit expression in (2.35) in terms of $\phi_1,\dots ,\phi_n:$
$$
q_n= W( \phi_1, \dots ,\phi_n)/W( \phi_1, \dots ,\phi_{n-1}); \leqno (5.18) $$
hence (5.9) can be rewritten as
$$
\int_T^{\to x_0}\frac{W\big( \phi_1(t), \dots ,\phi_{n-1}(t)\big)}{ W\big( \phi_1(t), \dots ,\phi_n(t)\big)}
 L_{\phi_1,\dots,\phi_n}[f(t)] dt  \ \ \ \textit{converges.}\leqno (5.19) $$

For $n=2$ the ratio inside the integral equals $\phi_1/W(\phi_1,\phi_2)$ and we reobtain the result in [7; condition (5.15), p. 265].

3.\ In Theorem 4.5, generally speaking, no such estimates as in (5.16)-(5.17) can be obtained due to the divergence of all the improper integrals in (4.33) if the innermost integral is factored out.

4.\ Referring to the elementary characterizations in (1.3) of the coefficients $a_k$   Theorem 5.1 changes the perspective: in (1.3) the $a_k$'s are defined recursively whereas in (5.7) each $a_k$ has its own independent expression and, moreover, the existence of $a_n$, as the limit in (5.8), implies the existence of $a_1,\dots ,a_{n-1}.$ 

\vspace{10pt}
 In the following result about incomplete expansions formal differentiation is in general legitimate a number of times less than the "length" of the expansion (see Remark 2 after the statement).
 
 \vspace{10pt}
{\bf Theorem 5.2}\ (A result on incomplete asymptotic expansions). \  \textit{Let our T.A.S. be such that all the Wronskians in $(2.24)$ are strictly positive and let $f\in AC^{n-1}[T,x_0[$.}

 (I) \textit{For a fixed $i\in \{2,\dots ,n-1\}$ the following are equivalent properties}:
 
(i)\   \textit{There exist $i$ real numbers $a_1,\dots ,a_i$ such that}:
 $$
 f(x)=a_1\phi_1(x)+\dots +a_i\phi_i(x)+o\big(\phi_i(x)\big),\ x\to x_0^-; 
  \leqno (5.20) $$
 $$\begin{cases}\begin{aligned}
   M_k[f(x)]&=a_{k+1}M_k[\phi_{k+1}(x)]+\dots +a_iM_k[\phi_i(x)]+\\ 
   &+o\big(M_k[\phi_i(x)]\big),\ x\to x_0^-;\ 1\leq k\leq i-1.
   \end{aligned}\end{cases}\leqno (5.21)$$
 
(ii)\   \textit{All the following limits exist as finite numbers}:
$$
\lim_{x\to x_0^-} M_k[f(x)]\equiv a_{k+1}, \ 0\leq k\leq i-1, \leqno (5.22) $$
\textit{where the $a_k$'s coincide with those in} {\rm (5.20)-(5.21).}

(iii)\   \textit{The single last limit in $(5.22)$  exists as a  finite number, i.e. }
$$
\lim_{x\to x_0^-} M_{i-1}[f(x)]\equiv a_i, \leqno (5.23) $$
\textit{and $(5.23)$ coincides with the relation in $(5.21)$ for $k=i-1.$}

(iv)\ \textit{The improper integral}
$$
\int_T^{\to x_0}\frac{dt_i}{q_i(t_i)}\int_T^{ t_i} \frac{dt_{i+1}}{q_{i+1}(t_{i+1})} \dots \int_T^{t_{n-1}} \frac{L_{\phi_1,\dots,\phi_n}[f(t)]}{q_n(t)}dt \ \ \ \textit{converges,}\leqno(5.24)$$
\textit{and automatically also the iterated improper integral}
$$
\int_T^{x_0}\frac{1}{q_1}\dots \int_{t_{i-1}}^{x_0}\frac{1}{q_i}\int_T^{t_i} \frac{1}{q_{i+1}}\dots \int_T^{t_{n-1}}\frac{L_{\phi_1,\dots,\phi_n}[f(t)]}{q_n(t)}dt 
\ \ \  \textit{converges.}\leqno(5.25)$$

 (II) \textit{For $i=1$  the theorem simply asserts that the asymptotic relation}
$$
f(x)=a_1\phi_1(x)+o\big(\phi_1(x)\big),\  x\to x_0^-, \leqno (5.26) $$
\textit{holds true for some real number $a_1$ iff the improper integral}
$$
\int_T^{\to x_0}\frac{1}{q_1}\int_T^{t_1} \frac{1}{q_2}\dots \int_T^{t_{n-1}}\frac{L_{\phi_1,\dots,\phi_n}[f(t)]}{q_n(t)}dt 
\ \ \  \textit{converges.}\leqno(5.27)$$

 \vspace{10pt}
 {\it Remarks. 1.}\ We shall see in the proof of Theorem 5.2, formula (7.44), that the representations of the quantities $M_k[f(x)], 0\leq k\leq i-1,$ contain some unspecified constants not determinable through the sole condition (5.24) which, for this reason, grants neither explicit representations nor numerical estimates of the remainders of the expansions in (5.20)-(5.21). 
 
 2.\ As concerns estimates of the quantities $M_k[f(x)]$ for $i\leq k\leq n-1$, the situation is as follows. Formula (7.43), given in the proof, reads:
 $$
 M_i[f(x)]=\underset{I(x)} {\underbrace{\int_T^{x} \frac{dt_{i+1}}{q_{i+1}(t_{i+1})} \dots \int_T^{t_{n-1}} \frac{L_{\phi_1,\dots,\phi_n}[f(t)]}{q_n(t)}dt }}+c+o(1), \leqno (5.28) $$
 for some constant $c$. If, as $x\to x_0^-,\ I(x)$ converges to a real number then we may apply Theorem 5.2 with $i$ replaced by $i+1$; but if it is unbounded and oscillatory no asymptotic relation more expressive than (5.28) can be obtained generally speaking. On the contrary a favourable situation occurs when it is known a priori that $I(x)$ either converges or diverges to $\pm \infty$ and the corrresponding estimates are reported in Theorem 6.3.

 \vspace{10pt}
{\bf Theorem 5.3}\ (The analogue of Theorems 5.1-5.2 with ``$O$"-estimates). \  \textit{Let our T.A.S. be such that all the Wronskians in $(2.24)$ are strictly positive,  let $f\in AC^{n-1}[T,x_0[$ and let $i\in \{2,\dots ,n\}$ be fixed. The following are equivalent properties:}
 
(i)\   \textit{There exist $(i-1)$ real numbers $a_1,\dots ,a_{i-1}$ such that}:
 $$
 f(x)=a_1\phi_1(x)+\dots +a_{i-1}\phi_{i-1}(x)+O\big(\phi_i(x)),\ x\to x_0^-; 
  \leqno (5.29) $$
 $$\begin{cases}\begin{aligned}
   M_k[f(x)]&=a_{k+1}M_k[\phi_{k+1}(x)]+\dots +a_{i-1}M_k[\phi_{i-1}(x)]+\\
   &+O\big(M_k[\phi_i(x)]\big),\ x\to x_0^-;\ 1\leq k\leq i-1.
   \end{aligned}\end{cases}\leqno (5.30)$$
 
(ii)\   \textit{All the following relations hold true}:
$$\begin{cases}
\lim_{x\to x_0^-} M_k[f(x)]\equiv a_{k+1}, \ 0\leq k\leq i-2; \\ \\
 M_{i-1}[f(x)]=O(1),\ x\to x_0^-;
 \end{cases}\leqno (5.31) $$
\textit{where the $a_k$'s coincide with those in} (5.29)-(5.30).

(iii)\   \textit{It holds true the single last relation in $(5.31)$, i.e. }
$$
M_{i-1}[f(x)]=O(1),\ x\to x_0^-; \leqno (5.32)$$
(iv)\ \textit{We have the following estimate instead of condition $(5.24)$:}
$$
\int_T^x\frac{dt_i}{q_i(t_i)}\int_T^{ t_i} \frac{dt_{i+1}}{q_{i+1}(t_{i+1})} \dots \int_T^{t_{n-1}} \frac{L_{\phi_1,\dots,\phi_n}[f(t)] }{q_n(t)}dt=O(1),\ x\to x_0^-.\leqno(5.33)$$

\textit{For $i=n$ condition $(5.32)$ reads}
$$
 \int_T^x\frac{L_{\phi_1,\dots,\phi_n}[f(t)] }{q_n(t)}dt=O(1),\ x\to x_0^-,\leqno(5.34)$$
\textit{and representation {\rm (5.11)-(5.12)} must be replaced by}

$$f(x)=\sum_{i=1}^{n-1}a_i\phi_i(x)
+\frac{(-1)^{n-1}}{q_0(x)}\int_x^{x_0} \frac{1}{q_1}\int_{t_1}^{x_0}\frac{1}{q_2}\dots \int_{t_{n-2}}^{x_0} \frac{1}{q_{n-1}}
\int_T^{t_{n-1}}\frac{L_{\phi_1,\dots,\phi_n}[f(t)]}{q_n(t)}dt.\leqno (5.35)$$

\textit{For $i=1$  the theorem simply asserts that the asymptotic relation}
$$
f(x)=O(\phi_1(x)),\  x\to x_0^-, \leqno (5.36) $$
\textit{holds true iff }
$$
\int_T^x\frac{1}{q_1}\int_T^{t_1} \frac{1}{q_2}\dots \int_T^{t_{n-1}}\frac{L_{\phi_1,\dots,\phi_n}[f(t)]}{q_n(t)}dt=O(1),\  x\to x_0^-. \leqno(5.37)$$
   
  \vspace{10pt}
  {\it  An outlook on the theory developed so far}.\ We suggest a way of visualizing what our theory is all about. Referring, say, to the situations characterized in Theorem 5.2 we have an asymptotic scale of the type:
  $$
   \frac{1}{q_0(x)}\gg \frac{1}{q_0(x)} \int_x^{x_0} \frac{1}{q_1}\gg\ldots\gg
    \frac{1}{q_0(x)} \int_x^{x_0} \frac{1}{q_1} \dots \int _{t_{i-2}}^{x_0} \frac{1}{q_{i-1}}, \ x\to x_0^-, \leqno (5.38) $$
    where the $q_i$'s are continuous and everywhere-nonzero functions on some interval $[T,x_0[$, and are interested in the validity of an  asymptotic expansion of the type:
    $$\begin{cases}\displaystyle
    f(x)=\frac{a_1}{q_0(x)}+ \frac{a_2}{q_0(x)} \int_x^{x_0} \frac{1}{q_1}+\frac{a_3}{q_0(x)} \int_x^{x_0} \frac{1}{q_1}\int_x^{x_0} \frac{1}{q_2}+\ldots+\\ \\ \displaystyle
    +\frac{a_i}{q_0(x)}\left[ \int_x^{x_0} \frac{1}{q_1} \ldots \int _{t_{i-2}}^{x_0} \frac{1}{q_{i-1}}\right]\cdot\big[1+o(1)\big], \ x\to x_0^-. \end{cases}\leqno (5.39) $$
    
    Theorem 5.2  gives characterizations of the set formed by (5.39) and the following expansions obtained in a quite natural way:
    $$\begin{cases}\displaystyle
    M_1\big[f(x)\big]\equiv q_1(x)\big(q_0(x)f(x)\big)'=-\frac{a_2}{q_1}-\frac{a_3}{q_1(x)} \int_x^{x_0} \frac{1}{q_2}+\ldots+\\ \\ \displaystyle
    -\frac{a_i}{q_1(x)}\left[ \int_x^{x_0} \frac{1}{q_2} \dots \int _{t_{i-2}}^{x_0} \frac{1}{q_{i-1}}\right]\cdot\big[1+o(1)\big], \ x\to x_0^-;
    \\ \ldots \ldots \ldots \\  \displaystyle
    M_{i-1}\big[f(x)\big]\equiv q_i(x)\big(q_{i-1}(x)(\ldots(q_0(x)f(x))'\ldots)'\big)'=\\ =(-1)^{i-1}a_i+o(1),\  x\to x_0^-.
\end{cases}\leqno (5.40) $$
    
   The formal derivations of (5.40) from (5.39) may seem a triviality but it is not an automatic fact and we have tied up our theory with the concepts of Chebyshev systems and canonical factorizations, useful in other contexts. Moreover in some applications the asymptotic scale is explicitly given whereas the expressions of the coefficients $q_i$ or $p_i$ of the canonical factorizations, as given by formulas (2.35) or (2.43), are unmanageable  even for small values of $n$ and only some properties of them can be detected and used. In other applications, e.g. when a function $f$ is defined as a solution of a functional equation, it may happen that the asymptotic scale is implicitly defined and only the principal parts of the $\phi_i$'s are known; in such cases there is no searching out the expressions of  the $q_i$'s and $p_i$'s, but $q_n$, as the ratio in (5.18), might be indirectly known and this would let us decide whether or not Theorem 5.1 applies.

\vspace{20pt}
\centerline{\textbf {6. Absolute convergence and solutions of differential inequalities}}

\vspace{10pt}
The foregoing theory becomes particularly simple when the involved improper integrals are absolutely convergent and still more expressive for a function $f$ satisfying the $n$th-order differential inequality:
$$
L_{\phi_1,\dots,\phi_n}[f(x)]\geq 0\ \ \textit{a.e. on}\ \ [T,x_0[\ . \leqno (6.1) $$

If\ \  $W\big(\phi_1(x),\dots,\phi_i(x)\big)>0$ on $[T,x_0[, \ 1\leq i\leq n$, this is a subclass of the so-called ``\textit{generalized convex functions with respect to the $($complete extended Chebyshev$)$ system} $(\phi_1,\dots,\phi_n)$", and we make this assumption, as in the preceding section,  to simplify relations involving the operators $M_k$ and to state precise inequalities for the remainders. The nice result stated in the next theorem claims that: if such a function  admits of an asymptotic expansion (1.1) then this expansion is automatically differentiable $(n-1)$ times in the senses of both relations (4.31) and (5.6). 

\vspace{10pt}
{\bf Theorem  6.1}\ (Complete asymptotic expansions).\ \textit{If all the Wronskians in $(2.24)$ are strictly positive and if $f\in AC^{n-1}[T,x_0[$ satisfies $(6.1)$ then the following are equivalent properties:}
 
 (i) \textit{There exist $(n-1)$ real numbers $a_1,\dots ,a_{n-1}$ such that:}
$$
 f(x)=a_1\phi_1(x)+\dots +a_{n-1}\phi_{n-1}(x)+O\big(\phi_n(x)\big),\ x\to x_0^-. \leqno (6.2) $$
 
 (ii) \textit{There exist $n$ real numbers $a_1,\dots ,a_n$ such that:}
 $$
f(x)=a_1\phi_1(x)+\dots +a_{n-1}\phi_{n-1}(x)+a_n\phi_n(x)+o\big(\phi_n(x)\big),\ x\to x_0^-. \leqno (6.3) $$

(iii) \textit{The following set of asymptotic expansions holds true:}
$$ \begin{cases}\begin{aligned}
       L_k[f(x)]&= a_1L_k[\phi_1(x)]+\dots +
      a_{n-k}  \underset {constant} {\underbrace{L_k[\phi_{n-k}(x)]}}+  \\ 
      &+o(1),\ x\to x_0^-, \ 0\leq k\leq n-1;\ \  \text{see (4.31).} 
\end{aligned}\end{cases}\leqno (6.4)$$

(iv)  \textit{The following set of asymptotic expansions holds true:}
 $$\begin{cases}\begin{aligned}
   M_k[f(x)]&=a_{k+1}M_k[\phi_{k+1}(x)]+\dots +a_nM_k[\phi_n(x)]+\\
  &+o\big(M_k[\phi_n(x)]\big),\ x\to x_0^-;\ 0\leq k\leq n-1;\ \  \text{see (5.5)-(5.6)}.\end{aligned}\end{cases}\leqno (6.5)$$
   
   (v) \textit{The following integral condition is satisfied:}
   $$
   \int_T^{x_0} \frac{1}{p_1} \dots \int _{t_{n-2}}^{x_0} \frac{1}{p_{n-1}}\int_{t_{n-1}}^{x_0}
     \frac{1}{p_n(t)} L_{\phi_1,\dots,\phi_n}[f(t)] dt <+\infty;\ \ \text {see (4.32).} \leqno (6.6)$$
     
     (vi) \textit{The following integral condition is satisfied:}
   $$
 \int_T^{ x_0} \frac{1}{q_n(t)} L_{\phi_1,\dots,\phi_n}[f(t)] dt < + \infty, \ \text{see (5.9) and (5.19).} 
 \leqno (6.7)  $$
 
 \textit{To this list we may obviously add the other properties in Theorem $5.1$.}
 
 \textit{If this is the case the remainder $R_0(x)$ of the expansion in $(6.3)$ admits of the two representations on $[T,x_0[$:
 $$ \begin{cases}
    R_0(x)&=\displaystyle\frac{(-1)^n}{p_0(x)} \int_x^{x_0} \frac{1}{p_1}\dots \int _{t_{n-2}}^{x_0} \frac{1}{p_{n-1}}\int_{t_{n-1}}^{x_0}\frac{L_{\phi_1,\dots,\phi_n}[f(t)]}{p_n(t)}dt,\\ \\   
 &= \displaystyle\frac{(-1)^n}{q_0(x)}\int_x^{x_0} \frac{1}{q_1}\dots \int_{t_{n-2}}^{x_0} \frac{1}{q_{n-1}}\int_{t_{n-1}}^{x_0}\frac{L_{\phi_1,\dots,\phi_n}[f(t)]}{q_n(t)}dt\ , \end{cases}\leqno (6.8)  $$
 whence it  follows that 
 $$
 (-1)^nR_0(x)\ge0\ \ \forall\ x\in [T,x_0[\ , \leqno(6.9)$$
 and that both of the following two functions are decreasing on $[T,x_0[$}:
 $$
 (-1)^nR_0(x)p_0(x)\equiv -R_0(x)\phi_n(x),\ \ \ (-1)^nR_0(x)q_0(x)\equiv -R_0(x)\phi_1(x).\leqno(6.10)$$
 
 \vspace{10pt}
 In addition to  the equivalence (iii)$\Leftrightarrow$(iv) stated in Theorem 6.1, there is another remarkable circumstance wherein  the two types of formal differentiations are simultaneously admissible namely when the convergence of the pertinent improper integrals is absolute.     
 
  \vspace{10pt} 
  {\bf Theorem  6.2.}\ \textit{For} $f\in AC^{n-1}[T,x_0[$ \textit{ the following integral conditions are equivalent}: 
    $$
\int_T^{x_0} \frac{1}{p_1} \dots \int _{t_{n-2}}^{x_0} \frac{1}{p_{n-1}}\int_{t_{n-1}}^{x_0}
     \frac{1}{p_n(t)}\big| L_{\phi_1,\dots,\phi_n}[f(t)] \big| dt <+\infty ; \leqno (6.11)$$
   $$\begin{cases}  \displaystyle
     \int_T^{x_0}P(t) \big| L_{\phi_1,\dots,\phi_n}[f(t)] \big| dt <+\infty ,\ \ \   \textit{where}   \\  \\ \displaystyle 
          P(t):=\frac{1}{p_n(t)}\int_T^t\frac{dt_{n-1}}{p_{n-1}(t_{n-1})}\dots \int_T^{t_3}\frac{dt_2}{p_2(t_2)} \int_T^{t_2} \frac{dt_1}{p_1(t_1)}\ ;
          \end{cases}\leqno (6.12)$$
\vspace{5pt}
  $$
  \int_T^{ x_0} \frac{\big|L_{\phi_1,\dots,\phi_n}[f(t)]\big|}{q_n(t)} dt \equiv  \int_T^{x_0}\left|\frac{W\big( \phi_1(t), \dots ,\phi_{n-1}(t)\big)}
  {W\big( \phi_1(t), \dots ,\phi_n(t)\big)} L_{\phi_1,\dots,\phi_n}[f(t)]\right| dt < + \infty. \leqno (6.13) $$
  
  \textit{Hence each of these three conditions implies both sets of asymptotic expansions $(4.31)$ and} (5.5)-(5.6). (Here the signs of the Wronskians are 
 immaterial.)
 
 \vspace{10pt}
 {\bf Open problem 1.}\ In \S7 we give an indirect proof of the equivalence   ``(6.12)$\Leftrightarrow$(6.13)"  based on Theorem 6.1; a more refined statement would be:
 $$
 P(x)\sim c\ \frac{W\big( \phi_1(x), \dots ,\phi_{n-1}(x)\big)}{ W\big(\phi_1(x), \dots ,\phi_n(x)\big)},\ x\to x_0^-, \leqno(6.14)$$
 for some constant $c\ne 0$. For $n=2$ this is quite elementary and we also found a proof for $n=3$; but for the time being we leave this minor question as an open problem. 
  
   \vspace{5pt}
 Using Theorems 4.4 and 5.2  we can also get the analogues of Theorems 6.1-6.2 for incomplete asymptotic expansions left apart the integral representations of the remainders but with meaningful estimates for weighted derivatives of orders $\geq i$. We give here a simplified statement wherein all asymptotic relations refer to $x\to x_0^-$ of course.
 
 \vspace{10pt}
{\bf Theorem 6.3}\ (Incomplete asymptotic expansions).\ \textit{Let the Wronskians in $(2.24)$ be strictly positive, let $f\in AC^{n-1}[T,x_0[$ satisfy $(6.1)$ and let $i\in\{1,\dots ,n-1\}$ be fixed. Then the following are equivalent properties:}
$$
 f(x)=a_1\phi_1(x)+\dots +a_{i-1}\phi_{i-1}(x)+O\big(\phi_i(x)\big); \leqno (6.15) $$
 $$
 f(x)=a_1\phi_1(x)+\dots +a_{i-1}\phi_{i-1}(x)+a_i\phi_i(x)+o\big(\phi_i(x)\big); \leqno (6.16) $$
 $$\begin{cases} \begin{aligned}
            L_k[f(x)] &=a_1L_k[\phi_1(x)]+\dots +a_iL_k[\phi_i(x)]+ \\  
     &+o( L_k[\phi_i(x)]),\ 0\leq k \leq n-i;\\ \\
            L_{n-i+h}[f(x)]&= a_1L_{n-1+h}[\phi_1(x)]+\dots + \\
       &  +a_{i-h}L_{n-i+h}[\phi_{i-h}(x)]+o(1),\ 0\leq h \leq i-1,\end{aligned} \end{cases} \leqno(6.17) $$
 {\rm (which last relations are written in $(4.28)$ in an expanded form);}
  
$$ \begin{cases}
\begin{aligned}
   M_k[f(x)] &=a_{k+1}M_k[\phi_{k+1}(x)]+\dots +a_iM_k[\phi_i(x)]+\\ 
   & +o\big(M_k[\phi_i(x)]\big),\ 0\leq k\leq i-1; \\ \\
   M_k[f(x)] &=O\left(\int_T^x\frac{1}{q_{k+1}}\dots \int_T^{t_{n-1}}\frac{L_{\phi_1,\dots,\phi_n}[f(t)]}{q_n(t)}dt\right), i\leq k\leq n-2; \\ \\
   M_{n-1}[f(x)]&=O\left(\int_T^x\frac{L_{\phi_1,\dots,\phi_n}[f(t)]}{q_n(t)}dt\right); \end{aligned} \end{cases} \leqno(6.18) $$
 $$
   \int_T^{x_0}\frac{1}{p_{n-i+1}}\int_{t_{n-i+1}}^{x_0}\frac{1}{p_{n-i+2}}\dots
   \int_{t_{n-1}}^{x_0} \frac{L_{\phi_1,\dots,\phi_n}[f(t)]}{p_n(t)}dt<+\infty;
    \leqno (6.19)$$
 $$ \begin{cases}   \displaystyle 
     \int_T^{x_0}P(t)  L_{\phi_1,\dots,\phi_n}[f(t)]  dt <+\infty, \ \  \ \textit{where} 
     \\ \\ \displaystyle 
          P(t):=\frac{1}{p_n(t)}\int_T^t\frac{dt_{n-1}}{p_{n-1}}\dots \int_T^{t_{n-i+2}} \frac{dt_{n-1+1}}{p_{n-1+1}}\ \  \textit{if}\ \ i\geq 2;\end{cases}\leqno (6.20)$$

 $$
 \int_T^{x_0}\frac{dt_i}{q_i}\int_T^{t_i}\frac{dt_{i+1}}{q_{i+1}} \dots \int_T^{t_{n-1}} \frac{L_{\phi_1,\dots,\phi_n}[f(t)]}{q_n(t)}dt<+\infty. \leqno(6.21)$$

 \textit{To the foregoing list we may obviously add property {\rm (ii)} or property {\rm (v)} in Theorem $4.4$ and properties {\rm (ii)-(iii)} in Theorem $5.1$. For $i=1$ the first group of expansions in $(6.18)$ reduces to relation in $(5.26)$.}
 
 \vspace{5pt}
 As pointed out in Remark $1$ after Theorem 5.2 the ``$O"$-estimates in (6.18) are meaningful whenever all the involved integrals diverge as $x\to x_0$ i.e. whenever the asymptotic expansion in (6.16) cannot be improved by adding more meaningful terms of the form $a_{i+j}\phi_{i+j}(x)$.  As soon as one of these integrals converges to a real number as $x\to x_0$ then we may apply the theorem with a greater value of $i$.  
 
 \vspace{10pt}
 {\it Remark.} \ In Theorem 6.1 the two types of formal differentibility $1,2,\ldots,n-1$ times are equivalent facts whereas it is not so for a generic $f$ such that 
 $L_{\phi_1,\dots,\phi_n}[f(x)]$ changes sign on each deleted left neighborhood of $x_0$. The equivalence has been proved for polynomial expansions [4] and for real-power expansions [6] in an indirect way by expressing the two sets of differentiated expansions as suitable sets of expansions involving the standard operators $d^k/dx^k$; these new sets of expansions make evident that what we called ``weak formal differentiability'' indeed is a weaker property  than  what we called ``strong formal differentiability''. The same circumstance occurs for a general two-term expansion [7; Remarks, p. 261] but is not a self-evident fact. In each of these three cases direct proofs could be also provided working on the corresponding integral conditions. Hence in these cases the locutions of ``weak or strong formal differentiation'' are legitimate. But in the general theory for $n\ge3$ we face a nontrivial situation and state
 
 \vspace{5pt}
 {\bf Open problem 2.}\ \ For $n\ge3$  consider the two types of formal differentiability characterized in Theorems 4.5 and 5.1. Investigate whether or not property in Theorem 5.1 always implies that in Theorem 4.5 for any T.A.S. .
 
  \vspace{20pt}  
   \centerline{\bf 7. Proofs} 
 
   \vspace{10pt}  
   {\it Proof of Proposition} 2.1.\  For the equivalence of the two properties in (i) see Coppel [1; Prop.  3, p. 82]. \ "(i)$\Leftrightarrow$(ii)"   is proved in Levin [9; Thm.  2.1, p .66] where the interval $I$ is explicitly stated to be open not in the statement of the cited theorem but at the outset of \S 2 on p. 58;  "(ii)$\Leftrightarrow$(iii)" is the classical result by P\'olya [12]; "(i)$\Leftrightarrow$(iv)'' is the fundamental result by Trench [14]; "(i)$\Rightarrow$(v)" is to be found in [2; Thm. 2.2, p. 162] whereas the converse rests on the trivial fact that disconjugacy on $]a,b[$ is equivalent to disconjugacy on every compact subinterval of $]a,b[$.   \hfill $\Box$
 
 \vskip5pt
 {\it Proof of Proposition} 2.2.\ Part (I) is contained in [9; Th. 2.1, p. 66] with reverse numbering of the $\phi_i$'s whereas part (II) follows from [9; Lemma 2.6, pp. 63-64, and remarks on p. 67 concerning the hierarchies of the Wronskians], here again with reverse numbering of the $\phi_i$'s. Levin's results are valid for an open interval and this is stated explicitly at the outset of \S2 in [9; p. 58]; moreover, the tacit assumption of strict positivity of the functions forming the scale is agreed in a long list of notations and terminology in [9; \S1, p. 57, item 20°].    \hfill $\Box$
 
 \vskip5pt
 {\it Proof of Proposition} 2.3. \ (i)$\Rightarrow$(ii). Let $(\widetilde\phi_1,\ldots,\widetilde\phi_n)$ be an extension of  $(\phi_1,\ldots,\phi_n)$ of class $C^{n-1}]T-\epsilon,x_0[,\ \epsilon>0$, such that:
 $$
 W\big(\widetilde\phi_1,\ldots,\widetilde\phi_i\big)\ne0\ on\ ]T-\epsilon,x_0[,\ 1\le i\le n.\leqno (7.1)$$
 
 In particular we have:
 $$\begin{cases}
 W\big(\widetilde\phi_n,\widetilde\phi_{n-1},\ldots,\widetilde\phi_1\big)\ne0\\
  W\big(\widetilde\phi_{n-1},\ldots,\widetilde\phi_1\big) \ne0\end{cases}
  \ on\ ]T-\epsilon,x_0[,\leqno (7.2)$$
 and we may apply part (II) of Proposition 2.2 (regardless of the signs) because the second condition in (7.2) coincides with the condition in (2.16) for $r=1$. So we infer the inequalities:
 $$
 W\big(\widetilde\phi_n,\widetilde\phi_{n-1},\ldots,\widetilde\phi_i\big)\ne0\ on\ ]T-\epsilon,x_0[,\ 1\le i\le n,\leqno (7.3)$$
 which imply (2.27). Proposition 2.2 also implies all the claims in part (II).
 
 (i)$\Leftrightarrow$(iii). We refer to the standard definition of the concept of ``extended complete Chebyshev system on a generic interval $J$'', based on the maximum number of zeros for their linear combinations: see, e.g., [8; Ch. I]. A classical result  states the equivalence between an ordered $n$-tuple $(u_1,\ldots,u_n)$ forming such a system on $J$ and the strict positivity of the Wronskians $W(u_1,\ldots,u_i),\ 1\le i\le\ n,$ on $J$. This is proved, e.g., in [8; Ch. XI, Th. 1.1, p. 376] for a compact interval $J$ but the argument is valid for any interval as observed, e.g., by Mazure [10; Prop. 2.6]. This equivalence is a general fact involving only inequalities (2.24).
 
 (ii)$\Rightarrow$(iv). Here we are merely retracing the steps of the proof in [8; Ch. XI, Th 1.2, pp. 379-380] in a way that includes in one proof the expressions given in (2.31). First, inequalities (2.24) grant that the functions $w_i,\ 0\le i\le n-1$, are well defined on $[T,x_0[$ and satisfy (2.29): the second expression for $w_i,\ i\ge2$, in (2.30) is a classical identity: see [1; Lemma 4, p. 87] for a syntetic proof under our regularity assumptions. Moreover inequalities (2.24) and (2.27) together grant, by Proposition 2.2-(II), the asymptotic relations (2.14) hence:
 $$\begin{cases}
 \phi_2(x)/\phi_1(x)=o(1),\ x\to x_0^-,\\ \\ \displaystyle
 \frac{W( \phi_1(x), \dots ,\phi_{i-1}(x),\phi_{i+1}(x))}{W( \phi_1(x), \dots ,\phi_{i-1}(x),\phi_i(x))}=o(1),\ x\to x_0^-,\ 2\le i\le n-1.\end{cases}\leqno(7.4)$$
 
 This implies: the convergence of the improper integrals
 $$
 \int^{x_0}(\phi_2/\phi_1)';\ \ \int^{x_0}\left[\frac{W( \phi_1(t), \dots ,\phi_{i-1}(t),\phi_{i+1}(t))}{W( \phi_1(t), \dots ,\phi_{i-1}(t),\phi_i(t))}\right]' dt,\ 2\le i\le n-1,\leqno(7.5)$$
 the representations for $\phi_1, \phi_2$ and the identity
 $$
 \frac{W( \phi_1(x), \dots ,\phi_{i-1}(x),\phi_{i+1}(x))}{W( \phi_1(x), \dots ,\phi_{i-1}(x),\phi_i(x))}=\int_x^{x_0}w_i(t)dt,\ 2\le i\le n-1.\leqno(7.6)$$
 
Before using induction we prove the representation of $\phi_3$ to highlight the role of  (7.6). We have:
$$\begin{cases}\displaystyle
W(\phi_1,\phi_3)\big/W(\phi_1,\phi_2)=\int_x^{x_0}w_2(t)dt,\\
W(\phi_1,\phi_2)=-(\phi_1)^2w_1,\ \ W(\phi_1,\phi_3)=(\phi_1)^2\big(\phi_3/\phi_1)',\end{cases}\leqno(7.7)$$
whence
$$\begin{cases}\displaystyle
\big(\phi_3/\phi_1)'(x)=-w_1(x)\int_x^{x_0}w_2(t)dt,\\ \\ \displaystyle
\phi_3(x)/\phi_1(x)=c+\int_x^{x_0}w_1\int_{t_1}^{x_0}w_2(t)dt\overset{by (2.23)}{=}
\int_x^{x_0}w_1\int_{t_1}^{x_0}w_2(t)dt,\end{cases}\leqno(7.8)$$
which implies the representation of $\phi_3$ in (2.28). To prove the representations of $\phi_i$ for $4\le i\le n-1$ we proceed by induction supposing to have proved our inference (ii)$\Rightarrow$(iv) for any $i$-tuple forming a T.A.S. on $[T,x_0[$; hence our representations hold true for $\phi_1,\ldots,\phi_i$ and we must prove it for $\phi_{i+1}$. Putting
 $$
 \psi_k(x):=\big(\phi_{k+1}/\phi_1(x)\big)',\ 1\le k\le i,\leqno(7.9)$$
 we immediately infer from (1.5) and from (2.14) referred to $x\to x_0^-$ that:
 $$
 (\phi_1(x))^{k+1} W(\psi_1,\ldots,\psi_k)\equiv W(\phi_1,\ldots,\phi_{k+1})\ne0\ on\ [T,x_0[\ ;\leqno(7.10)$$
 $$\begin{cases}
 \psi_k\equiv (\phi_1)^{-2} W(\phi_1,\phi_{k+1})\gg (\phi_1)^{-2} W(\phi_1,\phi_{k+2})\equiv \psi_{k+1},\ x\to x_0^-,\  {\rm i.e.} \\
 \psi_1(x)\gg\psi_2(x)\gg\ldots\gg\psi_i(x),\ x\to x_0^-.\end{cases}\leqno(7.11)$$
 
(The $n$-tuple $(\psi_1,\ldots,\psi_n)$ is sometimes called the ``reduced system''.) Moreover (7.6) and (7.10) imply:
 $$
 \frac{W( \psi_1(x), \dots ,\psi_{i-2}(x),\psi_i(x))}{W( \psi_1(x), \dots ,\psi_{i-2}(x),\psi_{i-1}(x))}=\int_x^{x_0}w_i(t)dt.\leqno(7.12)$$
 
 We may now apply our inductive hypothesis inferring that:
 $$
 \psi_i(x)=\widetilde {w}_0(x)\cdot \int_x^{x_0} \widetilde {w}_1 \dots \int _{t_{i-2}}^{x_0} \widetilde {w}_{i-1},\leqno(7.13)$$
 where the $\widetilde {w}_k$'s are defined by the expressions on the right of (2.30) with the $\phi_k$'s replaced by the $\psi_k$'s and (7.10) implies:
 $$\begin{cases}
 \widetilde {w}_0:=\psi_1\equiv (\phi_2/\phi_1)'=-w_1,\\ \\ \displaystyle
  \widetilde {w}_k:=- \frac{W( \psi_1, \dots ,\psi_{k-1},\psi_{k+1})}{W( \psi_1, \dots ,\psi_{k-1},\psi_k)}=-\frac{W( \phi_1, \dots ,\phi_k,\phi_{k+2})}{W( \phi_1, \dots ,\phi_k,\psi_{k+1})}=w_{k+1},\ 2\le k\le i-1, \end{cases}\leqno(7.14)$$
and (7.13) becomes:
 $$
 \big(\phi_{i+1}/\phi_1)'(x)=-w_1(x)\cdot \int_x^{x_0} w_2 \dots \int _{t_{i-2}}^{x_0} w_i, \leqno(7.15)$$ 
 which, by (2.23), gives the sought-for formula for $\phi_{i+1}$.  Formulas (2.32) may be proved quite simply, in alternative to the inductive argument suggested in [8; p. 380], using  the second expressions for the $w_i$'s given in (2.31); putting for brevity:
 $$
 W_i:=W( \phi_1, \dots ,\phi_i),\leqno(7.16)$$
 we have as in [1; p. 92]:
 $$\begin{cases}
 W_1=w_0,\ \ W_2/W_1=-w_0w_1,\\ W_{i+1}/W_i=-w_i\cdot W_i/W_{i-1}=+w_iw_{i-1}W_{i-1}/W_{i-2}=\ldots=\\
 =(-1)^iw_0w_1\ldots w_i,\ 2\le i\le n-1;\end{cases}\leqno(7.17)$$
 hence:
 $$\begin{cases}
 W_i=(-1)^{i-1}w_0w_1\ldots w_{i-1}\cdot W_{i-1}=\\ =(-1)^{i-1}w_0w_1\ldots w_{i-1}\cdot (-1)^{i-2}w_0w_1\ldots w_{i-2}\cdot W_{i-2}=\ldots =\\  
 =(-1)^{(i-1)+(i-2)+\ldots+2+1}[w_0w_1\ldots w_{i-1}][w_0w_1\ldots w_{i-2}]\ldots[w_0w_1]w_0=\ (2.32),
\end{cases}\leqno(7.18)$$
 and this shows the converse inference (iv)$\Rightarrow$(ii). \hfill $\Box$

 \vskip5pt
{\it Proof of Proposition} 2.4.\ (i)-(ii). Properties in (2.36) follow directly from the assumptions and relations in (2.37) are a standard fact as remarked in the preceding proof. As concerns (2.38) the continuity of the $q_i$'s at the endpoint $T$ implies $\int_T(1/q_i)<+\infty$ whereas from (2.37) we get:
$$ \begin{cases}\displaystyle
    \int_T^x1/q_1=constant+\frac{\phi_2(x)}{\phi_1(x)}\ \overset{\ by\ (2.23)\ }
    {\ \ convergent\ \ },\ x\to x_0^-;  \\ \\ \displaystyle
    \int_T^x1/q_i(x)=constant+ \frac{W( \phi_1(x),\dots, \phi_{i-1}(x) ,\phi_{i+1}(x))}
    {W( \phi_1(x), \dots ,\phi_{i-1}(x),\phi_{i}(x))}\overset{by\ (2.14)}
    {\ \ convergent\ \ },\\ as\ x\to x_0^- ,\ 2\le i\le n-1.\end{cases}  \leqno (7.19)$$

Factorization (2.39) is then the classical factorization arising from (2.35) and discovered for the first time by P\'olya [12].  Representations (2.40) are contained in Proposition 2.3 with different notations. In general, by (2.12), the calculations in (7.19) prove the existence of a C.F. of type (II) at $x_0$ valid on a suitable left neighborhood of $x_0$.

(iii). The very same reasonings prove the properties of the $p_i$'s; the proof of (2.45) is similar to that in (7.19):
$$ \begin{cases}\displaystyle
    \int_T^x1/p_1=constant+\frac{\phi_{n-1}(x)}{\phi_n(x)}\overset{\ by\ (2.23)\ }
    {\ \ divergent\ \ },\ x\to x_0^-;  \\ \\ \displaystyle
    \int_T^x1/p_i(x)=constant+ \frac{W( \phi_n(x), \dots ,\phi_{n-i+2}(x),\phi_{n-i}(x))}{W( \phi_n(x), \dots ,\phi_{n-i+2}(x),\phi_{n-i+1}(x))}\overset{by\ (2.14)}
    {\ \ divergent\ \ },\\ as\ x\to x_0^- ,\ 2\le i\le n-1;\end{cases}  \leqno (7.20)$$
and in general, by (2.10), these calculations  prove the existence of a C.F. of type (I) at $x_0$ valid on the whole open interval where the given operator is asssumed disconjugate. The claims in (iv) are trivial.    \hfill $\Box$
    
    \vskip5pt
{\it Proof of Proposition} 3.2.\  Relations (3.35) to (3.37)  are directly checked using representations (2.40). Relation (3.38) follows from the second relation in (3.29) replacing $u$ by $\phi_{k+1}$ and using (3.36). If (3.39) holds true for some sufficiently regular $f$ then (3.27) implies $M_k[f(x)]\equiv0$ and (3.40) follows from (3.35)-(3.36). The converse trivially follows again from (3.35)-(3.36). Now suppose (3.39)-(3.40) to be true on the left of $x_0$; relation (3.41) for $h=1$ is nothing but the obvious relation
 $a_1=\lim_{x\to x_0^-}f(x)/\phi_1(x).$
 
For $h\geq 2$ we use all relations (3.35), (3.36), (3.37) and get from (3.40):
 $$\begin{cases}\displaystyle
 M_{h-1}[f(x)]=\sum_{i=0}^{k-h}a_{h+i}M_{h-1}[\phi_{h+i}(x)]= \\  \displaystyle
 =\epsilon_{h-1}a_h+\sum_{i=1}^{k-h}a_{h+i}M_{h-1}[\phi_{h+i}(x)]= \epsilon_{h-1}a_h+o(1),
 \end{cases}\leqno (7.21)$$
 where the remainder ``$o(1)"$ is $\equiv 0$ for $h=k$.
\hfill $\Box$
       
       \vskip5pt
 {\it Proof of Lemma} 4.2.\ From the chain $P_{n-1}(x)\gg\ldots\gg P_0(x),\ 
 x\to x_0$, we get
      $$
      \phi(x)=cP_i(x)+\alpha_{i-1}P_{i-1}(x)+\ldots +\alpha_0P_0(x) \leqno (7.22)$$
for suitable constants $\alpha_k$, hence
 $$
    L_k[\phi(x)]=c L_k[P_i(x)]+\alpha_{i-1}L_k[P_{i-1}(x)]+\ldots +\alpha_0L_k[P_0(x)];
     \leqno (7.23)$$
now (4.10) follows from (4.7), and (4.11) follows from (4.4). If in (7.22) we replace $\phi$ by $\phi_{n-i}$ we have $c=b_{n-i}$ and the identities in (4.12) follow from (4.4) and (4.5). The identity in (4.12) for $i=k$,  i.e. $L_k[\phi_{n-k}(x)]\equiv b_{n-k}$, together with the first relation in (3.28) imply (4.14).  \hfill{$\Box$}

\vskip5pt
{\it Proof of  Theorem} 4.4.\ \  Part (I). From (4.12) and (4.17), with $c_1=a_1$, we infer at once the equivalence ``(ii)$\Leftrightarrow$(iii)" as well as representation in (4.21). The inference ``(i)$\Rightarrow$(ii)" being obvious let us prove the converse simply denoting by $L$ our operator $L_{\phi_1,\dots,\phi_n}$. We shall repeatedly use the recursive formulas (3.25) in the form
$$
L_{k-1}u=\int_T^x\frac{1}{p_k(t)}L_k[u(t)] dt + \textit{constant},\ 1\leq k\leq n. 
\leqno (7.24) $$ 

If (4.19) holds true we have (4.21), and representations in (4.16) can be rewritten as
$$\begin{cases}\displaystyle
L_k[f(x)] =\big\{c_1L_k[\phi_1(x)]+\dots +c_{n-k}L_k[\phi_{n-k}(x)]\big\}+ 
\\ \\ \displaystyle
+ \left(\int_T^{x_0} \frac{L[f(t)]}{p_n(t)} dt\right)\cdot \int_T^x \frac{dt_{k+1}}{p_{k+1}(t_{k+1})} \dots \int_T^{t_{n-2}} \frac{dt_{n-1}}{p_{n-1}(t_{n-1})}\ +\\ \\ \displaystyle
   - \int_T^x \frac{dt_{k+1}}{p_{k+1}(t_{k+1})} \dots \int_T^{t_{n-2}} \frac{dt_{n-1}}{p_{n-1}(t_{n-1})}\int_{t_{n-1}}^{x_0} \frac{L[f(t)]}{p_n(t)} dt;\  \ 0\leq k\leq n-2. 
   \end{cases}\leqno (7.25)$$

 Now we have
   $$
   \int_T^x \frac{dt_{k+1}}{p_{k+1}(t_{k+1})} \dots \int_T^{t_{n-2}} \frac{dt_{n-1}}{p_{n-1}(t_{n-1})}\overset{\textit{by}\ (4.6)}{\equiv} L_k[P_{n-1}(x)]=\cdots 
   \leqno (7.26) $$
   by (4.1) and (4.10) \hskip20pt
 $ \cdots =(1/b_1)L_k[\phi_1(x)]+o\big(L_k[\phi_1(x)]\big).$
 
 After substituting into (7.25) we get: 
 $$\begin{cases}\displaystyle
 L_k[f(x)] =c_1L_k[\phi_1(x)]+o\left(L_k[\phi_1(x)]\right)+\bar{c}L_k[P_{n-1}(x)]+
 o\big(L_k[P_{n-1}(x)]\big)=\\ \\ \displaystyle
 =\left(c_1+\frac{\bar{c}}{b_1}\right)L_k[\phi_1(x)]+o\big(L_k[\phi_1(x)]\big),\ \ 0\leq k\leq n-2,\end{cases}\leqno (7.27) $$
 where $\bar{c}:=\int_T^{x_0}L[f(t)]/p_n(t)\ dt$, and the coefficient
 $$
 c:=c_1+(\bar{c}/b_1) \leqno (7.28) $$
 is independent of $k$. To show that $c$ coincides with the $a_1$ appearing in (4.19) we may suitably integrate (4.19) to obtain, by (3.24),
 $$\begin{cases}\displaystyle
 L_{n-2}[f(x)] =\int_T^x\frac{L_{n-1}[f(t)]}{p_{n-1}(t)}dt + \textit{constant} \overset{\textit{by}\ (4.19)}{=} \\ \\ \displaystyle
 =a_1b_1\int_T^x\frac{dt}{p_{n-1}(t)}+o\left(\int_T^x\frac{dt}{p_{n-1}(t)}\right)
 \overset{\textit{by}\ (4.6)}
 \equiv a_1b_1L_{n-2}[P_{n-1}(x)]+o\big(L_{n-2}[P_{n-1}(x)] \big)=\dots
 \end{cases}\leqno(7.29) $$
 by (4.1) and (4.10)\hskip20pt
 $ \dots =a_1L_{n-2}[\phi_1(x)]+o\left(L_{n-2}[\phi_1(x)]\right). $
 
 \vskip5pt
 Part (II). Case $i=2$. We must prove the equivalence of the following three contingencies:
 $$\begin{cases}
 f(x)=a_1\phi_1(x)+a_2\phi_2(x)+o(\phi_2(x)), \\
 L_k[f(x)]=a_1L_k[\phi_1(x)]+a_2L_k[\phi_2(x)]+o\big(L_k[\phi_2(x)]\big),\ 1\leq k\leq n-2, \\
  L_{n-1}[f(x)]=a_1L_{n-1}[\phi_1(x)]+o\big(L_{n-1}[\phi_1(x)]\big);
  \end{cases} \leqno(7.30) $$
  $$\begin{cases}
   L_{n-2}[f(x)]=a_1L_{n-2}[\phi_1(x)]+a_2L_{n-2}[\phi_2(x)]+o\big(L_{n-2}[\phi_2(x)]\big), \\
  L_{n-1}[f(x)]=a_1L_{n-1}[\phi_1(x)]+o\big(L_{n-1}[\phi_1(x)]\big);
  \end{cases} \leqno(7.31) $$
  $$
  \int_T^{x_0} \frac{dt}{p_{n-1}(t)}\int_t^{x_0} \frac{L[f(\tau)]}{p_n(\tau)} d\tau \ \ 
  \textit{convergent.} \leqno(7.32) $$
  
  First we prove ``(7.31)$\Leftrightarrow$(7.32)". If (7.32) holds true then , by part (I) of our theorem, we have all relations in (4.18) and in particular the second relation in (7.31). Moreover we can rewrite representation in (4.16) for $k=n-2$ in the form:
  $$
   L_{n-2}[f(x)]=a_1L_{n-2}[\phi_1(x)]+a_2L_{n-2}[\phi_2(x)]
   +\int_x^{x_0} \frac{dt}{p_{n-1}(t)}\int_t^{x_0} \frac{L[f(\tau)]}{p_n(\tau)} d\tau, 
   \leqno (7.33) $$
   where $a_1$ is just the same as in the second relation in (7.31) and $a_2$ is a suitable constant. This yields the first relation in (7.31) because $L_{n-2}[\phi_2(x)]$ is a nonzero constant by (4.12). Viceversa if relations in (7.31) hold true then, by part (I), we have representation in (4.21) by which  we replace the quantity $L_{n-1}[f(t)]$ in the first equality in (7.33). Denoting by $c_{n-2}, \bar{c}_{n-2}$ suitable constants we get:
   $$\begin{cases}\displaystyle
   L_{n-2}[f(x)] =c_{n-2}+\int_T^x\frac{L_{n-1}[f(t)]}{p_{n-1}(t)}dt = \\ \\ \displaystyle
    =c_{n-2}+\int_T^x\frac{1}{p_{n-1}(t)}\left[a_1L_{n-1}[\phi_1(t)]-
    \int_t^{x_0} \frac{L[f(\tau)]}{p_n(\tau)} d\tau\right] dt= \\ \\ \displaystyle
    =\bar{c}_{n-2}+a_1L_{n-2}[\phi_1(x)]-\int_T^x\frac{dt}{p_{n-1}(t)}\int_t^{x_0} \frac{L[f(\tau)]}{p_n(\tau)} d\tau. \end{cases}\leqno (7.34) $$
    
    By comparison with the first relation in (7.31) we get (7.32) because $L_{n-2}[\phi_2(x)]$ is a  constant. As the inference ``(7.30)$\Rightarrow$(7.31)" is obvious it remains to prove the converse. Using (7.24) and integrating the first relation in (7.31) we get (with suitable constants $c_{n-3}, \bar{c}_{n-3}$):
    $$\begin{cases}\displaystyle
    L_{n-3}[f(x)] \equiv c_{n-3}+\int_T^x\frac{L_{n-2}[f(t)]}{p_{n-2}(t)}dt =c_{n-3}+ \\ \\ \displaystyle
 +a_1\int_T^x\frac{L_{n-2}[\phi_1(t)]}{p_{n-2}(t)}dt+a_2\int_T^x\frac{L_{n-2}[\phi_2(t)]}{p_{n-2}(t)}dt 
 +\int_T^xo\left(\frac{L_{n-2}[\phi_2(t)]}{p_{n-2}(t)}\right)dt=\dots
 \end{cases}\leqno (7.35) $$  
 as $L_{n-2}[\phi_2(x)]$ is a nonzero constant and $\int^{\to x_0}1/p_{n-2}$ diverges
 $$
 \dots =\bar{c}_{n-3}+a_1L_{n-3}[\phi_1(x)]+a_2L_{n-3}[\phi_2(x)]+o\big(L_{n-3}[\phi_2(x)]\big). $$
 
 Here the constant  $\bar{c}_{n-3}$  is meaningless as the comparison functions are divergent as $x\to x_0^-$. Iterating the procedure we get all relations in (7.30). By induction on $i$ and the same kind of reasonings our theorem is proved for each $i\leq n$. \hfill $\Box$
 
 \vskip5pt
 {\it Proof of Theorem} 5.1.\  (i)$\Rightarrow$(ii). Relation (5.5) implies the existence of $a_1\equiv \lim_{x\to x_0^-}f(x)/\phi_1(x)\equiv \lim_{x\to x_0^-}M_0[f(x)]$, and each relation in (5.6) implies the relation in (5.7) with the same value of $k$ because of (3.36)-(3.37). (ii)$\Rightarrow$(iii) is obvious. (iii)$\Leftrightarrow$(iv). It follows from (5.3) that the limit in (5.8) exists in $\mathbb{R}$  iff (5.9) holds true and, in this case, (5.3) can be written as
 $$
 M_{n-1}[f(x)]=a_n-\int_x^{x_0} \frac{L[f(t)]}{q_n(t)} dt , \leqno(7.36) $$
 where, as above, $L\equiv L_{\phi_1,\dots,\phi_n}$. 
 
 (iv)$\Rightarrow$(i). We have already proved (7.36) which is (5.6) for $k=n-1$ together with an integral representation of the remainder. For $k=n-2$ the recursive formulas (3.27) give
 $$
 (M_{n-2}[f(x)])'=\frac{1}{q_{n-1}(x)}M_{n-1}[f(x)], \leqno (7.37) $$
 whence, by (7.36) and (2.38), we get:
 $$
 M_{n-2}[f(x)]=a_{n-1}-a_n\int_x^{x_0}\frac{1}{q_{n-1}}+\int_x^{x_0}\frac{d\tau}{q_{n-1}(\tau)}\int_{\tau}^{x_0} \frac{L[f(t)]}{q_n(t)} dt,  \leqno(7.38) $$ 
 for a suitable constant $a_{n-1}$. By  (3.36)-(3.37) this is nothing but
 $$
 M_{n-2}[f(x)]=a_{n-1}M_{n-2}[\phi_{n-1}(x)]+a_nM_{n-2}[\phi_n(x)]
 +\int_x^{x_0}\frac{d\tau}{q_{n-1}(\tau)}\int_{\tau}^{x_0} \frac{L[f(t)]}{q_n(t)} dt ,
 \leqno (7.39) $$
 which is the relation in (5.6) for $k=n-2$ with a representation of the remainder. In a similar way 
 for $k=n-3$ we start from
 $$
  (M_{n-3}[f(x)])'=\frac{1}{q_{n-2}(x)}M_{n-2}[f(x)], \leqno (7.40) $$
 and integrate (7.38) after dividing by $1/q_{n-2}$, so getting
 $$\begin{cases}\displaystyle
 M_{n-3}[f(x)]=a_{n-2}-a_{n-1}\int_x^{x_0}\frac{1}{q_{n-2}}+a_n\int_x^{x_0}\frac{dt_{n-2}}{q_{n-2}}\int_{t_{n-1}}^{x_0}\frac{dt_{n-2}}{q_{n-1}}+\\ \\ \displaystyle
 -\int_x^{x_0}\frac{dt_{n-2}}{q_{n-2}}\int_{t_{n-2}}^{x_0}\frac{dt_{n-1}}{q_{n-1}}\int_{t_{n-1}}^{x_0}\frac{L[f(t)]}{q_n(t)}dt \end{cases} \leqno (7.41)$$
 for a suitable constant $a_{n-2}$. By  (3.36)-(3.37) this can be rewritten as
 $$\begin{cases}\displaystyle
 M_{n-3}[f(x)]=a_{n-2}M_{n-3}[\phi_{n-2}(x)]+a_{n-1}M_{n-3}[\phi_{n-1}(x)]+
 \\ \\ \displaystyle
 +a_nM_{n-3}[\phi_n(x)]-\int_x^{x_0}\frac{1}{q_{n-2}}\int^{x_0}\frac{1}{q_{n-1}}\int^{x_0}\frac{L[f(t)]}{q_n(t)}dt ,\end{cases}  \leqno (7.42)$$
 with a suitable constant $a_{n-2}$. An iteration of the procedure gives all relations in  (5.6) together with the representation formulas (5.11)-(5.12) for $R_0(x)$ and (5.15) for $R_k(x), k\geq 1.$: ``(i)$\Rightarrow$(v)" has been proved. The last inference ``(v)$\Rightarrow$(i)" and (5.12) are trivially proved by applying the operators $M_k$ to (5.11). \hfill $\Box$
 
 \vskip5pt
  {\it Proof of Theorem} 5.2.\  (i)$\Rightarrow$(ii) follows from (3.36)-(3.37). (ii)$\Rightarrow$(iii) is obvious. (iii)$\Leftrightarrow$(iv): by (3.36)-(3.37)     	  the representation in (5.2) for $k=i-1$ has the form
  $$
  M_{i-1}[f(x)]=c_i+\int_T^x\frac{1}{q_i}\int_T^{t_i}\frac{1}{q_{i+1}}\dots \int_T^{t_{n-2}} \frac{1}{q_{n-1}}\int_T^{t_{n-1}}\frac{L[f(t)]}{q_n(t)}dt+o(1),\leqno(7.43)$$
  whence our equivalence follows at once. If this is the case (5.2) can be rewritten as
  $$\begin{cases}\displaystyle
   M_{i-1}[f(x)]=a_i-\int_x^{x_0} \frac{1}{q_i} \int_T^{t_i} \frac{1}{q_{i+1}} \dots \int_T^{t_{n-1}} \frac{L[f(t)]}{q_n(t)} dt+\\ \\ \displaystyle
   +c_{i+1}M_{i-1}[\phi_{i+1}(x)]+\ldots +c_nM_{i-1}[\phi_n(x)],\ x\in [T, x_0[,
   \end{cases} \leqno(7.44)$$
   where $a_i$ is uniquely determined by (5.23) but $c_{i+1},\dots ,c_n$ are non-better specified constants not determinable by the sole condition (5.23).
   
   (iv)$\Rightarrow$(i). This is proved like the corresponding inference in Theorem 5.1 by successive integrations of (7.44) starting from
  $$
  (M_{i-2}[f(x)])'=\frac{1}{q_{i-1}(x)}M_{i-1}[f(x)], \leqno (7.45) $$ 
 whence, by (2.38),(3.27) and (3.37), we get:
 $$\begin{cases}\displaystyle
  M_{i-2}[f(x)]=a_{i-1}-a_i\int_x^{x_0}\frac{1}{q_{i-1}}+\int_x^{x_0}\frac{1}{q_{i-1}}\int_{t_{i-1}}^{x_0} \frac{1}{q_i} \int_T^{t_i} \frac{1}{q_{i+1}} \dots \int_T^{t_{n-1}} \frac{L[f(t)]}{q_n(t)} dt+\\ \\ \displaystyle
   +c_{i+1}M_{i-2}[\phi_{i+1}(x)]+\ldots +c_nM_{i-2}[\phi_n(x)]=\ \ \ {\rm by\ (3.31)
   \  and\ (3.37)}\ \ \ =\\  \displaystyle
   =a_{i-1}+a_iM_{i-2}[\phi_i(x)]+o(M_{i-2}[\phi_i(x)]),\end{cases} \leqno(7.46)$$
   where the constant $a_{i-1}$, which includes all the constants from integration of the various terms, is uniquely determined by (5.22). By iteration of the procedure we get all relations in (5.20)-(15.21).  Relation (5.28) easily follows from (7.44) by (3.36)-(3.37). \hfill $\Box$
 
 \vskip5pt
  {\it Proof of Theorem} 5.3.\  This is almost a word-for word repetition of the proofs of Theorems 5.1-5.2. (i)$\Rightarrow$(ii). For $0\leq k\leq i-2$ this is included in the same inference in Theorems 5.1-5.2; whereas the relation in (5.30) for $k=i-1$ just reads $M_{i-1}[f(x)]=O(1)$. (ii)$\Rightarrow$(iii) is obvious.  (iii)$\Leftrightarrow$(iv) follows from (7.43). To show (iv)$\Rightarrow$(i) we use (7.45) and the representation in (5.2) for $k=i-1$ instead of (7.44) as in the proof of  Theorem 5.2. Due to the convergence of $\int^{\to x_0}1/q_{i-1}$\ we may still apply the operator $\int_x^{x_0}$ so getting, instead of (7.46),
  $$\begin{cases}\displaystyle
   M_{i-2}[f(x)]=a_{i-1}-\int_x^{x_0}\frac{1}{q_{i-1}}\int_T^{t_{i-1}}\frac{1}{q_i} \int_T^{t_i} \frac{1}{q_{i+1}} \dots \int_T^{t_{n-1}} \frac{L[f(t)]}{q_n(t)} dt+ 
   \\ \\ \displaystyle
   -\int_x^{x_0}\frac{1}{q_{i-1}}\left[\sum_{j=i+1}^nc_jM_{i-2}[\phi_j(x)]\right]=a_{i-1}+O\left(\int_x^{x_0}\frac{1}{q_{i-1}}\right)+o\left(\int_x^{x_0}\frac{1}{q_{i-1}}\right)=\\ \\ \displaystyle
   =a_{i-1}+O\left(\int_x^{x_0}\frac{1}{q_{i-1}}\right)\equiv a_{i-1}+O\big(M_{i-2}[\phi_i(x)]\big).\end{cases}\leqno(7.47)$$
   
   By iteration we get all relations in (5.28)-(5.29). \hfill $\Box$
  
  \vskip5pt
 {\it Proof of Theorem}  6.1.\  The only thing to be proved is the inference
 ``(i)$\Rightarrow$(v)$\land$(vi)",  the other properties being included in Theorems 4.5 and 5.1. We use a procedure already used in [4; p. 193] and in [6; p. 213]. From representation in (4.15) we get (using the simplified notation $L\equiv L_{\phi_1,\dots,\phi_n}$):
 $$
   \frac{f(x)}{\phi_1(x)}-c_1+o(1)
     =\frac{ 1/p_0(x) }{\phi_1(x)} \int_T^x \frac{1}{p_1} \dots \int_T^{t_{n-2}} \frac{1}{p_{n-1}}\int_T ^{t_{n-1}}\frac{L[f(t)]}{p_n(t)} dt, \ x\in [T,x_0[.   \leqno (7.48) $$  
    
     By the assumption (6.2) the left-hand side has a finite limit as $x\to x_0$, and for the limit of the right-hand side we have:
     $$
     \lim_{x\to x_0^-}\frac{ 1/p_0(x) }{\phi_1(x)} \int_T^x \frac{1}{p_1} \dots \int_T^{t_{n-2}} \frac{1}{p_{n-1}}\int_T ^{t_{n-1}}\frac{L[f(t)]}{p_n(t)} dt \overset{\ by\ (4.1)\ and\ (2.42)\ }{\ =\ } \leqno(7.49)$$
     $$
    = \frac{1}{b_1} \lim_{x\to x_0^-}\frac{\displaystyle\int_T^x \frac{1}{p_1} \dots \int_T^{t_{n-2}} \frac{1}{p_{n-1}}\int_T ^{t_{n-1}}L[f(t)]/p_n(t)\ dt }{P_{n-1}(x)/P_0(x)}=$$
     $$
     =\frac{1}{b_1} \lim_{x\to x_0^-}\frac{\displaystyle\int_T^x \frac{1}{p_1} \dots \int_T^{t_{n-2}} \frac{1}{p_{n-1}}\int_T ^{t_{n-1}}L[f(t)]/p_n(t)\ dt }{\displaystyle\int_T^x \frac{1}{p_1} \dots \int_T^{t_{n-2}} \frac{1}{p_{n-1}}}=\dots =\frac{1}{b_1} \lim_{x\to x_0^-}\int_T^xL/p_n,$$
   after applying L'Hospital's rule $(n-1)$ times (which is legitimate as all the denominators diverge to $+\infty$).   
     By the positivity of the integrand this last limit exists in $\overline{\mathbb{R}}$ and coincides with the limit of the left-hand side in (7.48) hence it must be a real number and (4.15) can take the form:
     $$\begin{cases}\displaystyle
      f(x)=a_1\phi_1(x)+c_2\phi_2(x)+\dots +c_n\phi_n(x)+\\ \\ \displaystyle
     - \frac{1}{p_0(x)}  \int_T^x \frac{1}{p_1} \dots \int_T^{t_{n-2}} \frac{1}{p_{n-1}}
     \int_{t_{n-1}}^{x_0}
     \frac{L[f(t)]}{p_n(t)}  dt, \ x\in [T,x_0[\ ,\end{cases}   \leqno (7.50) $$  
  with suitable constants $c_2,\dots ,c_n$. From this we get:
$$\begin{cases}\displaystyle
\frac{f(x)-a_1\phi_1(x)}{\phi_2(x)}-c_2+o(1)=\\ \\ \displaystyle
 =- \frac{1/p_0(x)}{\phi_2(x)}  \int_T^x \frac{1}{p_1} \dots \int_T^{t_{n-2}} \frac{1}{p_{n-1}} \int_{t_{n-1}}^{x_0}
     \frac{L[f(t)]}{p_n(t)}  dt, \ x\in [T,x_0[\ .\end{cases}  \leqno (7.51) $$

Here again the left-hand side has a finite limit as $x\to x_0$ whereas  the limit of the right-hand side, by (4.1), equals:
$$
 -\frac{1}{b_2} \lim_{x\to x_0^-}\frac{\displaystyle\int_T^x \frac{1}{p_1} \dots \int_T^{t_{n-2}} \frac{1}{p_{n-1}}\int_{t_{n-1}}^{x_0}L[f(t)]/p_n(t)\ dt }{P_{n-2}(x)/P_0(x)}=$$
 $$
     =-\frac{1}{b_2} \lim_{x\to x_0^-}\frac{\displaystyle\int_T^x \frac{1}{p_1} \dots \int_T^{t_{n-2}} \frac{1}{p_{n-1}}\int_{t_{n-1}}^{x_0}L[f(t)]/p_n(t)\ dt }{\displaystyle\int_T^x \frac{1}{p_1} \dots \int_T^{t_{n-3}} \frac{1}{p_{n-2}}}=\dots =$$
$$
      =-\frac{1}{b_2} \lim_{x\to x_0^-}\int_T^x\frac{1}{p_{n-1}}\int_{t_{n-1}}^{x_0}L/p_n,$$
 after applying L'Hospital's rule $(n-2)$ times.     
     Hence this last limit, which exists in $\overline{\mathbb{R}}$, must be a real number and (7.50) can be rewritten as:
     $$\begin{cases}
       f(x)=a_1\phi_1(x)+a_2\phi_2(x)+c_3\phi_3(x)+\dots +c_n\phi_n(x)+ 
       \\ \\ \displaystyle
          + \frac{1}{p_0(x)}  \int_T^x \frac{1}{p_1} \dots \int_T^{t_{n-3}} \frac{1}{p_{n-2}}
     \int_{t_{n-2}}^{x_0} \frac{1}{p_{n-1}} \int_{t_{n-1}}^{x_0}
     \frac{L[f(t)]}{p_n(t)}  dt,  x\in [T,x_0[\ , \end{cases}  \leqno (7.52) $$  
with suitable constants $c_3,\dots ,c_n$.  It is now clear how this procedure works and by induction one can prove the validity of representation:
$$\begin{cases}\displaystyle
 f(x)=a_1\phi_1(x)+\ldots +a_{n-1}\phi_{n-1}(x)+c_n\phi_n(x)+\\ \\ \displaystyle
+ \frac{(-1)^{n-1}}{p_0(x)}  \int_T^x \frac{1}{p_1}  \int _x^{x_0} \frac{1}{p_2}\dots \int_{t_{n-3}}^{x_0} \frac{1}{p_{n-2}}
     \int_{t_{n-2}}^{x_0} \frac{1}{p_{n-1}} \int_{t_{n-1}}^{x_0} \frac{L[f(t)]}{p_n(t)}  dt, 
      \ \ x\in [T,x_0[\ , \end{cases}\leqno (7.53) $$
with a suitable constant $c_n$. As a last step we observe that (6.2) implies:
$$
 \big[f(x)-a_1\phi_1(x)-\dots -a_{n-1}\phi_{n-1}(x)\big]\big/\phi_n(x)=O(1),\ x\to x_0^-, 
  \leqno (7.54)  $$
  and (7.53) in turn implies:
  $$\begin{cases}\displaystyle
   \frac{1/p_0(x)}{\phi_n(x)}  \int_T^x \frac{1}{p_1}  \int _x^{x_0} \frac{1}{p_2}\dots      \int_{t_{n-2}}^{x_0} \frac{1}{p_{n-1}} \int_{t_{n-1}}^{x_0} \frac{L[f(t)]}{p_n(t)}  dt
     \overset{\ by (4.2)\ }{\equiv} \\ \\ \displaystyle
     \equiv  \int_T^x \frac{1}{p_1}  \int _x^{x_0} \frac{1}{p_2}\dots      \int_{t_{n-2}}^{x_0} \frac{1}{p_{n-1}} \int_{t_{n-1}}^{x_0} \frac{L[f(t)]}{p_n(t)}  dt=O(1),\ x\to x_0^-.\end{cases}\leqno(7.55)$$
     
By the positivity of the integrand this last relation implies (6.6) and the first representation in (6.8) for $R_0(x)$.  
To prove (6.7) we apply the same ideas starting from representation (5.1) and dividing by $\phi_1$; recalling that $\phi_1=1/q_0$ we get:
$$
    \frac{f(x)}{\phi_1(x)}-c_1+o(1)= 
     \int_T^x \frac{1}{q_1} \dots \int_T^{t_{n-2}} \frac{1}{q_{n-1}}\int_T ^{t_{n-1}}
     \frac{ L[f(t)]}{q_n(t)} dt, \ x\in [T,x_0[\ .  \leqno (7.56) $$  
     
     This implies:
     $$
     \int_T^{x_0} \frac{1}{q_1}\int_T^{t_1} \frac{1}{q_2} \dots \int_T^{t_{n-2}} \frac{1}{q_{n-1}}\int_T ^{t_{n-1}}\frac{L[f(t)]}{q_n(t)}dt<+\infty, \leqno (7.57) $$
and (5.1) can be rewritten as:
$$\begin{cases}\displaystyle
 f(x)=a_1\phi_1(x)+c_2\phi_2(x)+\dots +c_n\phi_n(x)+\\ \\ \displaystyle
     - \frac{1}{q_0(x)}  \int_x^{x_0} \frac{1}{q_1}\int_T^{t_1} \frac{1}{q_2} \dots \int_T^{t_{n-2}} \frac{1}{q_{n-1}}\int_T ^{t_{n-1}} \frac{L[f(t)]}{q_n(t)}dt,\ x\in [T,x_0[\ , 
     \end{cases}\leqno (7.58) $$  
    with suitable constants $c_2,\dots ,c_n$. From this we get:
$$
\frac{f(x)-a_1\phi_1(x)}{\phi_2(x)}-c_2+o(1)
=-\frac{\phi_1(x)}{\phi_2(x)} \int_x^{x_0} \frac{1}{q_1}\int_T^{t_1} \frac{1}{q_2} \dots \int_T ^{t_{n-1}} \frac{L[f(t)]}{q_n(t)}dt. \leqno (7.59) $$

Evaluating the limit of the right-hand side by L'Hospital's rule and using formula in  (2.37), $1/q_1=(\phi_2/\phi_1)'$, we get:
$$
\lim_{x\to x_0^-}\frac{ \displaystyle-\int_x^{x_0} \frac{1}{q_1}\int_T^{t_1} \frac{1}{q_2} \dots \int_T ^{t_{n-1}}L[f(t)]/q_n(t)\ dt} {\phi_2(x)/\phi_1(x)}\overset{\ H}{=}
\lim_{x\to x_0^-}\int_T^x \frac{1}{q_2} \dots \int_T ^{t_{n-1}}L[f(t)]/q_n(t)\ dt ,$$
and  this last limit, which exists in $\overline{\mathbb{R}}$, must be a real number. This means that
$$
 \int_T^{x_0} \frac{1}{q_1}\int_{t_1}^{x_0} \frac{1}{q_2}\int_T^{t_2} \frac{1}{q_3} \dots \int_T^{t_{n-2}} \frac{1}{q_{n-1}}\int_T ^{t_{n-1}}\frac{L[f(t)]}{q_n(t)}dt<+\infty, 
 \leqno (7.60) $$
and (7.50) can be rewritten as:
$$\begin{cases}\displaystyle
  f(x)=a_1\phi_1(x)+a_2\phi_2(x)+c_3\phi_3(x)+\dots +c_n\phi_n(x)+ 
  \\ \\ \displaystyle
+\frac{1}{q_0(x)} \int_x^{x_0} \frac{1}{q_1}\int_{t_1}^{x_0} \frac{1}{q_2}\int_T^{t_2} \frac{1}{q_3} \dots \int_T^{t_{n-2}} \frac{1}{q_{n-1}}\int_T ^{t_{n-1}}\frac{L[f(t)]}{q_n(t)}dt,\end{cases}\leqno (7.61) $$
 with suitable constants $c_3,\dots ,c_n$. For the clarity's sake we make explicit the steps of this second part of our proof. Assume by induction that the following two conditions hold true:
$$
\int_T^{x_0} \frac{1}{q_i}\int_T^{t_i} \frac{1}{q_{i+1}} \dots \int_T^{t_{n-2}} \frac{1}{q_{n-1}}\int_T ^{t_{n-1}}\frac{L[f(t)]}{q_n(t)}dt<+\infty; \leqno (7.62) $$
$$\begin{cases}\displaystyle
 f(x)=a_1\phi_1(x)+\dots +a_i\phi_i(x)+c_{i+1}\phi_{i+1}(x)+\dots +c_n\phi_n(x)+ 
 \\ \\ \displaystyle
+\frac{(-1)^i}{q_0(x)} \int_x^{x_0} \frac{1}{q_1}\int_{t_1}^{x_0}\frac{1}{q_2}\dots \int_{t_{i-1}}^{x_0} \frac{1}{q_i}\int_T^{t_i} \frac{1}{q_{i+1}} \dots \int_T ^{t_{n-1}}\frac{L[f(t)]}{q_n(t)}dt, \end{cases}\leqno (7.63) $$  
  for some $i,1\leq i\leq n-2$. and suitable constants $c_{i+1},\dots ,c_n$. 
 Dividing both sides of (7.63) by $\phi_{i+1}$ and taking account of (6.2) we infer that the limit of the quantity
 $$
 \left[\frac{(-1)^i}{q_0(x)} \int_x^{x_0} \frac{1}{q_1}\dots \int_{t_{i-1}}^{x_0} \frac{1}{q_i}\int_T^{t_i} \frac{1}{q_{i+1}} \dots \int_T ^{t_{n-1}}\frac{L[f(t)]}{q_n(t)}dt\right] \Big/\phi_{i+1}(x) \overset{\ (2.34)\ }{\equiv} \leqno(7.64)$$
 
 $$
 \equiv \left[\int_x^{x_0} \frac{1}{q_1}\dots \int_{t_{i-1}}^{x_0} \frac{1}{q_i}\int_T^{t_i} \frac{1}{q_{i+1}} \dots \int_T ^{t_{n-1}}\frac{L[f(t)]}{q_n(t)}dt\right]\left/\int_x^{x_0} \frac{1}{q_1}\dots \int_{t_{i-1}}^{x_0} \frac{1}{q_i}\right. $$
 exists in $\mathbb{R}$. Applying L'Hospital's rule $i$ times to evaluate this  limit we get the new limit
 $$
 \lim_{x\to x_0^-}\int_T^x \frac{1}{q_{i+1}} \dots \int_T ^{t_{n-1}}\frac{L[f(t)]}{q_n(t)}dt $$
 which, by the positivity of the integrand, exists in $\overline{\mathbb{R}}$ hence it must be a real number. We infer that condition (7.62) holds true with $i$ replaced by $i+1$ and this implies representation (7.63) with $i$ replaced by $i+1$ and suitable constants $c_{i+2},\dots ,c_n$. By this inductive procedure we arrive at representation:
 $$\begin{cases}\displaystyle
 f(x)=a_1\phi_1(x)+\dots +a_{n-1}\phi_{n-1}(x)+c_n\phi_n(x)+\\ \\ \displaystyle
+\frac{(-1)^{n-1}}{q_0(x)} \int_x^{x_0} \frac{1}{q_1}\dots \int_{t_{n-2}}^{x_0} \frac{1}{q_{n-1}}\int_T ^{t_{n-1}}\frac{L[f(t)]}{q_n(t)}dt,\ x\in [T,x_0[\ , \end{cases} \leqno (7.65) $$
with some constant $c_n$. Dividing by $\phi_n$ and using (6.2) we may now conclude that
$$
\left[\int_x^{x_0} \frac{1}{q_1}\dots \int_{t_{n-2}}^{x_0} \frac{1}{q_{n-1}} \int_T ^{t_{n-1}}\frac{L[f(t)]}{q_n(t)}dt\right]\left/\int_x^{x_0} \frac{1}{q_1}\dots \int_{t_{n-2}}^{x_0} \frac{1}{q_{n-1}}\right. =O(1),\leqno (7.66) $$
and if we try to evaluate the limit of the ratio on the left applying L'Hospital's rule $(n-1)$ times we get the $\lim_{x\to x_0^-}\int_T^x L[f(t)]/q_n(t)dt$, 
 which exists in $\overline{\mathbb{R}}$ and must be a finite number. This is condition (6.7) which allows the second representation in (6.8).
\hfill $\Box$

 \vskip5pt
 {\it Proof of Theorem} 6.2.\ \ The equivalence between (6.11) and (6.12) easily follows from Fubini's theorem by interchanging the order of integrations in (6.11) whereas the equivalence between (6.12) and (6.13) is by no means an obvious fact. We give a concise proof based on Theorem 6.1. Putting
 $$
 F(x):= \frac{1}{p_0(x)}  \int_T^x \frac{1}{p_1} \dots \int_T^{t_{n-2}} \frac{1}{p_{n-1}}\int_T ^{t_{n-1}}\frac{1}{p_n(t)}|L_{\phi_1,\dots,\phi_n}[f(t)]|dt, \leqno (7.67)$$
 we have
 $$
 F\in AC^{n-1}[T,x_0[;\  L_{\phi_1,\dots,\phi_n}[F(x)]=|L_{\phi_1,\dots,\phi_n}[f(x)]|\ a.\ e.\ on\ [T,x_0[\ ;\leqno(7.68)$$
 hence $F$ satisfies  $L_{\phi_1,\dots,\phi_n}[F(x)]\geq 0$\ a. e. on $[T,x_0[$ and Theorem 6.1 implies the equivalence between (6.12) and (6.13). \hfill $\Box$
 
  \vskip5pt
 {\it Proof of Theorem} 6.3.\ \ The only thing to prove is the $O$-estimates in (6.18). From representation (7.44) we get:
 $$
    M_i[f(x)]=\int_T^x \frac{1}{q_{i+1}} \dots \int_T^{t_{n-1}}\frac{L[f(t)]}{q_n(t)} dt+ 
   c_{i+1}M_i[\phi_{i+1}(x)]+\dots +c_nM_i[\phi_n(x)]=\cdots\leqno(7.69)$$
   by (3.36) and (3.37)
$$
\dots=\int_T^x \frac{1}{q_{i+1}} \dots \int_T^{t_{n-1}}\frac{L[f(t)]}{q_n(t)} dt+ c+o(1)$$
for some constant $c$ whence it follows the estimate for $M_i[f(x)]$. The other estimates are similarly obtained.  \hfill $\Box$
        
 \vspace{20pt}
  \centerline {\bf 8. \ Appendix: algorithms for constructing canonical factorizations}
    
   \vspace{10pt}
 The original  procedure used by Trench [14] to construct a C.F. of type (I) for a disconjugate operator is not an intuitive one.  Here we exhibit two easier algorithms to construct both types of C.F.'s starting from an explicit fundamental system of solutions which is also an asymptotic scale at one endpoint. The so-obtained factorizations will be proved to coincide with those obtainable by P\'olya's procedure when applied either to the asymptotic scale $(\phi_1,\dots,\phi_n)$ or to the inverted $n$-tuple $(\phi_n,\dots,\phi_1)$. As each step in the algorithms  has an asymptotic meaning they provide asymptotic intepretations of P\'olya's procedure and they may sometimes be quicker to apply than P\'olya's procedure, especially for small values of $n$, avoiding the explicit use of Wronskians. The algorithm for a C.F. of type (II) is particularly meaningful as it highlights how the operators $M_k$ naturally arise from an asymptotic expansion with an identically-zero remainder when one attempts to find out independent expressions for each of its coefficients: see \S3-C and \S3-D. Moreover this algorithm provides an asymptotic interpretation of formulas (2.31) related to representation (2.29). 
 
 Let us consider a generic element $u\in$ span $(\phi_1,\dots,\phi_n)$ of the type
$$
 u=a_1\phi_1(x)+\dots+a_n\phi_n(x),\ a_i\neq 0\ \ \forall \ i, \leqno (8.1)$$
  which we interpret as an asymptotic expansion at $x_0$ (with a zero remainder). 
  
  We shall first present the algorithm for a  $C.F.$ of type (II) as it is more simple to describe.

     \vspace{10pt}
{\bf Proposition 8.1}\ (The algorithm for a special $C.F.$ of type (II)).
  \textit{Let  $(\phi_1,\dots,\phi_n)$ satisfy conditions {\rm (2.25),(2.23)} and {\rm(2.24)} with all the Wronskians strictly positive; then the following algorithm yields the special  global  
$C.F.$ of  $L_{\phi_1,\dots,\phi_n}$ of type {\rm(II)} at $x_0$ in $(2.39)$ together with  $(n-1)$ asymptotic expansions which, after dividing by the first meaningful term on the right, concide with the expansions obtained by applying to $(8.1)$ the operators $M_k$ defined in $(3.27)$. Formulas for the coefficients $a_k$ in $(3.41)$, with $\epsilon_{h-1}=1$, are reobtained.}
     
     (A)\  \textit{\underline{Verbal description of the algorithm.}}
  
   $1^{st} \ step.$ \ \  Divide both sides of (8.1) by the first term on the right, which is the term with the largest growth-order at $x_0$, and then take derivatives so obtaining 
$$
  \left(  \frac{u(x)}{\phi_1(x)}  \right )'=a_2 \left(  \frac{\phi_2(x)}{\phi_1(x)} \right )'+\dots+
  a_n \left(  \frac{\phi_n(x)}{\phi_1(x)} \right )' . \leqno (8.2)$$
 
 Notice that division of both sides by the first term on the right just yields the expansion obtained by applying to $(8.1)$ the operator $M_1$; a similar remark applies to each of the subsequent expansions both in this and in the next  proposition.
 
 $2^{nd}\  step.$ \ \ Divide both sides of (8.2) by the first term on the right and take derivatives so obtaining
 $$
  \left[ \frac{1}{(\phi_2(x)/\phi_1(x))' } \left( \frac{u(x)}{\phi_1(x)} \right)' \right]'=
 a_3\left( \frac{(\phi_3(x)/\phi_1(x))'}{(\phi_2(x)/\phi_1(x))' } \right)'+\dots 
 \leqno (8.3) $$
 $$ 
 \dots+a_n\left( \frac{(\phi_n(x)/\phi_1(x))'}{(\phi_2(x)/\phi_1(x))' } \right)'.
 $$
 
$3^{ rd}\ step.$ \ \  Repeat the procedure on (8.3) dividing by the first term on the right and then taking derivatives so getting
    $$
   \left[ \frac{1}{ \left( \frac{(\phi_3/ \phi_1)'}{(\phi_2/\phi_1)' } \right)' }
   \left(  \frac{1}{(\phi_2/\phi_1)' } \left( \frac{u}{\phi_1}\right)' \right)' \right]'= \leqno(8.4)$$
   $$
  =a_4\left( \frac{(\phi_4/ \phi_1)'}{(\phi_2/\phi_1)' } \right)' \left /
   \left( \frac{(\phi_3/ \phi_1)'}{(\phi_2/\phi_1)' } \right)' +\dots +
   a_n\left( \frac{(\phi_n/ \phi_1)'}{(\phi_2/\phi_1)' } \right)' \right /
   \left( \frac{(\phi_3/ \phi_1)'}{(\phi_2/\phi_1)' } \right)' . $$
   
   \vskip5pt
   \textit{Iterating the procedure each  of the obtained relation is an identity on $[T,x_0[$ and is an asymptotic expansion at $x_0$, hence at each step we are dividing by the term on the right with the largest growth-order at $x_0$.  Notice that at each step the asymptotic expansion loses its first meaningful term and this is the same  phenomenon occurring in differentiation of Taylor's formula.} \textit{After $n$ steps we arrive at an identity}:
 $$
  [q_{n-1}(\dots(q_0u)'\dots)']'\equiv 0 \ \textit{on} \ [T,x_0[\ , \leqno (8.5) $$
 \textit{where the $q_i$'s coincide with those in $(2.35)$.}

 (B)\  \underline{\textit{Schematic description of the algorithm.}}
   $$ \begin{aligned}
    &  u\ =\ a_1\underbrace{\phi_1}+\dots+a_n\phi_n \\
    &  \underset{\textnormal{ d \   \& \  d  }} {  \uparrow_{\longleftarrow\longleftarrow\longleftarrow} \downarrow }
     \end{aligned} \leqno   \textit{Step\ \ "1":} $$
     $$ \begin{aligned}
  	&  (u/\phi_1)'=a_2\underbrace{(\phi_2/\phi_1)'}+\dots +a_n(\phi_n/\phi_1)' \\
	& \underset{\textnormal{ d \   \& \  d  }} { \ \   \uparrow_{\longleftarrow \longleftarrow \longleftarrow \longleftarrow \leftarrow} \downarrow }
	 \end{aligned} \leqno   \textit{Step\ \ "2":} $$
 $$  \begin{aligned}
 & \left( \frac{(u/ \phi_1)'}{(\phi_2/\phi_1)' } \right)' =
a_3 \underbrace{\left( \frac{(\phi_3/ \phi_1)'}{(\phi_2/\phi_1)' }\right)'} +\dots +
a_n \left ( \frac{(\phi_n/ \phi_1)'}{(\phi_2/\phi_1)' }\right)' \\
& \underset{\textnormal{ d \   \& \  d  }} {\ \ \ \ \ \  \uparrow_{\longleftarrow \longleftarrow\longleftarrow\longleftarrow\longleftarrow\longleftarrow\leftarrow} \downarrow}
 \end{aligned}  \leqno   \textit{Step\ \ "3":} $$
   $$\begin{aligned}
       &  \left[ \frac{1}{ \left( \frac{(\phi_3/ \phi_1)'}{(\phi_2/\phi_1)' } \right)' }
   \left(  \frac{1}{(\phi_2/\phi_1)' } \left( \frac{u}{\phi_1}\right)' \right)' \right]'= \\
   & \underset{\textnormal{ d \   \& \  d  }} {  \hspace {75pt} \uparrow_{\longleftarrow \longleftarrow\longleftarrow \longleftarrow\longleftarrow\longleftarrow \longleftarrow\longleftarrow}  } 
\end{aligned}  \leqno \textit{Step "4":}$$
   
   $$\begin{aligned}
   &=a_4\underbrace{\left( \frac{(\phi_4/ \phi_1)'}{(\phi_2/\phi_1)' } \right)'   \left/ 
   \left( \frac{(\phi_3/ \phi_1)'}{(\phi_2/\phi_1)' } \right)' \right.}+\dots +
   a_n\left( \frac{(\phi_n/ \phi_1)'}{(\phi_2/\phi_1)' } \right)'  \left /
  \left( \frac{(\phi_3/ \phi_1)'}{(\phi_2/\phi_1)' } \right)' \right., \\
&  \longleftarrow \longleftarrow\longleftarrow \longleftarrow\longleftarrow \downarrow
  \end{aligned} $$
 \
\newline
 \textit{and so on, where the symbol ``d $\&$ d" stands for the two operations  ``\underline{divide}'' both sides by the underbraced term on the right and then ``\underline{differentiate}'' both sides} (the equation in each step being the result of the preceding step).

  \vspace{10pt}
  {\bf Proposition} 8.2\ (The algorithm for the $C.F.$ of type (I)).
  \textit{Let  $(\phi_1,\dots,\phi_n)$ satisfy conditions {\rm (2.23), (2.24), (2.25)}; then the following algorithm yields ``the" global  
$C.F.$ of  $L_{\phi_1,\dots,\phi_n}$ of type {\rm(I)} at $x_0$ and  $(n-1)$ asymptotic expansions which, after dividing by the last meaningful term on the right,  coincide $($apart from the signs of the coefficients$)$ with the expansions obtained by applying to $(8.1)$ the operators $L_k $ defined in $(3.24)$.}
     
  (A)\  \textit{\underline{Verbal description of the algorithm.}}
  
   $1^{st} \ step.$ \ \  Divide both sides of (8.1) by the last term on the right, which is the term with the smallest growth-order at $x_0$, and then take derivatives so obtaining 
$$
  \left(  \frac{u(x)}{\phi_n(x)}  \right )'=a_1 \left(  \frac{\phi_1(x)}{\phi_n(x)} \right )'+\dots+
  a_{n-1} \left(  \frac{\phi_{n-1}(x)}{\phi_n(x)} \right )' . \leqno (8.6)$$
 
 $2^{nd}\  step.$ \ \ Divide both sides of (8.6) by the last term on the right and then take derivatives so obtaining
 $$
  \left[ \frac{1}{(\phi_{n-1}(x)/\phi_n(x))' } \left( \frac{u(x)}{\phi_n(x)} \right)' \right]'=
 a_1\left( \frac{(\phi_1(x)/\phi_n(x))'}{(\phi_{n-1}(x)/\phi_n(x))' } \right)'+\dots 
 \leqno (8.7) 
 $$
 $$ 
 \dots+a_{n-2}\left( \frac{(\phi_{n-2}(x)/\phi_n(x))'}{(\phi_{n-1}(x)/\phi_n(x))' } \right)'.
 $$
 
$3^{ rd}\ step.$ \ \  Repeat the procedure on (8.7) dividing by the last term on the right and then taking derivatives so getting
    $$
   \left[ \frac{1}{ \left( \frac{(\phi_{n-2}/ \phi_n)'}{(\phi_{n-1}/\phi_n)' } \right)' }
   \left(  \frac{1}{(\phi_{n-1}/\phi_n)' } \left( \frac{u}{\phi_n}\right)' \right)' \right]'= \leqno(8.8)$$
   $$
  =a_1\left( \frac{(\phi_1/ \phi_n)'}{(\phi_{n-1}/\phi_n)' } \right)' \left /
   \left( \frac{(\phi_{n-2}/ \phi_n)'}{(\phi_{n-1}/\phi_n)' } \right)' +\dots +
   a_{n-3}\left( \frac{(\phi_{n-3}/ \phi_n)'}{(\phi_{n-1}/\phi_n)' } \right)' \right /
   \left( \frac{(\phi_{n-2}/ \phi_n)'}{(\phi_{n-1}/\phi_n)' } \right)' . $$
   
   \vskip5pt
   \textit{Iterating the procedure each  of the obtained relation is an identity on $[T,x_0[$ and is an asymptotic expansion at $x_0$, hence at each step we are dividing by the term on the right with the smallest growth-order at $x_0$. Also notice that at each step the asymptotic expansion loses its last meaningful term and this is a phenomenon different from that occurring in differentiation of Taylor's formula} (see the foregoing proposition). \textit{After $n$ steps we arrive at an identity}:
 $$
  [p_{n-1}(\dots(p_0u)'\dots)']'\equiv 0 \ \textit{on} \ [T,x_0[\ , \leqno (8.9) $$
 \textit{where the $p_i$'s coincide, signs apart, with those in $(2.43)$.} 
  
  (B)\  \underline{\textit{Schematic description of the algorithm.}}
    $$   \begin{aligned}
   &  u=a_1\phi_1+\dots+a_n\underbrace{\phi_n}  \\
    &   \underset{\textnormal{ d \   \& \  d  }} { \uparrow_{\longleftarrow \longleftarrow\longleftarrow \longleftarrow\longleftarrow\longleftarrow \leftarrow} \downarrow }
    \end{aligned} \leqno   \textit{Step\ \ "1":} $$
$$ \begin{aligned}
  (u/ & \phi_n)'=a_1(\phi_1/\phi_n)'+\dots +a_{n-1}\underbrace{(\phi_{n-1}/\phi_n)'} \\
& \underset{\textnormal{ d \   \& \  d  }} {  \uparrow_{\longleftarrow \longleftarrow\longleftarrow \longleftarrow\longleftarrow\longleftarrow \longleftarrow\longleftarrow\longleftarrow\longleftarrow\longleftarrow\longleftarrow} \downarrow } 
 \end{aligned}  \leqno   \textit{Step\ \ "2":} $$
$$  \begin{aligned}
&  \left( \frac{(u/ \phi_n)'}{(\phi_{n-1}/\phi_n)' } \right)' =
a_1 \left( \frac{(\phi_1/ \phi_n)'}{(\phi_{n-1}/\phi_n)' }\right)' +\dots +
a_{n-2}\underbrace{ \left( \frac{(\phi_{n-2}/ \phi_n)'}{(\phi_{n-1}/\phi_n)' }\right)'} \\
& \underset{\textnormal{ d \   \& \  d  }} {  \hspace {30pt} \uparrow_{\longleftarrow \longleftarrow\longleftarrow \longleftarrow\longleftarrow\longleftarrow \longleftarrow\longleftarrow\longleftarrow\longleftarrow\longleftarrow\longleftarrow\longleftarrow
\longleftarrow\longleftarrow\longleftarrow\longleftarrow} \downarrow  }  \end{aligned}  \leqno   \textit{Step\ \ "3":} $$
     $$\begin{aligned}
       &  \left[ \frac{1}{ \left( \frac{(\phi_{n-2}/ \phi_n)'}{(\phi_{n-1}/\phi_n)' } \right)' }
   \left(  \frac{1}{(\phi_{n-1}/\phi_n)' } \left( \frac{u}{\phi_n}\right)' \right)' \right]'= \\
   & \underset{\textnormal{ d \   \& \  d  }} {  \hspace {75pt} \uparrow_{\longleftarrow \longleftarrow\longleftarrow \longleftarrow\longleftarrow\longleftarrow \longleftarrow\longleftarrow\longleftarrow\longleftarrow}  } 
\end{aligned}  \leqno \textit{Step "4":}$$
   
   $$\begin{aligned}
   &=a_1\left( \frac{(\phi_1/ \phi_n)'}{(\phi_{n-1}/\phi_n)' } \right)'   \left/ 
   \left( \frac{(\phi_{n-2}/ \phi_n)'}{(\phi_{n-1}/\phi_n)' } \right)' \right.+\dots +
   a_{n-3}\underbrace{\left( \frac{(\phi_{n-3}/ \phi_n)'}{(\phi_{n-1}/\phi_n)' } \right)'  \left /
  \left( \frac{(\phi_{n-2}/ \phi_n)'}{(\phi_{n-1}/\phi_n)' } \right)' \right.}, \\
&  \longleftarrow \longleftarrow\longleftarrow \longleftarrow\longleftarrow\longleftarrow \longleftarrow\longleftarrow\longleftarrow\longleftarrow\longleftarrow\longleftarrow\longleftarrow\longleftarrow\longleftarrow\longleftarrow\longleftarrow\longleftarrow\downarrow
  \end{aligned} $$
 \newline
 \textit{and so on  with the symbol ``d $\&$ d" reminding of the two operations ``\underline{divide}'' both sides by the underbraced term on the right and then ``\underline{differentiate}'' both sides} ( the equation in each step being the result of the preceding step).  

  \vspace{10pt} 
  {\it Remarks}.\ 1. In order to obtain any $C.F.$ by the above procedures  one may simply choose $a_1=\dots=a_n=1$. 
  
  2. If some operator $L_n$ of type $(2.1)_{1,2}$ is known to be disconjugate on a left neighborhood of $x_0, I_{x_0}$, and if 
 $$
 \begin{cases}
      & \phi_1,\dots ,\phi_n \in \text{ker} \ L_n, \\
      &  \phi_1(x)\gg \dots \gg \phi_n(x),\  x\to x_0^-, \end{cases} \leqno (8.10)$$
 then the algorithm in Proposition 8.2  yields the $C.F.$ of type (I) at $x_0$, valid on the whole given interval $I_{x_0}$ by Proposition 2.2-(ii). But the algorithm in Proposition 8.1 yields a $C.F.$ of type (II) at $x_0$ valid on some left neighborhood of $x_0$, the largest of them being characterized by the nonvanishingness of all the Wronskians  $W( \phi_1, \dots ,\phi_i), 1\leq i\leq n-1,$ which does not automatically follow from the nonvanishingness of $W( \phi_n,\phi_{n-1}, \dots ,\phi_i), 1\leq i\leq n,$: see remarks preceding Proposition 2.3.
 
 3. By applying the above algorithms to (8.1) with a random choice of the term to be factored out at each step one may well obtain, after $n$ steps, a factorization valid on a certain subinterval of the given interval but, in general for $n\geq 3$, it will not be a $C.F.$ at an endpoint. 
  
4. In practical applications of the algorithms there is a fatal pitfall to avoid, namely the temptation at each step of suppressing brackets, cancelling possible opposite terms and rearranging in an aestetically-nicier asymptotic scale. This in general gives rise to a factorization of an operator quite different from $L_{\phi_1,\dots,\phi_n}$. Hence it is essential that all the terms coming from a single term in the preceding step  be kept grouped together as a single term to the end of the procedure: see examples at the end  of this section.
 
\vskip5pt
    {\it Proof of Proposition} 8.1, that of  Proposition 8.2 being exactly the same after replacing $(\phi_1,\dots,\phi_n)$ by ($\phi_n,\dots,\phi_1)$. We have to prove that the $q_i$'s in (8.5) coincide with those in (2.35) and this does not seem to be an obvious fact though it is made explicit in the algorithm that the first three coefficients $q_0,q_1,q_2$ coincide with P\'olya's coefficients in (2.35). Now, known $q_i$, our algorithm constructs $q_{i+1}$, for $i\geq 2$, by the following rule:
    $$
    1/q_{i+1}=\big[q_i\times(\textit{the expression for $1/q_i$ with the one change:}
   \  \phi_{i+1}\  \textit{replaced by}\  \phi_{i+2})\big]', \leqno (8.11) $$ 
    hence it is enough to show that P\'olya's expression for $1/q_{i+1}$ is obtained by the same rule.  We present two different proofs, the first being based on the equivalent  representations (2.35) and (2.37). We have:
    $$
     \big[q_i\times(\textit{expression in $(2.29)$ for $1/q_i$ with}\ \phi_{i+1}\  \textit{replaced by}\  \phi_{i+2})\big]'= \leqno (8.12) $$
         $$
     =\left[\frac{\left[W(\phi_1,\dots ,\phi_{i})\right]^2}
     {W(\phi_1,\dots ,\phi_{i-1})W(\phi_1,\dots ,\phi_i,\phi_{i+1})}\cdot
     \frac{W(\phi_1,\dots ,\phi_{i-1})W(\phi_1,\dots ,\phi_i,\phi_{i+2})}{\left[W(\phi_1,\dots ,\phi_{i})\right]^2} \right]'=$$
     $$
     =\left[\frac{W(\phi_1,\dots ,\phi_i,\phi_{i+2})}{W(\phi_1,\dots ,\phi_i,\phi_{i+1})}\right]'  \overset{\ by (2.31)\ }{=} \frac{1}{q_{i+1}}\ . $$ 
     
     It is also clear that the various identities obtained are nothing but those obtained by applying to (8.1) the operators $M_k$ defined in (3.26) which, by (3.29), differ from  $L_{\phi_1,\dots ,\phi_k}$ by a factor which is a non-vanishing function. 
     
     The second proof is based on a nontrivial identity involving Wronskians of Wronskians, Karlin [7$_{bis}$; p. 60], which we report here in the version needed in our proof:
     $$
    W(g_1,\ldots ,g_n,f_1,f_2)\cdot  W(g_1,\ldots ,g_n)\equiv W\big(W(g_1,\ldots ,g_n,f_1),W(g_1,\ldots ,g_n,f_2)\big) .\leqno(8.13)$$
    
    Comparing the expressions in (2.37) and those given by our algorithm we see that the two procedures coincide if the following identity holds true:
    $$
    \left[{W(\phi_1,\dots ,\phi_i,\phi_{i+2})\over W(\phi_1,\dots ,\phi_i,\phi_{i+1})}\right]'\equiv \left\{{\big[W(\phi_1,\dots ,\phi_{i-1},\phi_{i+2})\big/ W(\phi_1,\dots ,\phi_{i-1},\phi_i)\big]'\over
    \big[W(\phi_1,\dots ,\phi_{i-1},\phi_{i+1})\big/ W(\phi_1,\dots ,\phi_{i-1},\phi_i)\big]'}\right\}'.\leqno(8.14)$$
    
    We shall show the validity of this identity even if the outer derivatives are suppressed. Using the elementary formula $(g_2/g_1)'=W(g_1,g_2)\cdot g_1^{-2}$ we have:
    $$\begin{cases}\displaystyle
    {\big[W(\phi_1,\dots ,\phi_{i-1},\phi_{i+2})\big/ W(\phi_1,\dots ,\phi_{i-1},\phi_i)\big]'\over
    \big[W(\phi_1,\dots ,\phi_{i-1},\phi_{i+1})\big/ W(\phi_1,\dots ,\phi_{i-1},\phi_i)\big]'}=\\ \\ \displaystyle
    ={ W\big(W(\phi_1,\ldots ,\phi_i),W(\phi_1,\ldots ,\phi_{i-1},\phi_{i+2})\big)\ 
    \big(W(\phi_1,\ldots ,\phi_i)\big)^{-2}\over
    W\big(W(\phi_1,\ldots ,\phi_i),W(\phi_1,\ldots ,\phi_{i-1},\phi_{i+1})\big)\ 
    \big(W(\phi_1,\ldots ,\phi_i)\big)^{-2}}=\cdots\end{cases}\leqno(8.15)$$
     by (8.13) with $g_1,\ldots ,g_n$ replaced by $\phi_1,\ldots ,\phi_{i-1}$
     $$
     \cdots={W(\phi_1,\dots  ,\phi_{i-1},\phi_i,\phi_{i+2})\cdot W(\phi_1,\dots ,\phi_{i-1})\over W(\phi_1,\dots  ,\phi_{i-1},\phi_i,\phi_{i+1})\cdot W(\phi_1,\dots ,\phi_{i-1})}=
     {W(\phi_1,\dots ,\phi_i,\phi_{i+2})\over W(\phi_1,\dots ,\phi_i,\phi_{i+1})}\ .$$
     
     \hfill $\Box$
     
  \vspace{15pt}
   \textit{An example illustrating the two algorithms.}  Consider the fourth-order operator $L$ of type $ (2.1)_{1,2}$ such that
    $$
    \text{ker}\ L=\text{span}\ (e^x,x,\log x,1),  \leqno (8.16) $$
    acting on $AC^3]0,+\infty)$ or even on $C^\infty]0.+\infty)$. Starting from the asymptotic scale
    $$
    e^x\gg x\gg \log x\gg 1 ,\ x\to +\infty, \leqno (8.17) $$
    the algorithm in Proposition 8.2 yields in sequence:
     $$ \underset{\uparrow_{ \leftarrow\longleftarrow\longleftarrow\longleftarrow\longleftarrow\longleftarrow\longleftarrow\longleftarrow}\downarrow} {u=e^x+x+\log x+1}; $$
    $$\underset{\uparrow_{ \longleftarrow\longleftarrow\longleftarrow\longleftarrow\longleftarrow}\downarrow}
    { u'=e^x+1+x}^{-1}; $$
    $$ \underset{\uparrow_{ \leftarrow\longleftarrow\longleftarrow\longleftarrow\longleftarrow\longleftarrow\longleftarrow\longleftarrow}\downarrow}{ (xu')'=\left[(x+1)e^x\right]+1;}$$
    $$ \underset{\uparrow_{ \longleftarrow\longleftarrow\longleftarrow\longleftarrow\longleftarrow}\downarrow}{(xu')''=[(x+2)e^x] };$$
    $$ \left[(x+2)^{-1}e^{-x}(xu')''\right]'\equiv 0.$$
    
    Hence
    $$Lu\equiv x^{-1}(x+2)e^x\left[(x+2)^{-1}e^{-x}(xu')''\right]', \leqno (8.18)$$
    where
    $$ p_1(x)=x;\ p_2(x)=1;\ p_3(x)=(x+2)^{-1}e^{-x};\leqno (8.19)$$
    and (8.18) is ``the " global $C.F.$ of $L$ of type (I) at  $+\infty$.
    
    On the other hand the algorithm in Proposition 8.1 yields in sequence:
    $$ \underset {\uparrow_{\longleftarrow\leftarrow}\downarrow}{u=e^x}
    +x+\log x+1;$$
    $$ 
    (e^{-x}u)'=\underbrace{[(1-x)e^{-x}]}+[(x^{-1}-\log x)e^{-x}]-e^{-x};$$
    
    $$ [(1-x)^{-1}e^x (e^{-x}u)']'= \underbrace {[(1-x)^{-2}(-\log x+1+x^{-1}-x^{-2})]} -(1-x)^{-2};$$
    $$\left[(1-x)^2{(-\log x+1+x^{-1}-x^{-2})}^{-1}[(1-x)^{-1}e^x (e^{-x}u)']'\right]'=$$ 
    $$=\underbrace{ {(-\log x+1+x^{-1}-x^{-2})}^{-2}x^{-3}(-x^2-x+2)}; $$ \par
    $$ \left[{(-\log x+1+x^{-1}-x^{-2})}^2 x^3{(-x^2-x+2)}^{-1} \right. \times $$
    $$\times \left. \left[(1-x)^2{(-\log x+1+x^{-1}-x^{-2})}^{-1}[(1-x)^{-1}e^x (e^{-x}u)']'\right]'\right]'\equiv 0.$$
   (The underbraced  terms on the right are those by which one must divide and then differentiate.) Hence:
    $$Lu\equiv (1-x)^{-1}x^{-3}(-x^2-x+2){(-\log x+1+x^{-1}-x^2)}^{-1}\times  \leqno (8.20)$$
    $$\times\big[ {(-\log x+1+x^{-1}-x^2)}^2 x^3{(-x^2-x+2)}^{-1} \times $$
    $$\times \left[(1-x)^2(-\log x+1+x^{-1}-x^{-2})^{-1}[(1-x)^{-1}e^x (e^{-x}u)']'\right]' \big]', $$
    where
    $$\begin{cases}
      & \overline{q}_0(x):=e^{-x}; \ \overline{q}_1(x):=(1-x)^{-1} e^x;\\
      & \overline{q}_2(x):=(1-x^2){(-\log x+1+x^{-1}-x^{-2})}^{-1} \sim  x^2/\log x,\ x\to +\infty; \\
      & \overline{q}_3(x):=x^3(-x^2-x+2)^{-1}{(-\log x+1+x^{-1}-x^{-2})}^2 \sim  \\ 
      & \sim -x(\log x)^2,\ x\to +\infty;\ \ \ \int^{+\infty}1/|\overline{q}_i|<+\infty,\ i=1,2,3.
\end{cases} \leqno (8.21) $$

Hence (8.20) is a $C.F.$ of $L$ of type (II) at $+\infty$ valid on the largest neighborhood of $+\infty$ whereon 
$$
1-x\neq0;\ 1-x^2\neq 0;\ (\log x-1-x^{-1}+x^{-2})\neq 0;\ x^2+x-2\neq 0, $$
which is easily seen to be the interval $]1,+\infty)$. In conclusion: changing the signs of the $\overline{q}_i$'s, if necessary, we get a P\'olya-Mammana factorization of $L$ on
$]1,+\infty)$ which is a $C.F.$ of type (II) at $+\infty$.
The standard non-factorized form of $L$ is 
$$
Lu\equiv u^{(4)}+x^{-1}(6-x^2)(x+2)^{-1}u^{(3)}-2x^{-1}(x+3)(x+2)^{-1}u''. \leqno (8.22) $$

In the various steps of the above procedures one must carefully avoid the temptation of rearranging the terms in the right-hand side in (supposedly) nicier asymptotic scales. For instance the first procedure involves quite simple terms and only at the last-but-one step we may split the remaining  term on the right by writing
$$
\underset{\uparrow_{ \longleftarrow\longleftarrow\longleftarrow\longleftarrow\longleftarrow}\downarrow}
{(xu')''=xe^x+2e^x} $$
and taking $e^x$ as the term  with the smallest growth-order. The procedure then yields
$$
(e^{-x}(xu')'')'=1\ \ \text{and}\ \ (e^{-x}(xu')'')''\equiv 0. $$

This gives a fifth-order operator
$$
\tilde{L}u:=x^{-1}e^x(e^{-x}(xu')'')'', $$
distinct from the given fourth-order operator.

On the contrary the second procedure offers a great number of temptations! For instance if one rewrites the result of the first step as
$$
(e^{-x}u)'=-xe^{-x}- \log x\cdot e^{-x}+x^{-1}e^{-x}, \leqno (8.23) $$
and then goes on applying the second algorithm to (8.23) as if the right-hand side would be an asymptotic expansion with three meaningful terms, one gets:
$$
(x^{-1}e^x(e^{-x}u)')'=\begin{cases}
      & \underbrace{(x^{-2}\log x-x^{-2})}-2x^{-3}, \\
      & \underbrace{x^{-2}\log x}-x^{-2}-2x^{-3},
\end{cases} \leqno (8.24) $$
 the only difference between the two expressions on the right being the term-grouping. 
From the upper relation in (8.24), considered as an asymptotic expansion at $+\infty$ with two meaningful terms, one gets:
$$
\left[ x^2(\log x-1)^{-1}(x^{-1}e^x(e^{-x}u)')'\right ]'=(-2x^{-1}(\log x-1)^{-1})'= $$
$$
=2x^{-2}(\log x-1)^{-1}+2x^{-2}(\log x-1)^{-2}=2x^{-2}(\log x-1)^{-2}\log x, $$
and then 
$$ 
\left \{x^2(\log x)^{-1}(\log x-1)^2\left [x^2(\log x-1)^{-1}(x^{-1}e^x(e^{-x}u)')'\right ]' \right\}'\equiv 0,$$
whose left-hand side is a fourth-order operator distinct from our operator.  

If, instead, one starts from the lower relation in (8.24), considered as an asymptotic expansion at $+\infty$ with three meaningful terms, one gets:
$$
\left (x^2(\log x)^{-1}(x^{-1}e^x(e^{-x}u)')'\right)'=x^{-1}(\log x)^{-2}+2x^{-2}(\log x)^{-1}
+2x^{-2}(\log x)^{-2} $$
and so forth in an endless process leading nowhere!!
    
  \vspace{10pt}
   \textit{An example showing that application of the procedure in the algorithms regardless of the relative growth-orders of the terms may yield a non-canonical factorization.} Let us consider $Lu:=u'''$ acting on $AC^2]0,+\infty)$ or even on $C^{\infty}]0,+\infty)$ and the tern $(1,x,x^2)$ which satisfies
    $$
    1\gg x\gg x^2,\ x\to 0^+;\ \ \ x^2\gg x\gg 1, x\to +\infty, $$
    and is such that all the possible Wronskians constructed with these three functions do not vanish on $]0,+\infty)$. We now apply the following two procedures:
    $$
      \underset{\uparrow_{ \longleftarrow\longleftarrow\longleftarrow\leftarrow}\downarrow}{u=x^2+x}+1$$
      $$ \swarrow \searrow $$
      $$\begin{tabular}{c|c}
 $ \underset{\uparrow_{ \longleftarrow\longleftarrow\longleftarrow}\downarrow}
 {(x^{-1}u)'=1-}x^{-2}$  \ \ \ \ \ \  & 
  $ \underset{\uparrow_{ \longleftarrow\longleftarrow\longleftarrow\longleftarrow}\downarrow}
  {(x^{-1}u)'=1-x^{-2}}$    \\ & \\
  $ (x^{-1}u)''=2x^{-3} $ \ \ \ \ \ \  &  $(x^2 (x^{-1}u)')'=2x$  \\ & \\
   $(x^3(x^{-1}u)'')'\equiv 0$ \ \ \ \ \ \  & $ \left( x^{-1} (x^2 (x^{-1}u)')'\right)' \equiv 0$
\end{tabular} $$   \newline
and obtain the two factorizations:
$$
u'''\equiv x^{-2}(x^3(x^{-1}u)'')';\ \ \ \ u'''\equiv (x^{-1}(x^2(x^{-1}u)')')' \leqno(8.25)$$
both valid on $]0,+\infty)$ but none of which is a $C.F.$ at any of the endpoints.
    
\vspace{20pt}
 \centerline{\bf References}

\begin{description}
 \item[\ {[1]}]
    {\small W. A. C}{\scriptsize OPPEL,} {\small \emph{Disconjugacy.} Lecture Notes in Mathematics, vol. 220. 
    Springer-Verlag, Berlin, 1971.}
    \item[\ {[2]}]
    {\small A. G}{\scriptsize RANATA,} {\small Canonical factorizations of disconjugate
    differential operators, \emph{SIAM J. Math. Anal.}, \textbf{11}(1980), 160-172.}
    \item[\ {[3]}]
    {\small A. G}{\scriptsize RANATA,} {\small Canonical factorizations of disconjugate
    differential operators -Part II, \emph{SIAM J. Math. Anal.}, \textbf{19}(1988),
    1162-1173.}
    \item[\ {[4]}]
    {\small A. G}{\scriptsize RANATA,} {\small Polynomial asymptotic expansions in the real domain: the geometric, the factorizational, and the stabilization approaches,
    \emph{Analysis Mathematica}, \textbf{33}(2007), 161-198.} DOI: 10.1007/s10476-007-0301-0.
  \item[\ {[5]}]
    {\small A. G}{\scriptsize RANATA,} {\small The problem of differentiating an asymptotic expansion in real powers. Part I: Unsatisfactory or partial results by classical approaches,
    \emph{Analysis Mathematica}, \textbf{36}(2010), 85-112.} DOI: 10.1007/s10476-010-0201-6.
  \item[\ {[6]}]
    {\small A. G}{\scriptsize RANATA,} {\small The problem of differentiating an asymptotic expansion in real powers. Part II: factorizational theory,
    \emph{Analysis Mathematica}., \textbf{36}(2010), 173-218.} DOI: 10.1007/s10476-010-0301-3.
 \item[\ {[7]}]
    {\small A. G}{\scriptsize RANATA,} {\small Analytic theory of finite asymptotic expansions in the real domain. Part I: two-term expansions of differentiable functions, \emph{Analysis Mathematica}, \textbf{37}(2011), 245-287.} DOI: 10.1007/s10476-011-0402-7. For an enlarged version with corrected misprints see  arXiv:1405.6745v1 [math.CA] 26 May 2014.
    \item[{[7]$_{bis}$}]
       {\small S. K}{\scriptsize ARLIN,} {\small \emph{Total positivity}, Vol I. Stanford University Press, Stanford, California, 1968.}
    \item[\ {[8]}]
       {\small S. K}{\scriptsize ARLIN} {\small and W. S}{\scriptsize TUDDEN,}
    {\small \emph{Tchebycheff systems: with
    applications in analysis and statistics}. Interscience, New York, 1966.}
    \item[\ {[9]}]
    {\small A. Y}{\scriptsize u.} {\small L}{\scriptsize EVIN,} {\small Non-oscillation of solutions of the equation
    $x^{(n)}+p_1(t)x^{(n-1)}+\dots+p_n(t)x=0$, \emph{Uspehi Mat. Nauk}, \textbf{24}
    (1969), no. 2 (146), 43-96; \emph{Russian Math. Surveys}, \textbf{24}(1969), no.
    2, 43-99.}
     \item[{[10]}]
     {\small M.-L. M}{\scriptsize AZURE,} {\small Quasi Extended Chebyshev spaces and weight functions, \emph{Numer. Math.}, \textbf{118}(2011), 79-108. DOI: 10.1007/s00211-010-0312-9.}
      \item[{[11]}]
        {\small A. M. O}{\scriptsize STROWSKI,} {\small Note on the Bernoulli-L'Hospital rule, \emph{Amer. Math. Monthly}, \textbf{83}(1976), 239-242.}
       \item[{[12]}]
    {\small G. P}{\scriptsize \'OLYA,} {\small On the mean-value theorem corresponding to a given linear homogeneous differential equation, \emph{Trans. Amer. Math. Soc.}, \textbf{24}(1922), 312-324.}
        \item[{[13]}]
    {\small I. J. S}{\scriptsize CHOENBERG,} {\small Two applications of approximate differentiation formulae: an extremum problem for multiply monotone functions and the differentiation of asymptotic expansions, \emph{J. Math. Anal. Appl.}, \textbf{89}(1982), 251-261.}
     \item[{[14]}]
    {\small W. F. T}{\scriptsize RENCH,} {\small Canonical forms and principal systems for general disconjugate equations, \emph{Trans. Amer. Math. Soc.}, \textbf{189}(1974), 139-327.}
    \end{description}

\end{document}